\documentclass[11pt]{article}
\usepackage{pdfsync}

\usepackage{array} 
\usepackage{amsfonts}

\usepackage{amssymb}
\usepackage{graphicx}
\usepackage{color}
\usepackage{hyperref}
\usepackage{amsmath}
\newtheorem{thm}{Theorem}[section]
\newtheorem{cor}[thm]{Corollary}
\newtheorem{lem}[thm]{Lemma}
\newtheorem{prop}[thm]{Proposition}
\newtheorem{conj}[thm]{Conjecture}

\numberwithin{equation}{section}

\newcommand{\RR}{\mathbb{R}}
\newcommand{\NN}{\mathbb{N}}

\def\qed{\,\unskip\kern 6pt \penalty 500d
\raise -2pt\hbox{\vrule \vbox to8pt{\hrule width 6pt
\vfill\hrule}\vrule}\par}

\begin{document}
\title{\textbf{The  Fisher-KPP problem \\[1mm]
with doubly nonlinear diffusion}\\[5mm]}

\author{{\Large Alessandro Audrito\footnote{Also affiliated with Universit\`a degli Studi di Torino, Italy.} ~and~ Juan Luis V\'azquez}\\ [4pt]
Departamento de Matem\'{a}ticas \\ [4pt] Universidad
Aut\'{o}noma de Madrid
} 
\date{\vspace{-5ex} }

\maketitle

\begin{abstract}
The famous Fisher-KPP reaction-diffusion model combines linear diffusion with the typical KPP reaction term, and appears in a number of relevant applications in biology and chemistry. It is remarkable as a mathematical model since it possesses a family of travelling waves that describe the asymptotic behaviour of a large class solutions $0\le u(x,t)\le 1$  of the problem posed in the real line. The existence of propagation waves with finite speed has been confirmed in some related models and disproved in others.
We investigate here the corresponding theory when the linear diffusion is replaced by the ``slow'' doubly nonlinear diffusion and we find travelling waves that represent the wave propagation of more general solutions even when we extend the study to several space dimensions. A similar study is performed in the critical case that we call ``pseudo-linear'', i.e., when the operator is still nonlinear but has homogeneity one. With respect to the classical model and the ``pseudo-linear'' case, the ``slow'' travelling waves exhibit free boundaries.
\end{abstract}
%
%

%
%
\section{Introduction}\label{SECTIONINTRODUCTION}
In this paper we study the doubly nonlinear reaction-diffusion problem posed in the whole Euclidean space
\begin{equation}\label{eq:REACTIONDIFFUSIONEQUATIONPLAPLACIAN}
\begin{cases}
\begin{aligned}
\partial_tu = \Delta_p u^m + f(u) \quad &\text{in } \RR^N\times(0,\infty) \\
u(x,0) = u_0(x) \quad\qquad &\text{in } \RR^N.
\end{aligned}
\end{cases}
\end{equation}
We first discuss the problem of the existence of travelling wave solutions and, later, we use that information to establish the asymptotic behaviour for large times of the solution $u = u(x,t)$ with general initial data and for different ranges of the parameters $m > 0$ and $p > 1$. We recall that the $p$-Laplacian is a nonlinear operator defined for all $1 \leq p < \infty$ by the formula
\[
\Delta_p v := \nabla\cdot(|\nabla v|^{p-2}\nabla v)
\]
and we consider the more general diffusion term \ $
\Delta_p u^m := \Delta_p(u^m) = \nabla\cdot(|\nabla (u^m)|^{p-2}\nabla (u^m))$,
that we call ``doubly nonlinear''operator since it presents a double power-like nonlinearity. Here,  $\nabla$ is the spatial gradient while $\nabla\cdot$ is the spatial divergence. The doubly nonlinear operator (which can be thought as the composition of the $m$-th power and the $p$-Laplacian) is much used in the elliptic and parabolic literature (see for example \cite{C-D-D-S-V:art, DB:book, EstVaz:art, Kal:survey, Leib:art, Lindq:art, V1:book,V2:book} and their references) and allows to recover the Porous Medium operator choosing $p = 2$ or the $p$-Laplacian operator choosing $m = 1$. Of course, choosing $m=1$ and $p = 2$ we obtain the classical Laplacian.

\begin{figure}[h!]
\centering
  \includegraphics[scale=0.4]{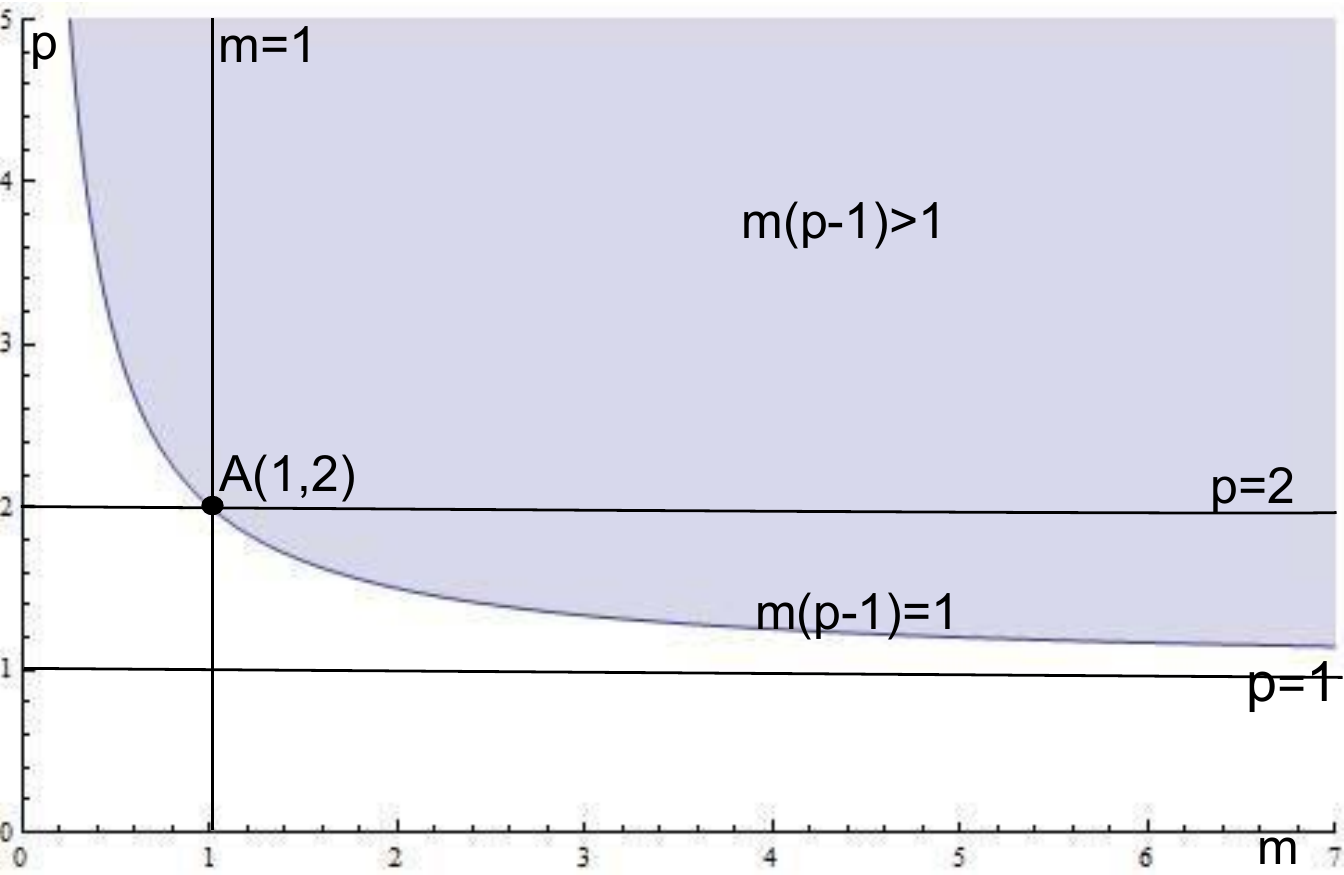}
  \caption{The ``slow diffusion'' area and the ``pseudo-linear'' line in $(m,p)$-plane.}\label{fig:PARAMETERSMP}
\end{figure}

In order to fix the notations and avoid cumbersome expressions in the rest of the paper, we define the constants
\[
\mu := (m-1)(p-1) \qquad \text{and} \qquad \gamma := m(p-1)-1,
\]
which will play an important role in our study. For example, note that the doubly nonlinear operator can be written in the form $\Delta_p u^m = m^{p-1}\nabla\cdot(u^{\mu}|\nabla u|^{p-2}\nabla u)$. The importance of the constant $\gamma$ is related to the properties of the fundamental solutions of the ``pure diffusive'' doubly nonlinear parabolic equation and we refer the reader to \cite{V1:book} and Section \ref{PRELIMINARIESINTRO}. From the beginning, we consider parameters $m > 0$ and $p > 1$ such that
\[
\gamma \geq 0.
\]
This is an essential restriction. We refer to the assumption $\gamma > 0$ (i.e. $m(p-1)>1$) as the ``slow diffusion'' assumption, while ``pseudo-linear'' assumption when we consider $\gamma = 0$ (i.e. $m(p-1)=1$). In Figure \ref{fig:PARAMETERSMP} the corresponding ranges in the $(m,p)$-plane are reported.

The function $f(\cdot)$ is a reaction term modeled on the famous references by Fisher \cite{Fisher:art}, and Kolmogorov-Petrovski-Piscounoff \cite{K-P-P:art} in their seminal works on the existence of travelling wave propagation. The classical example is the logistic term $f(u) = u(1-u)$, $0 \leq u \leq 1$. More generally, we will assume that
\begin{equation}\label{eq:ASSUMPTIONSONTHEREACTIONTERM}
\begin{cases}
f : [0,1] \to \RR \text{ and } f \in C^1([0,1]) \\
f(0) = 0 = f(1) \; \text{ and } \; f(u) > 0 \text{ in } (0,1) \\
f \text{ is concave in } [0,1].
\end{cases}
\end{equation}
see  \cite{Aro-Wein1:art, Aro-Wein2:art, Fife:book, Fisher:art, K-P-P:art, Murray:book} and Section \ref{SUBSECTIONMODELINTRO} for a more complete description of the model. Note that the positivity and the concavity of $f(\cdot)$ guarantee that $f'(0)>0$ while $f'(1)<0$. Moreover, typical assumptions on the initial datum are
\begin{equation}\label{eq:ASSUMPTIONSONTHEINITIALDATUM}
\begin{cases}
u_0 : \RR^N \to \RR \text{ is continuous with compact support: } u_0 \in \mathcal{C}_c(\RR^N) \\
u_0 \not \equiv 0 \; \text{ and } \; 0 \leq u_0 \leq 1.
\end{cases}
\end{equation}
We anticipate that our main result, Theorem \ref{NTHEOREMCONVERGENCEINNEROUTERSETS}, can be proved even if we relax assumptions \eqref{eq:ASSUMPTIONSONTHEINITIALDATUM} on the initial data. In particular, we will show in Corollary \ref{CORASYMBEHAMOREGENDATUM} that the statement of Theorem \ref{NTHEOREMCONVERGENCEINNEROUTERSETS} holds even when we take initial data with a certain spatial exponential decay for $|x| \sim \infty$. We ask the reader to keep in mind that the assumptions on $f(\cdot)$ and $u_0(\cdot)$ combined with the Maximum Principle allow us to deduce that all solutions of problem \eqref{eq:REACTIONDIFFUSIONEQUATIONPLAPLACIAN} are trapped between the steady states $u = 0$ and $u = 1$, i.e., $0 \leq u(x,t) \leq 1$ in $\RR^N\times[0,\infty)$.

We point out that the proof the asymptotic result follows a very detailed analysis of the existence of admissible travelling waves which is done in Section \ref{SECTIONEXISTENCEOFTWS} and Section \ref{CLASSIFICATIONEXISTENCETW}. Actually, we are going to extend the well-known results on the Cauchy problem \eqref{eq:REACTIONDIFFUSIONEQUATIONPLAPLACIAN} with $m=1$ and $p=2$ (see Subsection \ref{SECTIONRESULTSFORLINEARDIFFUSION}) to the case $m > 0$ and $p > 1$ (such that $\gamma > 0$), trying to stress the main differences between the linear and the nonlinear case. In particular, it will turn out that the nonlinearity of our diffusion operator causes the appearance of nonnegative solutions with {\sl free boundaries} which are not admitted in the case $m = 1$ and $p = 2$. Moreover, in Section \ref{CLASSIFICATIONEXISTENCETW}, we generalize the linear theory to the ranges of parameters $m > 0$ and $p > 1$ such that $\gamma = 0$, which can be seen as a limit case of the choice $\gamma > 0$. Finally, we point out that the methods used in the PDEs part (to study the asymptotic behaviour of the solution of problem \eqref{eq:REACTIONDIFFUSIONEQUATIONPLAPLACIAN}) are original and do not rely on the proofs given for the linear case, as \cite{Aro-Wein1:art, Aro-Wein2:art,K-P-P:art}.

%
\subsection{The Fisher-KPP model}\label{SUBSECTIONMODELINTRO}
The reaction term $f(u) = u(1-u)$ was introduced in the same year, by Fisher (\cite{Fisher:art}) in the context of the dynamics of populations, and by Kolmogorov, Petrovsky and Piscounoff (\cite{K-P-P:art}) who studied a similar model with a more general function $f(\cdot)$. It was followed by a wide literature, using both probabilistic and analytic methods (see  \cite{Aro-Wein1:art, Aro-Wein2:art, Fife:book, McKean:art, Murray:book,Uchiyama:art}).

Let us recall the model presented in \cite{Fisher:art}. Fisher considered a density of population distributed uniformly on the real line and supposed that it was divided in two sub-densities or sub-groups of individuals $\rho_1 = \rho_1(x,t)$ and $\rho_2 = \rho_2(x,t)$ such that $\rho_1 + \rho_2 = 1$ (i.e., no other sub-groups were admitted). The idea was to modeling the evolution in time of the entire population when genetic mutations occurred in some individuals. So, in order to fix the ideas, it is possible to suppose that $\rho_1$ represents the part of the population with mutated genes, while $\rho_2$ be the old generation. According to the evolution theories, it was supposed that mutations improved the capacity of survival of individuals and generated a consequent competition between the new and the old generation. Then, a ``rate'' of gene-modification was introduced, i.e., a coefficient $\kappa > 0$ describing the ``intensity of selection in favour of the mutant gene'' (\cite{Fisher:art}, p.355). In other words, $\kappa$ represents the probability that the coupling of two individuals from different groups generated a mutant. Then, the author of \cite{Fisher:art} supposed that the individuals $\rho_1$ sprawled on the real line with law $-D\partial_x\rho_1,$ where $D > 0$. This formula represents the typical flow of diffusive type caused by the movement of the individuals. Hence, according to the previous observations Fisher obtained the equation
\[
\partial_t\rho_1 = D\partial_{xx}\rho_1 + \kappa\rho_1\rho_2 \quad \text{in } \RR\times(0,\infty).
\]
Recalling that $\rho_1 + \rho_2 = 1$ and renaming $\rho_1 = u$, he deduced the evolution equation of the mutated individuals
\begin{equation}\label{eq:ONEDIMENSIONALKPPEQUATIONINTRO}
\partial_tu = D\partial_{xx} u + \kappa u(1 - u) \quad \text{in } \RR\times(0,\infty),
\end{equation}
which is the one-dimensional version of \eqref{eq:REACTIONDIFFUSIONEQUATIONPLAPLACIAN} with $m = 1$ and $p = 2$ and it explains the choice of the reaction term $f(u) = u(1-u)$. Note that we can suppose $D = \kappa = 1$, after performing a simple variable scaling.

Equation \eqref{eq:ONEDIMENSIONALKPPEQUATIONINTRO} and its generalizations have been largely studied in the last century. The first results were showed in \cite{Fisher:art} and \cite{K-P-P:art}, where the authors understood that the group of mutants prevailed on the other individuals, tending to occupy all the available space, and causing the extinction of the rest of the population. A crucial mathematical result was discovered: the speed of proliferation of the mutants is approximately constant for large times and the general solutions of \eqref{eq:ONEDIMENSIONALKPPEQUATIONINTRO} could be described by introducing particular solutions called Travelling Waves (TWs).
%
%
%
%
\subsection{Travelling waves}
They are special solutions with remarkable applications, and there is a huge mathematical literature devoted to them. Let us review the main concepts and definitions.

\noindent Fix $m > 0$ and $p \geq 1$ and assume that we are in space dimension 1 (note that when $N = 1$, the doubly nonlinear operator has the simpler expression $\Delta_pu^m = \partial_x(|\partial_xu^m|^{p-2}\partial_xu^m)$). A TW solution of the equation
\begin{equation}\label{eq:PLAPLACIANREACTIONDIFFUSIONEQUATIONINTRO}
\partial_tu = \Delta_pu^m + f(u) \quad \text{in }\RR\times(0,\infty)
\end{equation}
is a solution of the form $u(x,t) = \varphi(\xi)$, where $\xi = x + ct$, $c > 0$ and the \emph{profile} $\varphi(\cdot)$ is a real function. In our application to the Fisher-KPP problem, we will need the profile to satisfy
\begin{equation}\label{eq:CONDITIONONPHIADMISSIBLETWINTRO}
0 \leq \varphi \leq 1, \quad \varphi(-\infty) = 0, \; \varphi(\infty) = 1 \quad \text{and} \quad \varphi' \geq 0.
\end{equation}
In that case say that $u(x,t) = \varphi(\xi)$ is an \emph{admissible} TW solution. Similarly, one can consider admissible TWs of the form $u(x,t) = \varphi(\xi)$ with $\xi = x - ct$, $\varphi$ decreasing and such that $\varphi(-\infty) = 1$ and $\varphi(\infty) = 0$.  But it is easy to see that  these two options are equivalent, since the the shape of the profile of the second one can be obtained by reflection of the first one, $\varphi_{-c}(\xi)=\varphi_c(-\xi)$, and it moves in the opposite direction of propagation. In the rest of the paper, we will use both types of TWs, but if it is not specified, the reader is supposed to consider TW solutions defined as in \eqref{eq:CONDITIONONPHIADMISSIBLETWINTRO}. We recall that, since admissible TWs join the levels $u = 0$ and $u = 1$, they are called ``change of phase type'' and are really important in many physical applications.

\begin{figure}[h!]
  \centering
  \includegraphics[scale =0.4]{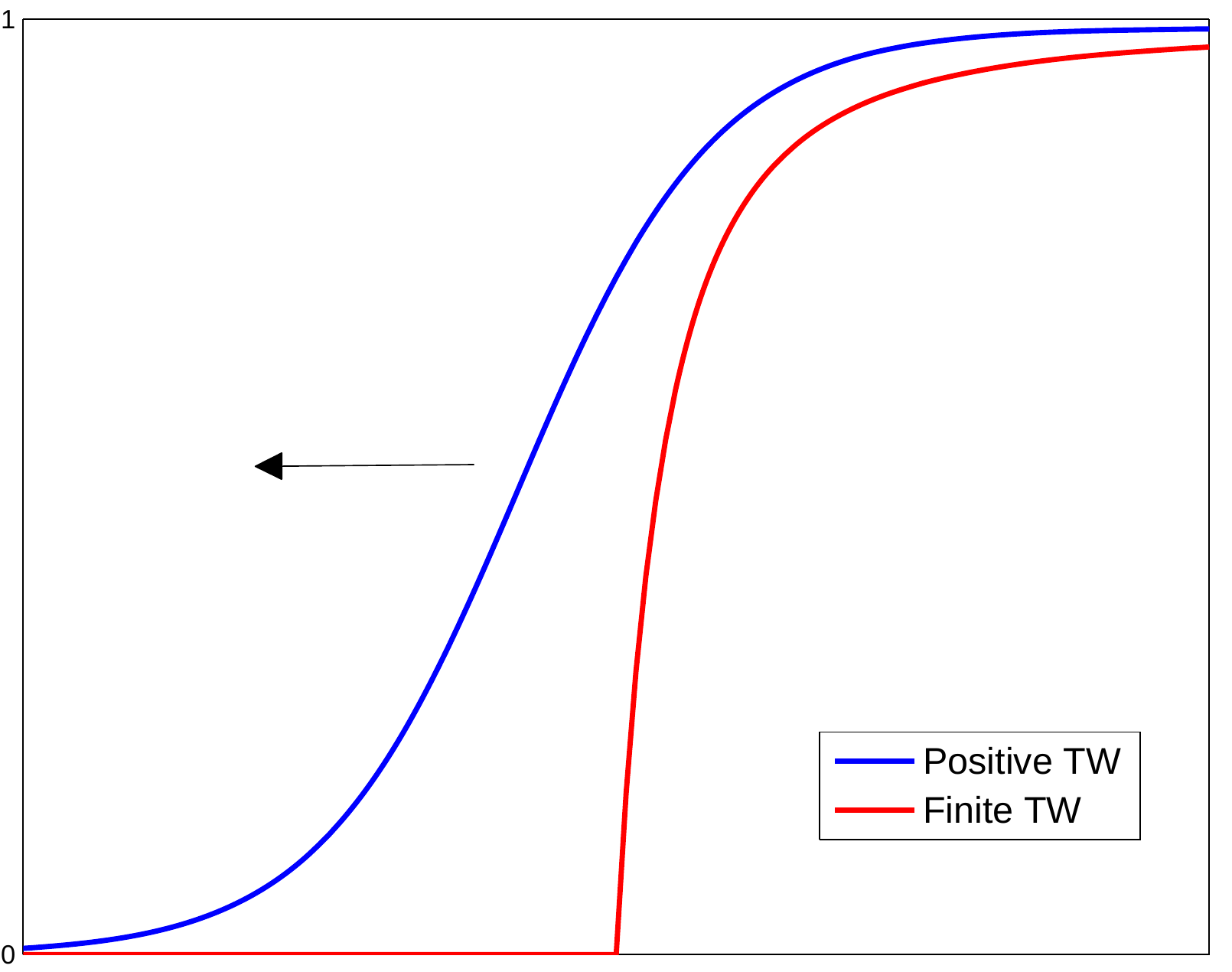}
  \includegraphics[scale =0.4]{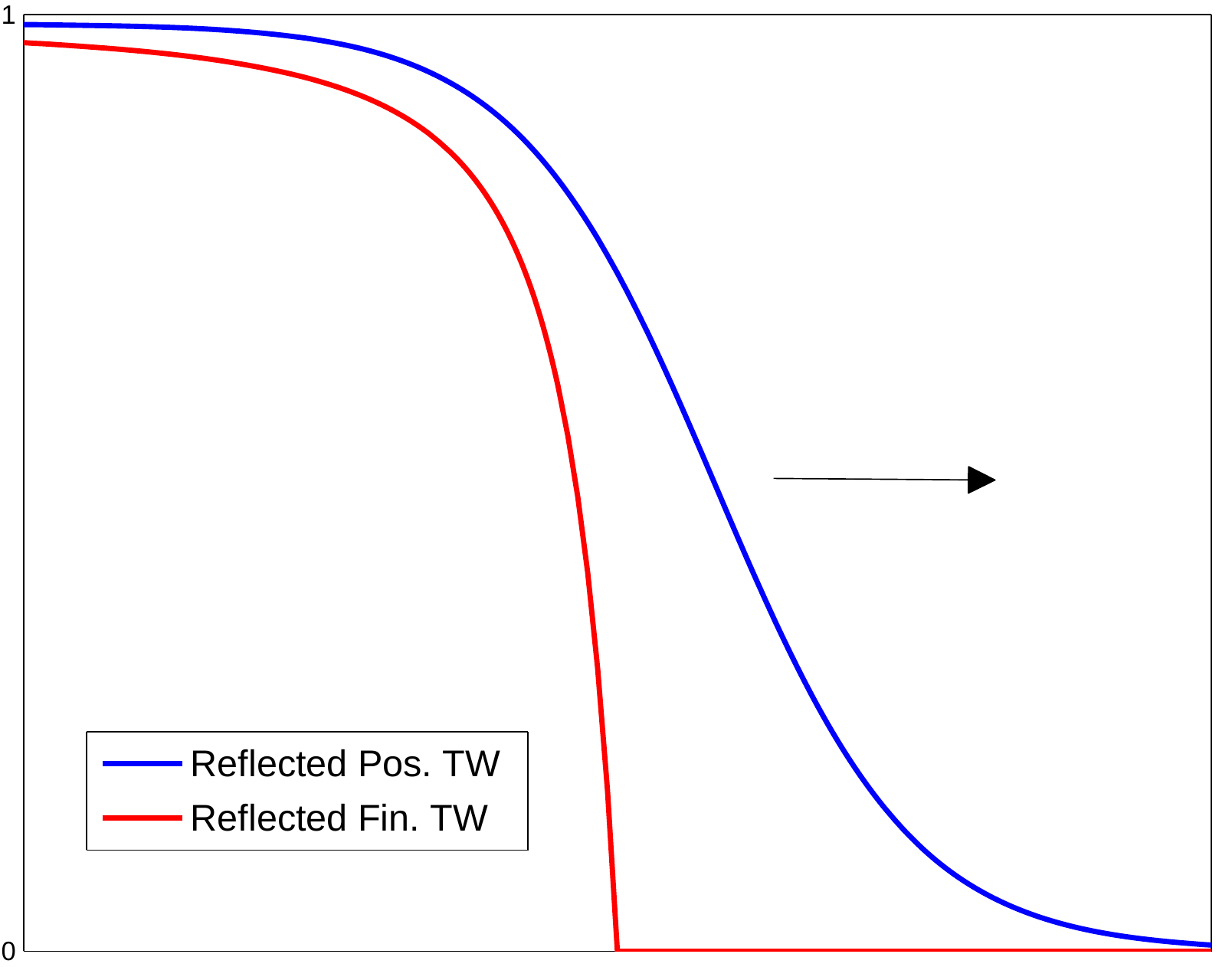}
  \caption{Examples of admissible TWs and their ``reflections'': Finite and Positive types}\label{fig:ADMISSIBLEANDREFLECTEDFINITETW}
\end{figure}

\noindent Finally, an admissible TW is said \emph{finite} if $\varphi(\xi) = 0$ for $\xi \leq \xi_0$ and/or $\varphi(\xi) = 1$ for $\xi \geq \xi_1$, or \emph{positive} if $\varphi(\xi) > 0$, for all $\xi \in \RR$. The line $x = \xi_0 - ct$ that separates the regions of positivity and vanishing of $u(x,t)$ is then called the \emph{free boundary}. Same name would be given to the line $x = \xi_1 - ct$ and $\varphi(\xi) = 1$ for $\xi \geq \xi_1$ with  $x_1$ finite, but this last situation will not happen.

\noindent Though the TWs are essentially one-dimensional objects, there is a natural extension to several dimensions. Indeed, we can consider equation \eqref{eq:PLAPLACIANREACTIONDIFFUSIONEQUATIONINTRO} in dimension $N \geq 1$ and look for solution in the form $u(x,t) = \varphi(x\cdot n + ct)$ where $x \in \RR^N$ and $n$ is a unit vector of $\RR^N$. Note that the direction of the vector $n$ coincides with the direction of propagation of the wave and we are taken back to the study of \eqref{eq:PLAPLACIANREACTIONDIFFUSIONEQUATIONINTRO} in spatial dimension 1.
%

\subsection{Results for linear diffusion}\label{SECTIONRESULTSFORLINEARDIFFUSION}
As we have said, the case of linear diffusion has been extensively studied. In \cite{Aro-Wein1:art, Aro-Wein2:art,K-P-P:art} it has been proved that the existence of admissible TWs strongly depends on the speed of propagation $c > 0$ of the wave and in particular, on the local behaviour of the solution near the steady states $u = 0$ and $u = 1$. In particular, they showed the following theorem, assuming the parameters $D = \kappa = 1$.
\begin{thm}\label{THMEXISTENCETWLINEARINTRO}
There exists a critical speed of propagation $c_{\ast} = 2$ such that for all $0 < c < c_{\ast}$, there are not admissible TW solutions for equation \eqref{eq:ONEDIMENSIONALKPPEQUATIONINTRO}, while for all $c \geq c_{\ast}$, there exists a unique positive TW solution (up to its horizontal reflection or horizontal displacement).
\end{thm}
In this case there are no free boundaries. Theorem \ref{THMEXISTENCETWLINEARINTRO} and suitable comparison principles allowed the authors of \cite{Aro-Wein1:art} and \cite{Aro-Wein2:art} to understand the asymptotic behaviour of the general solution of \eqref{eq:ONEDIMENSIONALKPPEQUATIONINTRO} stated in the following theorem.
\begin{thm}\label{THMASYMPTOTICSTLARGELINEARCASEINTRO}
The solution $u = u(x,t)$ of equation \eqref{eq:ONEDIMENSIONALKPPEQUATIONINTRO} with initial datum of type \eqref{eq:ASSUMPTIONSONTHEINITIALDATUM} satisfies
\[
\lim_{t\to\infty} u(\xi - ct, t) =
\begin{cases}
\begin{aligned}
1 \quad &\text{if }\quad |c| < c_{\ast} \\
0 \quad &\text{if }\quad |c| > c_{\ast},
\end{aligned}
\end{cases}
\]
for all $\xi \in \RR$. Equivalently, we can reformulate the previous limit saying that
\[
\begin{aligned}
&u(x,t) \to 1 \quad \text{uniformly in } \{|x| \leq ct \} \; \text{ as } t \to \infty, \quad \text{for all }\; 0 < c < c_{\ast}, \\
&u(x,t) \to 0 \quad \text{uniformly in } \{|x| \geq ct \} \; \text{ as } t \to \infty, \quad \text{for all }\; c > c_{\ast}.
\end{aligned}
\]
\end{thm}
We recall that in \cite{Aro-Wein1:art} and \cite{Aro-Wein2:art}, the authors worked with a more general reaction term $f(\cdot)$ with properties similar to \eqref{eq:ASSUMPTIONSONTHEREACTIONTERM}. It then turns out that the critical speed of propagation depends on the reaction: $c_{\ast} = 2\sqrt{f'(0)}$ (which evidently generalizes the case $f(u) = u(1-u)$). Finally, in \cite{Aro-Wein2:art}, they extended the previous theorem when the spatial dimension is greater than one.

The assertions of Theorem \ref{THMASYMPTOTICSTLARGELINEARCASEINTRO} mean that the steady state $u = 1$ is asymptotically stable while the null solution $u = 0$ is unstable and, furthermore, the asymptotic stability/instability can be measured in terms of speed of convergence of the solution which, in this case, is asymptotically linear in distance of the front location as function of time. From the point of view of the applications, it is possible to state that the density $u(x,t)$ tends to saturate the space and its resources, while the rate of propagation in space is linear with constant $c_{\ast}$. However, it is not clear what are the properties of the solution $u = u(x,t)$ on the moving coordinate $x = \xi - c_{\ast}t$ when we consider the critical speed and still it seems to represent a very challenging problem.  We will give some information in the last section for the reader's convenience.
%
%
%
%
%
\subsection{Nonlinear diffusion}\label{SUBSECTIONNONLINEARDIFFUSION}
Alongside the just described studies, researchers have paid attention to \emph{nonlinear} variants of the Fisher-KPP model which are interesting for the applications. The mathematical issue is whether the TW paradigm to explain the long-time propagation of the solution is valid or not.
\paragraph{Finite propagation. The Porous Medium case.} As we have just mentioned, an interesting class of such new models concerns equations with nonlinear diffusion. For instance, in \cite{Aronson1:art} and \cite{Aronson2:art}, Aronson considered a modified version of equation \eqref{eq:ONEDIMENSIONALKPPEQUATIONINTRO}, introducing a nonlinear Porous Medium diffusion:
\begin{equation}\label{eq:ONEDIMENSIONALKPPPOROUSMEDIUMEQUATIONINTRO}
\partial_tu =  \partial_{xx}(u^m) + f(u) \quad \text{in } \RR\times(0,\infty),
\end{equation}
with $f(\cdot)$ satisfying \eqref{eq:ASSUMPTIONSONTHEREACTIONTERM} and $m > 1$. In \cite{Aronson2:art}, he proved a nonlinear version of Theorem \ref{THMEXISTENCETWLINEARINTRO} which asserts: for all $m > 1$, there exists a critical speed of propagation $c_{\ast} = c_{\ast}(m) > 0$ such that there are not admissible TW solutions for equation \eqref{eq:ONEDIMENSIONALKPPPOROUSMEDIUMEQUATIONINTRO} for all $0 < c < c_{\ast}(m)$, whilst for all $c \geq c_{\ast}(m)$ there exists exactly one admissible TW solution. Moreover, the TW corresponding to the speed $c_{\ast}(m)$ is finite, while for $c > c_{\ast}(m)$ the TWs are positive.

\noindent Hence, it is clear that the nonlinear diffusion presents an interesting peculiarity: the TW corresponding to the value $c_{\ast}(m)$ is finite, which means that \emph{free boundaries} appear exactly as in the simpler Porous Medium equation (see for instance \cite{V1:book} and \cite{V2:book}). This observation suggests that the general solution of equation \eqref{eq:ONEDIMENSIONALKPPPOROUSMEDIUMEQUATIONINTRO} with initial datum \eqref{eq:ASSUMPTIONSONTHEINITIALDATUM} has a \emph{free boundary} too and this is evidently a really interesting difference respect to the linear case. Similar results were showed in a series of papers (\cite{DP-S:art, DePablo-Vazquez1:art, DePablo-Vazquez2:art}) for equation \eqref{eq:ONEDIMENSIONALKPPPOROUSMEDIUMEQUATIONINTRO} with a non-smooth reaction term $f(u) = a u^n(1-u)$ where $a>0$ and $n \in \RR$. The authors' concern was to study the the combination of a slow diffusion term given by $\partial_{xx}u^m$ ($m > 1$) and strong reaction term given by the singularity in the origin of the function $u^n$ when $n < 1$. They proved the following theorem (see Theorem 7.1 of \cite{DePablo-Vazquez1:art}).
\begin{thm}\label{THEOREMSLOWDIFFUSIONSTRONGREACTIONPOROUSMEDIUM}
Let $m > 1$, $n \in \RR$ and $q := m + n - 2$ and $a = 1/m$. Then there exist admissible TW solutions for equation \eqref{eq:ONEDIMENSIONALKPPPOROUSMEDIUMEQUATIONINTRO} with $f(u) = au^n(1-n)$ if and only if $q \geq 0$.

\noindent Moreover, for all $q \geq 0$, there exists a critical speed $c_{\ast} = c_{\ast}(q) > 0$ such that equation \eqref{eq:ONEDIMENSIONALKPPPOROUSMEDIUMEQUATIONINTRO} with $f(u) = au^n(1-n)$ possesses a unique admissible TW for all $c \geq c_{\ast}(q)$ and does not have admissible TWs for $0 < c < c_{\ast}(q)$.

\noindent Finally, the TW corresponding to the speed $c_{\ast}(q)$ is finite while the TWs with $c > c_{\ast}(q)$ are finite if and only if $0 < n < 1$.
\end{thm}
Theorem \ref{THEOREMSLOWDIFFUSIONSTRONGREACTIONPOROUSMEDIUM} is really interesting since it gives the quantitative ``combination'' between slow diffusion and strong reaction in order to have only finite TWs. We will sketch the proof of a generalization of this result at the end of the paper (see Section \ref{SECTIONSTRONGREACTION}).
%
%
%
%
%
\paragraph{The $\boldsymbol{p}$\,-Laplacian case.} Recently, the authors of \cite{Eng-Gav-San:art} and \cite{Gav-San:art} studied the existence of admissible TW solution for the $p$-Laplacian reaction-diffusion equation
\begin{equation}\label{eq:PLAPLACIANENGGAVSAN}
\partial_tu = \partial_x(|\partial_x u|^{p-2}\partial_xu) + f(u) \quad \text{in } \RR\times(0,\infty),
\end{equation}
for $p > 1$. They considered a quite general reaction term $f(\cdot)$ satisfying
\[
\begin{cases}
f(0) = 0 = f(1) \; \text{ and } \; f(u) > 0 \text{ in } (0,1) \\
\sup_{u \in (0,1)}u^{1-p'}f(u) = M_p < +\infty
\end{cases}
\]
where $p'$ is the conjugated of $p$, i.e., $1/p + 1/p' = 1$. They proved (see Proposition 2 of \cite{Eng-Gav-San:art}) the existence of a critical speed of propagation $c_{\ast} > 0$ depending only on $p$ and $f(\cdot)$ such that equation \eqref{eq:PLAPLACIANENGGAVSAN} possesses admissible TW solutions if and only if $c \geq c_{\ast}$ and, furthermore, they obtained the interesting bound $c_{\ast} \leq p^{1/p}(p'M_p)^{1/p'}$. Moreover, with the additional assumption
\[
\lim_{u \to 0^+}u^{1-p'}f(u) = \sup_{u \in (0,1)}u^{1-p'}f(u) = M_p < +\infty
\]
they proved that $c_{\ast} = p^{1/p}(p'M_p)^{1/p'}$. Finally, they focused on the study of TWs for equation \eqref{eq:PLAPLACIANENGGAVSAN} for different kinds of reaction terms (see \cite{Eng-Gav-San:art} and \cite{Gav-San:art} for a complete treatise).

The scene is set for us to investigate what happens in the presence of a doubly nonlinear diffusion. However, we have to underline we are not only interested in studying the existence of admissible speeds of propagation, but also in understanding the properties of the admissible TWs, i.e., if they are finite with free boundary or everywhere positive. Finally, our main goal is to apply the ODEs results to the study of the long-time behaviour of the PDEs problem \eqref{eq:REACTIONDIFFUSIONEQUATIONPLAPLACIAN}.
%
%
%
%
%
\section{Doubly nonlinear diffusion. Preliminaries and main results}\label{PRELIMINARIESINTRO}
Now we present some basic results concerning the Barenblatt solutions of the ``pure diffusive'' doubly nonlinear parabolic equation which are essential to develop our study in the next sections (the reference for this issue is \cite{V1:book}). Moreover, we recall some basic facts on existence, uniqueness, regularity and Maximum Principles for the solutions of problem \eqref{eq:REACTIONDIFFUSIONEQUATIONPLAPLACIAN}.
\paragraph{Barenblatt solutions.} Fix $m > 0$ and $p > 1$ such that $\gamma \geq 0$ and consider the ``pure diffusive'' doubly nonlinear problem:
\begin{equation}\label{eq:PARABOLICPLAPLACIANEQUATIONINTRO}
\begin{cases}
\begin{aligned}
\partial_tu = \Delta_p u^m \;\quad &\text{in } \RR^N\times(0,\infty) \\
u(t) \to M\delta_0 \quad\;\, &\text{in } \RR^N \text{ as } t \to 0,
\end{aligned}
\end{cases}
\end{equation}
where $M\delta_0(\cdot)$ is the Dirac's function with mass $M>0$ in the origin of $\RR^N$ and the convergence has to be intended in the sense of measures.
\paragraph{Case $\boldsymbol{\gamma > 0}$.} It has been proved (see \cite{V1:book}) that problem \eqref{eq:PARABOLICPLAPLACIANEQUATIONINTRO} admits continuous weak solutions in self-similar form $B_M(x,t) = t^{-\alpha}F_M(xt^{-\alpha/N})$, called Barenblatt solutions, where the \emph{profile} $F_M(\cdot)$ is defined by the formula:
\[
F_M(\xi) = \Big(C_M - k|\xi|^{\frac{p}{p-1}} \Big)_{+}^{\frac{p-1}{\gamma}}
\]
where
\[
\alpha = \frac{1}{\gamma + p/N}, \quad k = \frac{\gamma}{p}\Big(\frac{\alpha}{N}\Big)^{\frac{1}{p-1}}
\]
and $C_M>0$ is determined in terms of the mass choosing $M = \int_{\RR^N}B_M(x,t)dx$ (see \cite{V1:book} for a complete treatise). It will be useful to keep in mind that we have the formula
\begin{equation}\label{eq:RELATIONMASSESBARENBLATTSOLINTRO}
B_M(x,t) = MB_1(x,M^{\gamma}t)
\end{equation}
which describes the relationship between the Barenblatt solution of mass $M > 0$ and mass $M = 1$. We remind the reader that the solution has a \emph{free boundary} which separates the set in which the solution is positive from the set in which it is identically zero (``slow'' diffusion case).

\paragraph{Case $\boldsymbol{\gamma = 0}$.} Again we have Barenblatt solutions in self-similar form. The new profile can be obtained passing to the limit as $\gamma \to 0$:
\[
F_M(\xi) = C_M \exp \big(-k|\xi|^{\frac{p}{p-1}} \big),
\]
where $C_M > 0$ is a free parameter and it is determined fixing the mass, while now $k = (p-1)p^{-p/(p-1)}$. Note that, in this case the constant $\alpha = N/p$ and for the values $m=1$ and $p=2$, we have $\alpha = N/2$ and $F_M(\cdot)$ is the Gaussian profile. The main difference with the case $\gamma > 0$ is that now the Barenblatt solutions have no \emph{free boundary} but are always positive. This fact has repercussions on the shape of the TW solutions. Indeed, we will find finite TWs in the case $\gamma > 0$ whilst positive TWs in the case $\gamma = 0$.

\paragraph{Existence, Uniqueness, Regularity and Maximum Principles.} Before presenting the main results of this paper, we briefly discuss the basic properties of the solutions of problem \eqref{eq:REACTIONDIFFUSIONEQUATIONPLAPLACIAN}. Results about existence of weak solutions of the pure diffusive problem and its generalizations, can be found in the survey \cite{Kal:survey} and the large number of  references therein. The problem of uniqueness was studied later (see for instance \cite{Agueh:art, DBen-Her1:art, DBen-Her2:art, Li:art, Manf-Ves:art, Tsu:art, V2:book, Wu-Yin-Li:art}). The classical reference for the regularity of nonlinear parabolic equations is \cite{L-S-U:book}, followed by a wide literature. For the Porous Medium case ($p=2$) we refer to \cite{V1:book, V2:book}, while for the $p$-Laplacian case we suggest \cite{DB:book, Lindq:art} and the references therein. Finally, in the doubly nonlinear setting, we refer to \cite{Ivan:art, PorVes:art, Ves:art} and, for the ``pseudo-linear'' case, \cite{Kuusi-Sil-Urb:art}. The results obtained in these works showed the H\"{o}lder continuity of the solution of problem \eqref{eq:REACTIONDIFFUSIONEQUATIONPLAPLACIAN}. We mention \cite{DB:book, V2:book, Wu-Yin-Li:art, Zhao:art} for a proof of the Maximum Principle. Finally, we suggest \cite{Ag-Blan-Car:art} and \cite{S-V1:art} for more work on the ``pure diffusive'' doubly nonlinear equation and the asymptotic behaviour of its solutions.

\paragraph{Scaling Properties} In this paragraph we present a simple change of variable which allows us to work with reaction functions $f(\cdot)$ satisfying \eqref{eq:ASSUMPTIONSONTHEREACTIONTERM} and $f'(0) = 1$, without loss of generality. This reduction will result very useful to make the reading easier and the computations less technical.

\noindent Let $A > 0$. We define the function $v(y,s) = u(x,t)$, where $x = Ay$, $t = A^ps$, and $u = u(x,t)$ is a solution of the equation in \eqref{eq:REACTIONDIFFUSIONEQUATIONPLAPLACIAN}:
\[
\partial_t u = \Delta_{p,x}u^m + f(u) \quad \text{in } \RR^N\times(0,\infty),
\]
and $\Delta_{p,x}$ is the $p$-Laplacian with respect to the spacial variable $x \in \RR^N$. Thus, it is simple to see that $v = u(y,s)$ solves the same equation where we replace $\Delta_{p,x}$ with $\Delta_{p,y}$ (the $p$-Laplacian with respect to the spacial variable $y \in \RR^N$) and the function $f(\cdot)$ with the reaction  $\widetilde{f}(\cdot) = A^pf(\cdot)$. Hence, to have $\widetilde{f}'(0) = 1$, it is sufficient to choose $A = f'(0)^{-1/p}$. We then recover the properties of a solution of the equation with $f'(0) \not=1$ by means of the formula
\[
u(x,t) = v\,\big(f'(0)^{1/p}x,f'(0)t\big).
\]
We point out that our transformation changes both the space and the time variables. Consequently, we have to take into account the change of variable of the initial datum $u_0(x) = v_0(f'(0)^{1/p}x)$ when we consider the initial-value problem associated with our equation. Moreover, notice that the speed of propagation changes when we consider TW solutions $u(x,t) = \varphi(x + ct)$, $c > 0$. Indeed, if $c > 0$ is the speed of propagation of a TW of the equation in \eqref{eq:REACTIONDIFFUSIONEQUATIONPLAPLACIAN} with $f'(0) \not= 1$ and $\nu > 0$ is the same speed when $f'(0) = 1$, we have
\[
c = \frac{x}{t} = A^{1-p}\,\frac{y}{s} = f'(0)^{\frac{p-1}{p}}\nu.
\]
This last formula will be useful since it clearly shows how the propagation speed changes when we change the derivative $f'(0)$. Thus, from now on, we will try to make explicit the dependence of the propagation speeds on the quantity $f'(0)$ by using the previous formula. This fact allows us to state our results in general terms and, at the same time, in a standard and clear way.

\subsection{Main results and organization of the paper}\label{SECTIONPRINCIPALRESULTS}
The paper is divided in two main parts: the first one, is the ODE analysis of the existence of TW solutions joining the critical levels $u = 0$ and $u = 1$. These special solutions are then used in the second part to describe the asymptotic behaviour of general solutions of the PDE problem. Each of these two parts is divided into sections as follows:

\medskip

In Section \ref{SECTIONEXISTENCEOFTWS} we prove an extension of Theorem \ref{THMEXISTENCETWLINEARINTRO}, for the reaction-diffusion problem \eqref{eq:REACTIONDIFFUSIONEQUATIONPLAPLACIAN} in the case $m > 0$ and $p > 1$ such that $\gamma > 0$. We  study the equation
\begin{equation}\label{eq:REACDIFFPLAPLACIANTWS1}
\partial_tu = \partial_x(|\partial_xu^m|^{p-2}\partial_xu^m) + f(u) \quad \text{in } \RR\times(0,\infty).
\end{equation}
and show the following  result. Let's recall that $\gamma = m(p-1)-1$.

\begin{thm}\label{THEOREMEXISTENCEOFTWS}
Let $m > 0$ and $p > 1$ such that $\gamma > 0$. Then there exists a unique $c_{\ast} = c_{\ast}(m,p) > 0$ such that equation \eqref{eq:REACDIFFPLAPLACIANTWS1} possesses a unique admissible TW for all $c \geq c_{\ast}(m,p)$ and does not have admissible TWs for $0 < c < c_{\ast}(m,p)$. Uniqueness is understood as before, up to reflection or horizontal displacement.

\noindent Moreover, the TW corresponding to the value $c = c_{\ast}(m,p)$ is finite (i.e., it vanishes in an infinite half-line), while the TWs corresponding to the values $c > c_{\ast}(m,p)$ are positive everywhere.
\end{thm}
When we say unique we are in fact referring to the TWs with increasing profile, a second equivalent solution is obtained by reflection. Furthermore, even though there are no admissible TWs for values $c < c_{\ast}$, we construct solutions which are ``partially admissible'' and they will be used as sub-solutions in the PDE analysis of the last sections.

The study of TW solutions for equation \eqref{eq:REACDIFFPLAPLACIANTWS1} leads us to a second order ODEs problem (see equation \eqref{eq:PROFILEEQUATIONTWS}). The main difference respect to the linear case is that we use a nonstandard change of variables to pass from the second order ODE to a system of two first order ODEs (see formulas \eqref{eq:NONSTANDARDCHANGEOFVARIABLESTWS}). Moreover, we deduce the existence of the critical speed of propagation $c_{\ast}(m,p)$ via qualitative methods while, in the linear case (see \cite{Aronson2:art} and \cite{K-P-P:art}), the critical value of the parameter $c > 0$ was found by linearizing the ODE system around its critical points.

\medskip

In Section \ref{CLASSIFICATIONEXISTENCETW} we study the limit case $\gamma = 0$, that we call ``pseudo-linear'' (since  the homogeneity of the doubly nonlinear operator is one) and we prove the following theorem.
\begin{thm}\label{THEOREMEXISTENCEOFTWS1}
Let $m > 0$ and $p > 1$ such that $\gamma = 0$. Then there exists a unique critical speed
\begin{equation}\label{eq:DEFINITIONOFCASTPSEUDO}
c_{0\ast}(m,p) = p(m^2f'(0))^{\frac{1}{mp}}
\end{equation}
such that equation \eqref{eq:REACDIFFPLAPLACIANTWS1} possesses a unique admissible TW for all $c \geq c_{0\ast}(m,p)$ while it does not possess admissible TWs for $0 < c < c_{0\ast}(m,p)$. All the admissible TWs are positive everywhere. Uniqueness is understood as in the previous theorem.
\end{thm}
This case presents both differences and similarities respect to the case $\gamma > 0$ but we will follow the procedure used to show Theorem \ref{THEOREMEXISTENCEOFTWS} (with corresponding modifications). We ask the reader to note that the case $\gamma = 0$ covers and generalizes the linear case $m = 1$ and $p = 2$ and to observe that the most important difference between the case $\gamma > 0$ is that there are not finite TWs, i.e., when $\gamma = 0$ there are not admissible TWs with \emph{free boundary}. This turns out to be a really interesting fact, since it means that the positivity of the admissible TWs is not due to the linearity of the diffusion operator but depends on the homogeneity of the operator, i.e., the relations between the diffusive parameters. To conclude the differences between the cases $\gamma > 0$ and $\gamma = 0$, we will see that the the decay to zero of the TW with speed $c \geq c_{\ast}(m,p)$ (when $\gamma = 0$) differs from the decay of the positive TWs in the case $\gamma > 0$. On the other hand, note that formula \eqref{eq:DEFINITIONOFCASTPSEUDO} for $c_{0\ast}(m,p)$ agrees with the scaling for the critical speed in the case $\gamma > 0$:
\[
c_{\ast}(m,p) = f'(0)^{\frac{p-1}{p}}\nu_{\ast}(m,p),
\]
where $\nu_{\ast}(m,p)$ is the critical speed when $\gamma > 0$ and $f'(0) = 1$. This uses the fact that the critical speed when $\gamma = 0$ and $f'(0) = 1$ is $\nu_{0\ast} = pm^{2/(mp)}$, and $(p-1) \to 1/m$ as $\gamma \to 0$.

\noindent Finally, we prove that the critical speed $c_{\ast}(m,p)$ (when $\gamma > 0$), found in Theorem \ref{THEOREMEXISTENCEOFTWS}, converges to the critical speed $c_{0\ast}(m,p)$ (when $\gamma = 0$), found in Theorem \ref{THEOREMEXISTENCEOFTWS1}, as $\gamma \to 0$. In particular, we prove the following theorem.
\begin{thm}\label{THMCONTINUITYOFCASTINCLOSUREOFH}
Consider the region $\mathcal{H} = \{(m,p): \gamma = m(p-1)-1 > 0\}$. Then the function
\begin{equation}\label{eq:CASTEXTENDEDDEFINITION}
(m,p) \quad \to \quad
\begin{cases}
\begin{aligned}
c_{\ast}(m,p) \quad\; &\text{if } \gamma > 0 \\
c_{0\ast}(m,p) \quad  &\text{if } \gamma = 0
\end{aligned}
\end{cases}
\end{equation}
where $c_{\ast}(m,p)$ is the critical speed for the values $\gamma > 0$ (found in Theorem \ref{THEOREMEXISTENCEOFTWS}) while $c_{0\ast}(m,p)$ is the critical speed for the values $\gamma = 0$ (found in Theorem \ref{THEOREMEXISTENCEOFTWS1}), is continuous on the closure of $\mathcal{H}$, i.e., $\overline{\mathcal{H}} = \{(m,p): \gamma = m(p-1)-1 \geq 0\}$.
\end{thm}
The previous theorem is proved in two steps. We show the continuity of the function \eqref{eq:CASTEXTENDEDDEFINITION} in the region $\mathcal{H}$ in Lemma \ref{STABILITYOFTRAJECTORIESGGPOS} and Corollary \ref{COROLLARYCONTINUITYCASTINNERSET}. Then we extend the continuity to the closure $\overline{\mathcal{H}}$ in Lemma \ref{CONTINUITYOFCASTCOMPLETECLOSURELEM}.

We conclude this paragraph pointing out that we will use the separate notations $c_{\ast}(m,p)$ and $c_{0\ast}(m,p)$ for the cases $\gamma > 0$ and $\gamma = 0$ in Section \ref{SECTIONEXISTENCEOFTWS}, Section \ref{CLASSIFICATIONEXISTENCETW}, and in the statement and the proof of Corollary \ref{CORASYMBEHAMOREGENDATUM} while, later, we will unify them writing simply $c_{\ast}$ or $c_{\ast}(m,p)$ for all $(m,p) \in \overline{\mathcal{H}}$, thanks to the result stated in Theorem \ref{THMCONTINUITYOFCASTINCLOSUREOFH}.

\medskip

In Section \ref{SECTIONLINEAREXPANSIONSUPERLEVELSETS} we begin with the PDE analysis which is the main goal of this paper. As a first step, we study the Cauchy problem \eqref{eq:REACTIONDIFFUSIONEQUATIONPLAPLACIAN} with a particular choice of the initial datum
\begin{equation}\label{eq:INITIALDATUMTESTLEVELSETS1}
\widetilde{u}_0(x) :=
\begin{cases}
\begin{aligned}
\widetilde{\varepsilon} \qquad &\text{if } |x| \leq \widetilde{\varrho}_0 \\
0 \qquad &\text{if } |x| > \widetilde{\varrho}_0.
\end{aligned}
\end{cases}
\end{equation}
where $\widetilde{\varepsilon}$ and $\widetilde{\varrho}_0$ are positive real numbers. The solutions of problem \eqref{eq:REACTIONDIFFUSIONEQUATIONPLAPLACIAN}, \eqref{eq:INITIALDATUMTESTLEVELSETS1} are really useful since can be used as sub-solutions for the general problem \eqref{eq:REACTIONDIFFUSIONEQUATIONPLAPLACIAN}, \eqref{eq:ASSUMPTIONSONTHEINITIALDATUM} (see Lemma \ref{LEMMACOMPACTCONVERGENCETO1ASYMPTOTICBEHAVIOURPGREATER2}).

\noindent We first prove a technical lemma (see Lemma \ref{LEMMAEXPANDINGLEVELSETS}) which turns out to be essential in the study of the asymptotic behaviour of the solution of problem \eqref{eq:REACTIONDIFFUSIONEQUATIONPLAPLACIAN}. The main consequence of Lemma \ref{LEMMAEXPANDINGLEVELSETS} is the following proposition.
\begin{prop}\label{LEMMANONDISCRETEVERSIONLEVELSETS}
Let $m > 0$ and $p > 1$ such that $\gamma > 0$ and let $N \geq 1$. Then, for all $\widetilde{\varrho}_0 > 0$ and for all $\widetilde{\varrho}_1 \geq \widetilde{\varrho}_0$, there exists $\widetilde{\varepsilon} > 0$ and $t_0 > 0$, such that the solution $u(x,t)$ of problem \eqref{eq:REACTIONDIFFUSIONEQUATIONPLAPLACIAN} with initial datum \eqref{eq:INITIALDATUMTESTLEVELSETS1} satisfies
\[
u(x,t) \geq \widetilde{\varepsilon} \quad \text{in }  \{|x| \leq \widetilde{\varrho}_1/2\} \;\text{ for all } t \geq t_0.
\]
\end{prop}
This proposition asserts that for all initial datum small enough, the solution of problem \eqref{eq:REACTIONDIFFUSIONEQUATIONPLAPLACIAN} is strictly greater than a fixed positive constant on every compact set of $\RR^N$ for large times. This property will result essential  for the asymptotic study of general solutions in Section \ref{SECTIONASYMPTOTICBEHAVIOUR} (see Lemma \ref{LEMMACOMPACTCONVERGENCETO1ASYMPTOTICBEHAVIOURPGREATER2}).

\medskip

In Section \ref{SECTIONLINEAREXPANSIONSUPERLEVELSETSPSEUDO} we proceed with the PDE study done in the previous section taking $\gamma = 0$ (``pseudo-linear case''). In this case, we choose the initial datum
\begin{equation}\label{eq:INITIALDATUMTESTLEVELSETSPSEUDO1}
\widetilde{u}_0(x) :=
\begin{cases}
\begin{aligned}
\widetilde{\varepsilon} \;\;\qquad\qquad\qquad &\text{if } |x| \leq \widetilde{\varrho}_0 \\
a_0e^{-b_0|x|^{\frac{p}{p-1}}} \qquad &\text{if } |x| > \widetilde{\varrho}_0,
\end{aligned}
\end{cases}
\qquad a_0 := \widetilde{\varepsilon} e^{b_0\widetilde{\varrho}_0^{\,\frac{p}{p-1}}}
\end{equation}
where again $\widetilde{\varepsilon}$, $\widetilde{\varrho}_0$ and $b_0$ are positive numbers. The different choice of the datum is due to the different shape of the profile of the Barenblatt solutions in the case $\gamma = 0$ (see Section \ref{PRELIMINARIESINTRO}). In particular, the new datum has not compact support, but ``exponential'' tails.

\noindent Basically, we repeat the procedure followed in Section \ref{SECTIONLINEAREXPANSIONSUPERLEVELSETS}, but, as we will see in the proof of Lemma \ref{LEMMAEXPANDINGLEVELSETSPSEUDO}, the are several technical differences respect to the case $\gamma > 0$. We prove the following proposition.
\begin{prop}\label{LEMMANONDISCRETEVERSIONLEVELSETSPSEUDO}
Let $m > 0$ and $p > 1$ such that $\gamma = 0$ and let $N \geq 1$. Then, for all $\widetilde{\varrho}_0 > 0$ and for all $\widetilde{\varrho}_1 \geq \widetilde{\varrho}_0$, there exists $\widetilde{\varepsilon} > 0$ and $t_0 > 0$, such that the solution $u(x,t)$ of problem \eqref{eq:REACTIONDIFFUSIONEQUATIONPLAPLACIAN} with initial datum \eqref{eq:INITIALDATUMTESTLEVELSETSPSEUDO1} satisfies
\[
u(x,t) \geq \widetilde{\varepsilon} \quad \text{in }  \{|x| \leq \widetilde{\varrho}_1\} \;\text{ for all } t \geq t_0.
\]
\end{prop}
Even though Proposition \ref{LEMMANONDISCRETEVERSIONLEVELSETSPSEUDO} is a version of Proposition \ref{LEMMANONDISCRETEVERSIONLEVELSETS} when $\gamma = 0$, we decided to separate their proofs since there are significant deviances in the techniques followed to show them. The most important one respect to the case $\gamma > 0$ is that we work with \emph{positive} Barenblatt solutions and the free boundaries, which constitute the most technical ``trouble'', are not present.

\medskip

In Section \ref{SECTIONASYMPTOTICBEHAVIOUR}, \ref{SECTIONASYMPTOTICBEHAVIOURDIMENSION1}, and \ref{SECTIONASYMPTOTICBEHAVIOUR2} we address to the problem of the asymptotic behaviour of the solution of problem \eqref{eq:REACTIONDIFFUSIONEQUATIONPLAPLACIAN} for large times and we prove the following theorem.
\begin{thm}\label{NTHEOREMCONVERGENCEINNEROUTERSETS}
Let $m > 0$ and $p > 1$ such that $\gamma \geq 0$, and let $N \geq 1$.

\noindent (i) For all $0 < c < c_{\ast}(m,p)$, the solution $u(x,t)$ of the initial-value problem \eqref{eq:REACTIONDIFFUSIONEQUATIONPLAPLACIAN} with initial datum \eqref{eq:ASSUMPTIONSONTHEINITIALDATUM} satisfies
\[
u(x,t) \to 1 \text{ uniformly in } \{|x| \leq ct\} \;\text{ as } t \to \infty.
\]
\noindent (ii) Moreover, for all $c > c_{\ast}(m,p)$ it satisfies,
\[
u(x,t) \to 0 \text{ uniformly in } \{|x| \geq ct\} \;\text{ as } t \to \infty.
\]
\end{thm}
\begin{figure}[!ht]
  \centering
  \includegraphics[scale=0.4]{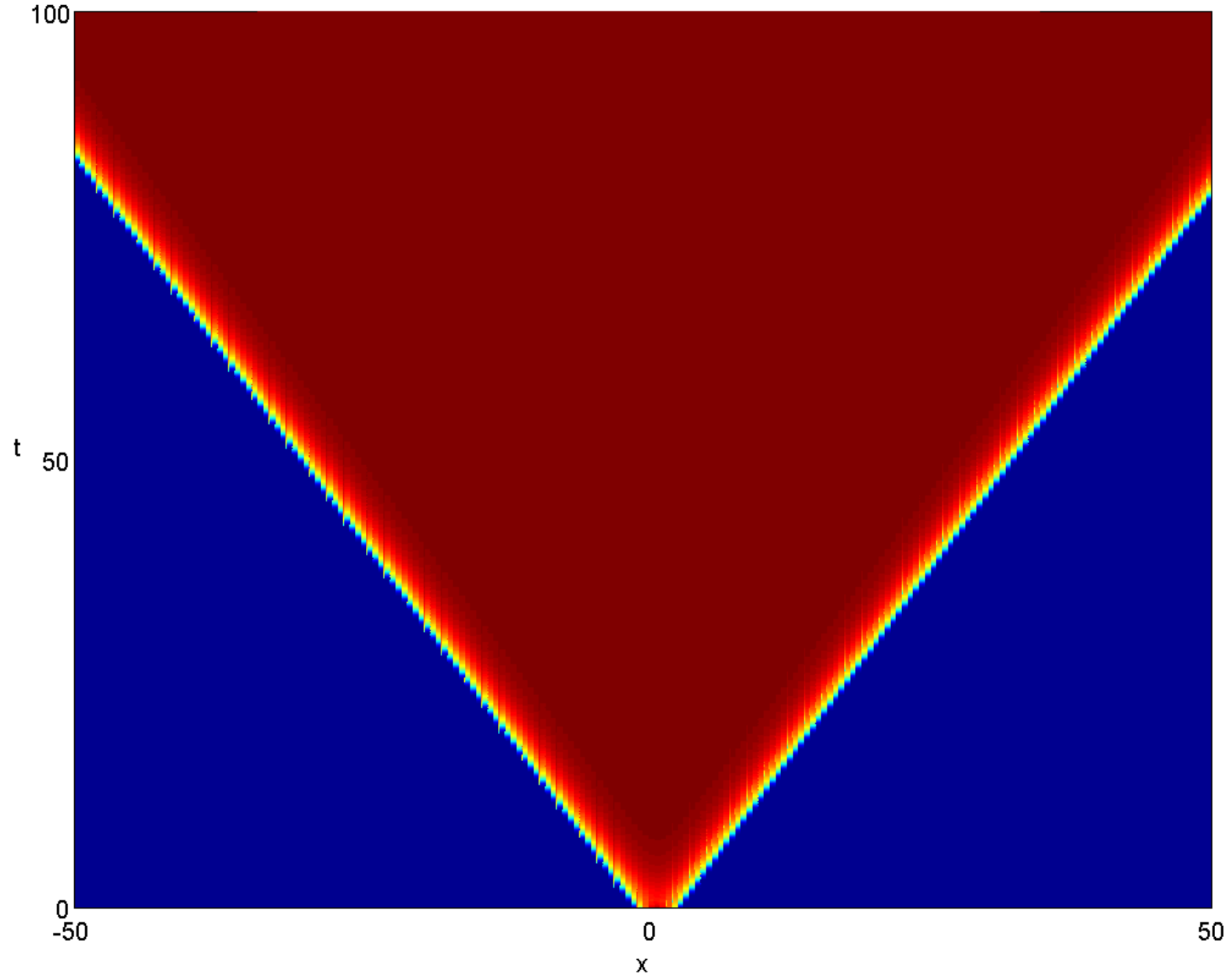} \quad
  \includegraphics[scale=0.45]{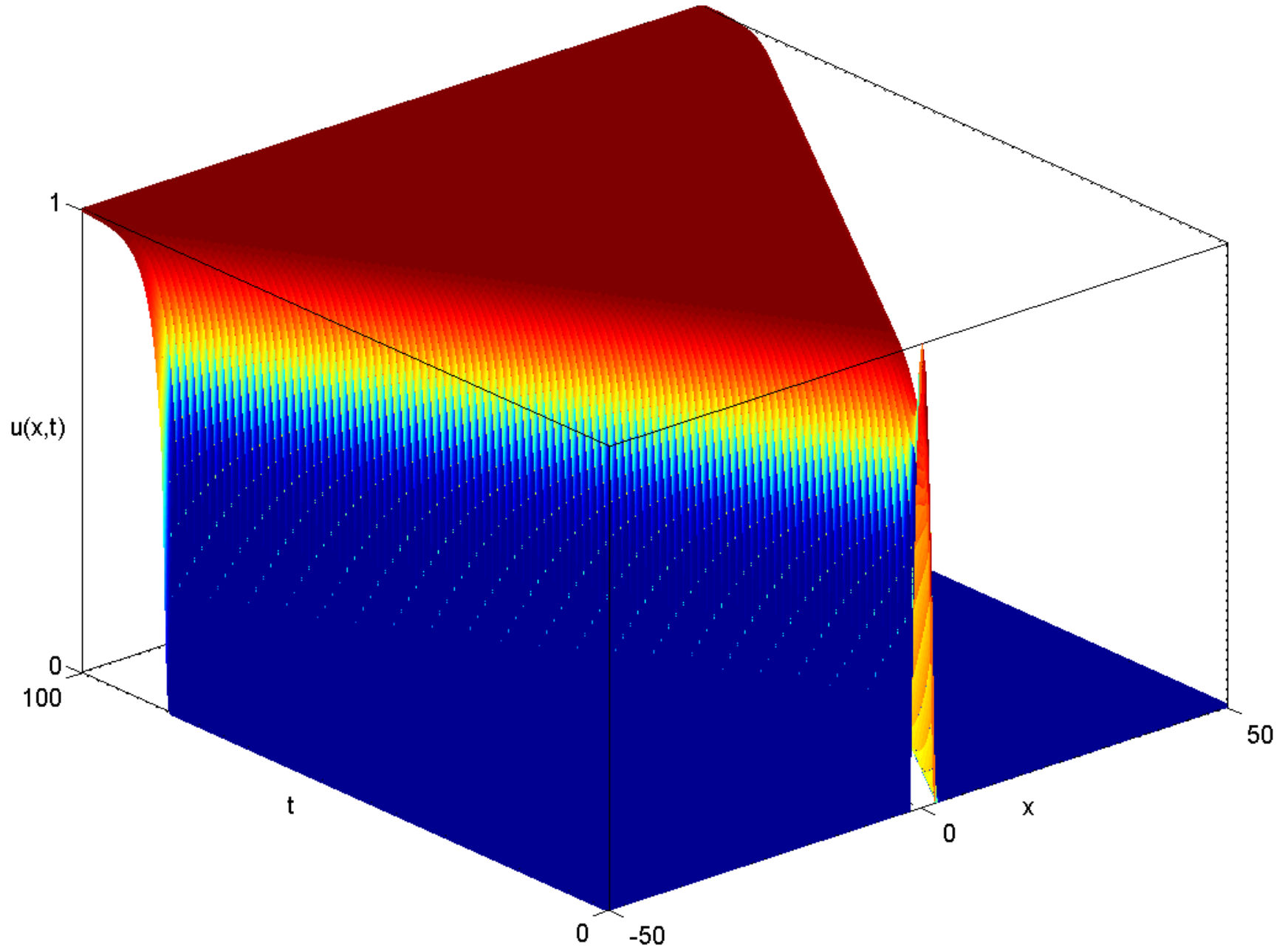}
  \caption{The slow diffusion and pseudo-linear case: convergence to 1 in the inner sets $\{|x| \leq ct\}$, for all $c < c_{\ast}$ (respectively $c < c_{0\ast}$ if $\gamma = 0$).}\label{fig:SIMULSLOWCASEEDPS}
\end{figure}
This result extends Theorem \ref{THMASYMPTOTICSTLARGELINEARCASEINTRO} to the all range of parameters $m > 0$ and $p > 1$ such that $\gamma \geq 0$ and it states that the part of the population with mutant gene invades all the available space with speed of propagation $c_{\ast}(m,p)$, driving to extinction the rest of the individuals. In the case $\gamma > 0$, we will prove that not only $u(x,t)$ converges to zero in the ``outer'' sets $\{|x| \geq ct\}$ for large values of $t > 0$ (for $c > c_{\ast}(m,p)$), but also it will turn out that $u(x,t)$ is identically zero in $\{|x| \geq ct\}$ for large times. This means that when $\gamma > 0$, the general solutions have free boundaries and it represents a significant difference respect to the ``pseudo-linear'' case ($\gamma = 0$), in which the general solutions are positive everywhere. We stress that this is due to the fact that when $\gamma > 0$ there exists a finite TW corresponding to the value $c = c_{\ast}(m,p)$.

\noindent Theorem \ref{NTHEOREMCONVERGENCEINNEROUTERSETS} is proved in three main steps. The first one consists in showing that the solution of problem \eqref{eq:REACTIONDIFFUSIONEQUATIONPLAPLACIAN} converges to 1 on every compact sets of $\RR^N$ for times large enough (see Lemma \ref{LEMMACOMPACTCONVERGENCETO1ASYMPTOTICBEHAVIOURPGREATER2}). Then, using the TW solutions found in Theorem \ref{THEOREMEXISTENCEOFTWS} as sub-solutions and super-solutions for the general solution of problem \eqref{eq:REACTIONDIFFUSIONEQUATIONPLAPLACIAN}, we prove a one dimensional version of Theorem \ref{NTHEOREMCONVERGENCEINNEROUTERSETS} (see Proposition \ref{THEOREMCONVERGENCEINNEROUTERSETS}). The last step consists in studying the radial solutions of problem \eqref{eq:REACTIONDIFFUSIONEQUATIONPLAPLACIAN} by using solutions of the one-dimensional problem as barriers from above and below and suitable comparison principles. For the reader's convenience, we decided to dedicate a separate section to each step.

\noindent Finally, we show that Theorem \ref{NTHEOREMCONVERGENCEINNEROUTERSETS} holds even when we relax the assumptions on the initial datum. In particular, we prove the following corollary.
\begin{cor}\label{CORASYMBEHAMOREGENDATUM}
Let $m > 0$ and $p > 1$ such that $\gamma \geq 0$, and let $N \geq 1$.
Consider a Lebesgue-measurable initial datum $u_0 : \RR^N \to \RR$ with $u_0 \not \equiv 0$ and $0 \leq u_0 \leq 1$ and satisfying
\begin{equation}\label{eq:ASSUMPTIONSONTHEINITIALDATUMEXPONENTIALGGGEQ0}
u_0(x) \leq a_0 \exp\Big(-\nu_{\ast}^{-1}f'(0)^{\frac{1}{p}}|x|\Big) \quad \text{for }\; |x| \sim \infty
\qquad \text{ if } \gamma > 0,
\end{equation}
or
\begin{equation}\label{eq:ASSUMPTIONSONTHEINITIALDATUMEXPONENTIALGG0}
u_0(x) \leq a_0 |x|^{\frac{2}{p}}\exp\Big(-m^{\frac{2-p}{p}}f'(0)^{\frac{1}{p}}|x|\Big) \quad \text{for }\; |x| \sim \infty
\qquad \text{ if } \gamma = 0,
\end{equation}
where $\nu_{\ast} = \nu_{\ast}(m,p)$ is the critical speed for $f'(0) = 1$, and $a_0 > 0$. Then the statements (i) and (ii) of Theorem \ref{NTHEOREMCONVERGENCEINNEROUTERSETS} hold for the solution of the initial-value problem \eqref{eq:REACTIONDIFFUSIONEQUATIONPLAPLACIAN}, \eqref{eq:ASSUMPTIONSONTHEINITIALDATUMEXPONENTIALGGGEQ0} if $\gamma > 0$, or \eqref{eq:REACTIONDIFFUSIONEQUATIONPLAPLACIAN}, \eqref{eq:ASSUMPTIONSONTHEINITIALDATUMEXPONENTIALGG0} if $\gamma = 0$, respectively.
\end{cor}

In Section \ref{SECTIONSTRONGREACTION} we present an extension of Theorem \ref{THEOREMSLOWDIFFUSIONSTRONGREACTIONPOROUSMEDIUM} to the doubly nonlinear case when $\gamma \geq 0$, studying the TWs ``produced'' by the combination of the doubly nonlinear diffusion term with a strong reaction term of the type $f(u) = u^n(1-u)$, $n \in \RR$.

\medskip

In Section \ref{Appendix1} we give a detailed proof of formula \eqref{eq:ASYMPTOTICTALESTWPSEUDOCAST} which describes the asymptotic behaviour of the TW solution with critical speed of propagation when $\gamma = 0$, and is linked to the choice of the initial datum \eqref{eq:ASSUMPTIONSONTHEINITIALDATUMEXPONENTIALGG0}.

\medskip

In Section \ref{SECTIONCOMMENTANDOPENPROBLEMS} we conclude the paper by discussing some issues related to our study. Firstly, as we have anticipated earlier, we briefly report the results obtained in \cite{B1:art,B2:art} and \cite{Hamel-N-R-R:art} on the behaviour of the general solutions of problem \eqref{eq:ONEDIMENSIONALKPPEQUATIONINTRO} on the moving coordinate $\xi = x - c_{\ast}t$. Then we present some results and open problems concerning the F-KPP theory when we deal with a ``fast diffusion'' term. Finally, we comment briefly what our results tell about the limit cases $m = 0$, $p = \infty$ and $m = \infty$, $p = 1$ with the condition $m(p-1) = 1$ ($\gamma = 0$).
%
%
%
%
%
%
%
\section{\texorpdfstring{\boldmath}{}Existence of Travelling Wave solutions: case \texorpdfstring{$\gamma > 0$}{gamma}}\label{SECTIONEXISTENCEOFTWS}
This section is devoted to the proof of Theorem \ref{THEOREMEXISTENCEOFTWS} by performing a very detailed analysis of the existence of travelling waves by ODE techniques. We work in one space dimension. We fix $m > 0$ and $p > 1$ such that $\gamma > 0$ and consider the nonlinear reaction-diffusion equation
\begin{equation}\label{eq:REACDIFFPLAPLACIANTWS}
\partial_tu = \partial_x(|\partial_xu^m|^{p-2}\partial_xu^m) + f(u) \quad \text{in } \RR\times(0,\infty),
\end{equation}
where the reaction term $f(\cdot)$ satisfies \eqref{eq:ASSUMPTIONSONTHEREACTIONTERM}.
In particular, we look for admissible TW solutions for equation \eqref{eq:REACDIFFPLAPLACIANTWS}, i.e., solutions of the form $u(x,t) = \varphi(\xi)$, where $\xi = x + ct$, $c > 0$ and $\varphi(\cdot)$ satisfying \eqref{eq:CONDITIONONPHIADMISSIBLETWINTRO}. Recall that we look for solution with profile $0 \leq \varphi \leq 1$, $\varphi(-\infty) = 0$, $\varphi(\infty) = 1$ and $\varphi' \geq 0$, and there is a second option in which $\varphi' \leq 0$ and the wave moves in the opposite direction, but we can skip reference to this case that is obtained from the previous one by reflection. Hence, we have to study the second order ODE
\begin{equation}\label{eq:PROFILEEQUATIONTWS}
c\varphi' = [|(\varphi^m)'|^{p-2}(\varphi^m)']' + f(\varphi) \quad \text{in } \RR,
\end{equation}
where the notation $\varphi'$ indicates the derivative of $\varphi$ with respect to the variable $\xi$.

The simplest approach consists in performing the standard change of variables $X = \varphi$ and $Y = \varphi'$, thus transforming the second-order ODE \eqref{eq:PROFILEEQUATIONTWS} into the system of two first-order ODEs
\[
\frac{dX}{d\xi} = Y, \quad\quad (p-1)m^{p-1}X^{\mu}|Y|^{p-2}\frac{dY}{d\xi} = cY - f(X) - \mu m^{p-1}X^{\mu-1}|Y|^p,
\]
which can be re-written after re-parametrization as the less singular system
\begin{equation}\label{SECONDSYSTEMNONSINGULARTWS}
\frac{dX}{ds} = (p-1)m^{p-1}X^{\mu}|Y|^{p-2}Y, \quad\quad \frac{dY}{ds} = cY - f(X) - \mu m^{p-1}X^{\mu-1}|Y|^p,
\end{equation}
where we used the re-parametrization $d\xi = (p-1)m^{p-1}X^{\mu}|Y|^{p-2}ds$. Note that both system are equivalent for $Y \ne 0$ but, at least in the case $\mu > 1$, the second one has two critical points $(0,0)$ and $(1,0)$.

\noindent This setting seems to be convenient since proving the existence of an admissible TW for \eqref{eq:REACDIFFPLAPLACIANTWS} corresponds to showing the existence of a trajectory in the region of the $(X,Y)$-plane where $0 \leq X \leq 1$ and $Y \geq 0$, and joining the critical points $(0,0)$ and $(1,0)$. More precisely, the desired trajectories must ``come out'' of $(0,0)$ and ``enter'' $(1,0)$. Let us point out that, contrary to the linear case (see for example \cite{Aronson2:art, Aro-Wein1:art} and \cite{K-P-P:art}), we face a more complicated problem. Indeed, it turns out that in the nonlinear diffusion case, the Lyapunov linearization method (\cite{Debnath:book}, Chapter 8) used to analyze the local behaviour of the trajectories near the critical points cannot be applied since system \eqref{SECONDSYSTEMNONSINGULARTWS} is heavily nonlinear.

System \eqref{SECONDSYSTEMNONSINGULARTWS} admits a large collection of trajectories and corresponding TWs. Hence, the very problem is to understand how to discern between different types of trajectories ``coming out'' of the steady state $(0,0)$. The flow around this point has a complicated structure, so, following the methods inspired in the study of the Porous Medium case (see \cite{DePablo-Vazquez1:art}), we introduce the new variables
\begin{equation}\label{eq:NONSTANDARDCHANGEOFVARIABLESTWS}
X = \varphi \qquad \text{ and } \qquad Z = \Bigg(\frac{m(p-1)}{\gamma}\varphi^{\frac{\gamma}{p-1}}\Bigg)' = mX^{\frac{\mu-1}{p-1}}X',
\end{equation}
using the fact that $\gamma=\mu+p-2$. These variables correspond to the density and the derivative of the pressure profile (see \cite{EstVaz:art}). Assuming only $X \geq 0$, we obtain the first-order ODE system
\begin{equation}\label{eq:SYSTEMNONSINGULARTWSXI}
\frac{dX}{d\xi} = (1/m)X^{\frac{1-\mu}{p-1}}Z, \quad\quad m(p-1)X^{\frac{\gamma}{p-1}}|Z|^{p-2}\, \frac{dZ}{d\xi} = cZ - |Z|^p - mX^{\frac{\mu-1}{p-1}}f(X),
\end{equation}
that again can be re-written as the non-singular system
\begin{equation}\label{eq:SYSTEMNONSINGULARTWS}
\frac{dX}{d\tau} = (p-1)X|Z|^{p-2} Z, \quad\quad \frac{dZ}{d\tau} = cZ - |Z|^p - mX^{\frac{\mu - 1}{p-1}}f(X)
\end{equation}
where we use the new re-parametrization  $d\xi = m(p-1)X^{\frac{\gamma}{p-1}}|Z|^{p-2}d\tau$. Note that, since $0 < \gamma = \mu + p - 2$, the function
\[
f_{m,p}(X) = mX^{\frac{\mu - 1}{p-1}}f(X) = mX^{\frac{\gamma}{p-1}-1}f(X)
\]
is continuous for all $0 \leq X \leq 1$ and satisfies $f_{m,p}(0) = 0 = f_{m,p}(1)$. Hence, the critical points are now three: $S = (1,0)$, $O = (0,0)$ and $R_c = (0,c^{1/(p-1)})$. The last one will play an important role.

\noindent It is clear that the change of variables \eqref{eq:SYSTEMNONSINGULARTWS} ``splits'' the critical point $(0,0)$ of system \eqref{SECONDSYSTEMNONSINGULARTWS} into two critical points: $O$ and $R_c$ of the new system. The idea is that formula \eqref{eq:NONSTANDARDCHANGEOFVARIABLESTWS} is the unique change of variables that allows us to separate the orbits corresponding to finite and positive TWs. In fact, we will show that the connections between $O$ and $S$ correspond to positive TW (they exist only if $c > c_{\ast}$ where $c_{\ast}$ is the critical speed of propagation in the statement of Theorem \ref{THEOREMEXISTENCEOFTWS}) while the connection between $R_c$ and $S$ corresponds to a finite TW ($c = c_{\ast}$).

As we have mentioned in the introduction, we stress that these heuristic explanations can be formally proved (see Theorem \ref{THEOREMEXISTENCEOFTWS}) and represent the main difference between the case $m = 1$, $p = 2$ and $m > 0$, $p > 1$ (with $\gamma > 0$). In the linear case (\cite{Aronson2:art, Aro-Wein1:art} and \cite{K-P-P:art}), all the TWs are positive and the analysis can be performed using the Lyapunov method while, in the nonlinear case (see also \cite{Aronson1:art} for a comparison with the Porous Medium equation), the slow diffusion term ``causes'' the apparition of finite TWs and free boundaries but the Lyapunov method cannot be applied.
\paragraph{Proof of Theorem \ref{THEOREMEXISTENCEOFTWS}.} In order to carry out the plan introduced before, we focus on the range $0 \leq X \leq 1$, we set $C = c^{1/(p-1)}$ and we note that, thanks to the assumption $\gamma > 0$ and the concavity of $f(\cdot)$, the function $f_{m,p}(\cdot)$ has a unique maximum point $X_{m,p} \in [0,1]$ (this fact is not completely trivial but we leave it as an exercise for the interested reader). Let $F_{m,p}$ the maximum of the function $f_{m,p}(\cdot)$ in $[0,1]$. We study the \emph{equation of the trajectories}
\begin{equation}\label{eq:EQUATIONTRAJECTORIESTWS}
\frac{dZ}{dX} = \frac{cZ - |Z|^p - f_{m,p}(X)}{(p-1)X|Z|^{p-2}Z} := H(X,Z;c)
\end{equation}
after eliminating the parameter $\tau$, and look for solutions defined for $0 \leq X \leq 1$ and linking the critical points $R_c$ and $S$  for some $c > 0$. The main difference with respect to the Porous Medium and the linear case (see \cite{Aro-Wein1:art} and \cite{DePablo-Vazquez1:art}) is that the critical points are all degenerate and it is impossible to describe the  local behaviour of the trajectories by linearizing the system \eqref{eq:SYSTEMNONSINGULARTWS}
around $S$. In what follows, we study some local and global properties of the equation \eqref{eq:EQUATIONTRAJECTORIESTWS} with other qualitative ODE methods in order to obtain a clear view of the graph of the trajectories.

\emph{Step1.} Firstly, we study the \emph{null isoclines} of system \eqref{eq:SYSTEMNONSINGULARTWS}, i.e., the set of the points $(X,Z)$ with $0 \leq X \leq 1$ and $Z \geq 0$ such that $H(X,Z;c) = 0$. Hence, we have to solve the equation
\begin{equation}\label{eq:ISOCLINESEQUATIONTWS}
c\widetilde{Z} - \widetilde{Z}^p = f_{m,p}(X) \quad \text{ in } [0,1]\times[0,\infty).
\end{equation}
Defining
\[
c_0 = c_0(m,p) := p\Bigg(\frac{F_{m,p}}{p-1}\Bigg)^{(p-1)/p},
\]
it is simple to see that for $c < c_0$ the graph of the isoclines is composed by two branches joining the point $O$ with $R_c$ and the point $(1,0)$ with $(1,C)$, while for $c > c_0$ the branches link the point $O$ with $S$ and the point $R_c$ with $(1,C)$. As $c$ approaches the value $c_0$ the branches move nearer and they touch in the point $(X_{m,p},(c_0/p)^{1/(p-1)})$ for $c = c_0$.
Note that the value $c_0(m,p)$ is critical in the study of the isoclines and it is found by imposing
\begin{equation}\label{eq:HOWTOFINDC0}
\max_{\widetilde{Z} \in [0,C]} \{c\widetilde{Z} - \widetilde{Z}^p \} = F_{m,p}.
\end{equation}
\begin{figure}[!ht]
  \centering
  \includegraphics[scale=0.3]{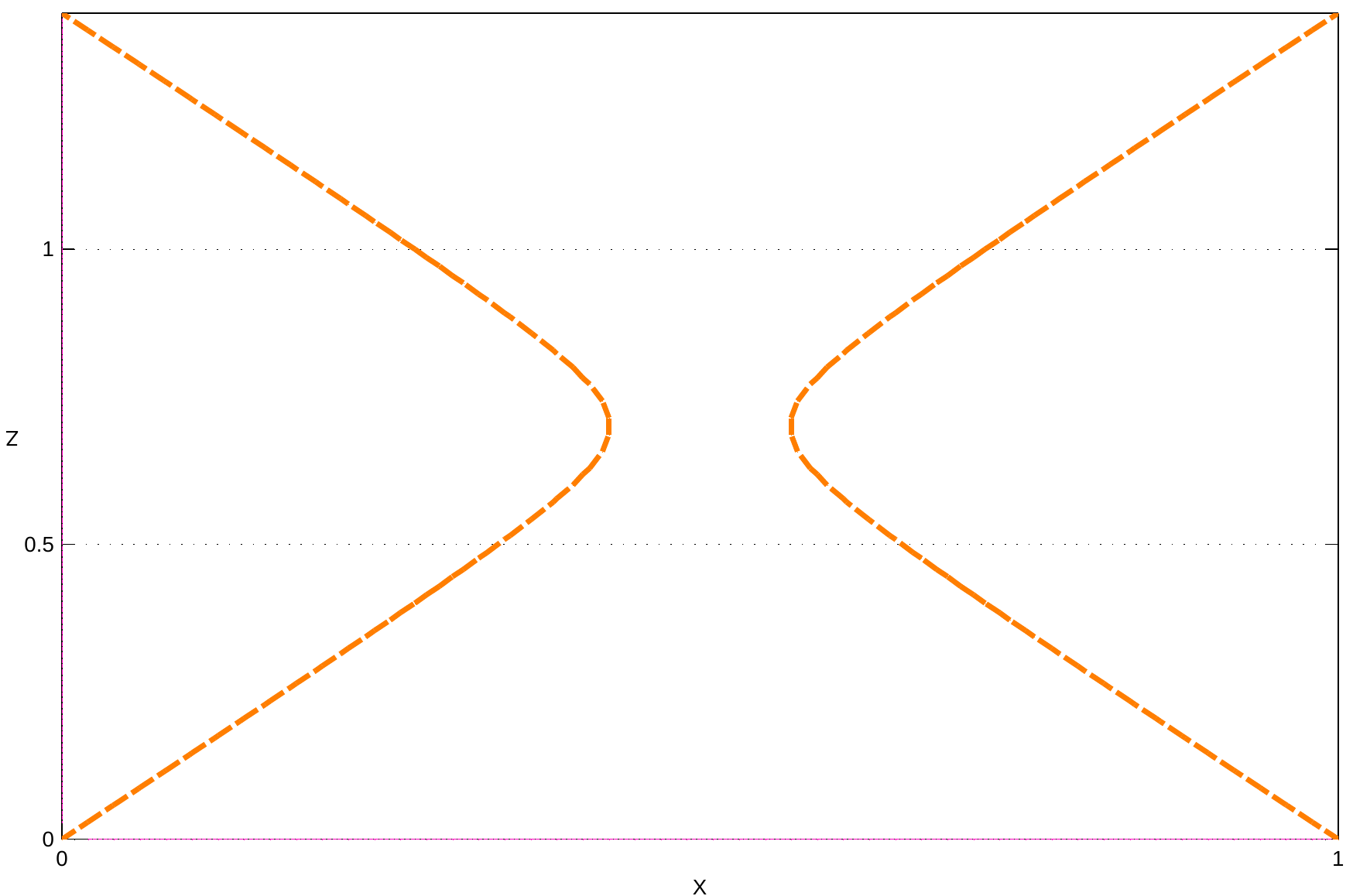}
  \includegraphics[scale=0.3]{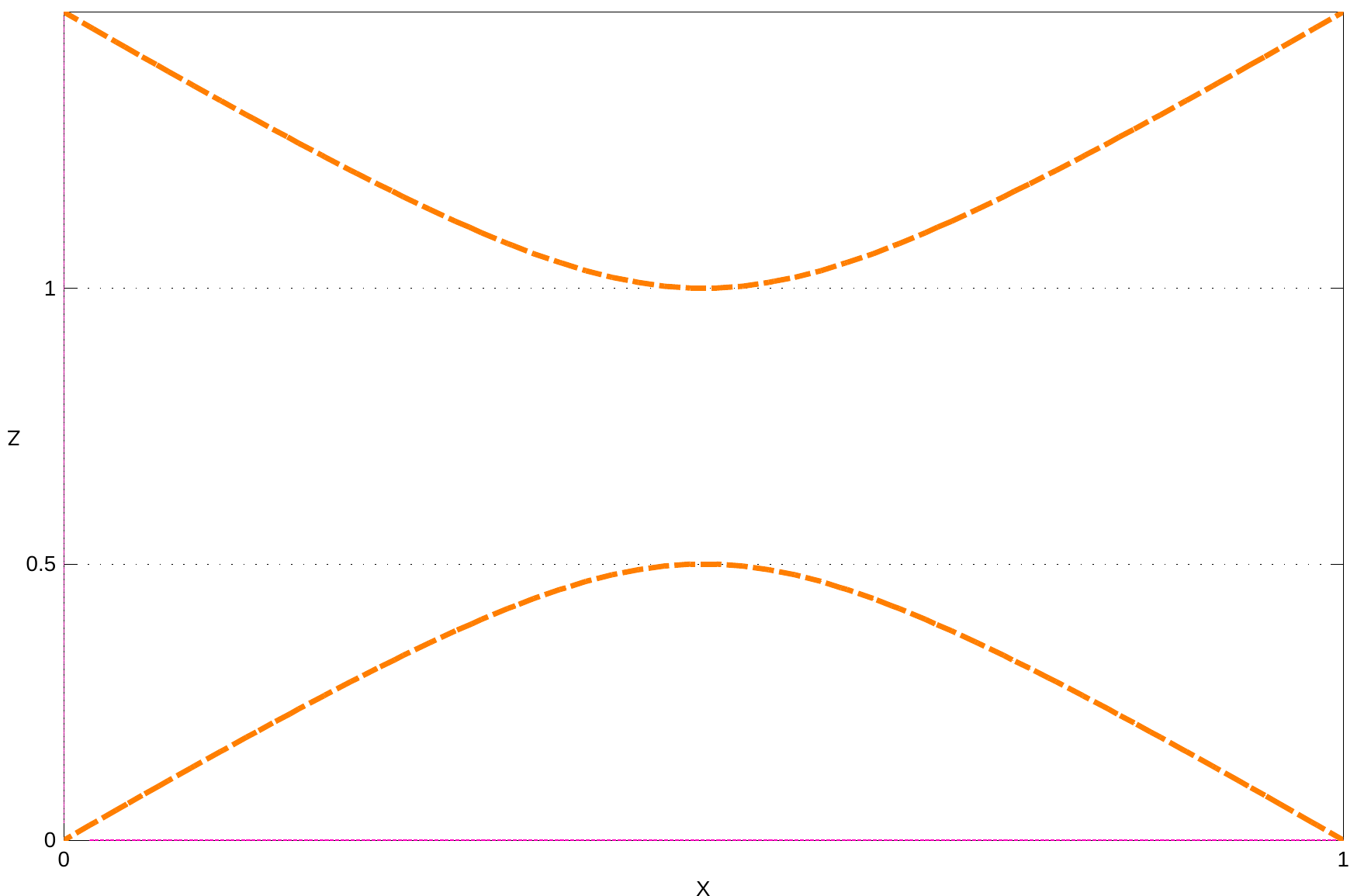}
  \caption{Null isoclines of system \eqref{eq:SYSTEMNONSINGULARTWS}: Case $c < c_0$ and $c > c_0$.}\label{fig:NULLISOCLINESDNL}
\end{figure}

For $c < c_0$, the trajectories in the region between the two branches have negative slopes, while those in the region between the Z-axis and the left branch and between the right branch and the line $X=1$ have positive slopes. Conversely, for $c > c_0$ the trajectories in the region between the two branches have positive slopes while, those in the region between the bottom-branch and the X-axis and between the line $Z=C$ and the top-branch have negative slopes. Finally, for $c = c_0$, it is simple to see that in the regions between the bottom-branch and the X-axis and between the line $Z=C$ and the top-branch the trajectories have negative slopes and they have positive slopes in the rest of the rectangle $[0,1]\times[0,C]$. We conclude this paragraph noting that for all $c > 0$ the trajectories have negative slopes for $Z > C$ and positive slopes for $Z < 0$.

\emph{Step2.} We prove the existence and the uniqueness of solutions of \eqref{eq:EQUATIONTRAJECTORIESTWS} ``entering'' in the point $S$.

\noindent \emph{Case $p = 2$.} If $p = 2$, it is not difficult to linearize system \eqref{eq:SYSTEMNONSINGULARTWS} through Lyapunov method and showing that the point $S$ is a saddle. Moreover, it follows that there exists exactly one locally linear trajectory in the region $[0,1]\times[0,\infty)$ ``entering'' in $S$ with slope $\lambda_S = (c - \sqrt{c^2 - 4mf'(1)})/2$.

\noindent \emph{Case $p > 2$.} Substituting the expression $Z = \lambda(1-X)$, $\lambda > 0$ in the equation of trajectories \eqref{eq:EQUATIONTRAJECTORIESTWS} and taking $X \sim 1$ we get
\[
-\lambda = H(X,\lambda(1-X)) \sim \frac{c\lambda(1-X) + mf'(1)(1-X)}{(p-1)\lambda^{(p-1)}(1-X)^{p-1}}, \quad \text{for } X \sim 1
\]
which can be rewritten as
\[
-(p-1)\lambda^{p}(1-X)^{p-2} \sim c\lambda + mf'(1), \quad \text{for } X \sim 1.
\]
Since the left side goes to zero as $X \to 1$, the previous relation is satisfied only if $\lambda = -mc^{-1}f'(1) := \lambda_S^+ > 0$ (note that this coefficient coincides with the slope of the null isocline near $X = 1$). Hence, for $p > 2$, there exists at least one trajectory $T_c$ which ``enters'' into the point $S$ and it is linear near this critical point. Note that the approximation $Z(X) \sim \lambda_S^+ (1 - X)$ as $X \sim 1$ can be improved with high order terms. However, we are basically interested in proving the existence of a trajectory ``entering'' in $S$ and we avoid to present technical computations which can be performed by the interested reader.

\noindent To prove the uniqueness of the ``entering orbit'', we begin by showing that the trajectory $T_c$ (which satisfies $Z(X) \sim \lambda_S^+ (1 - X)$ as $X \sim 1$) is ``repulsive'' near $X = 1$. It is sufficient to prove that the partial derivative of the function $H$ with respect to the variable $Z$ is strictly positive when it is calculated on the trajectory $T_c$ and $X \sim 1$. It is straightforward to see that
\[
\begin{aligned}
\frac{\partial H}{\partial Z}(X,\lambda_S^+(1-X)) & \sim \frac{m}{D(X)}[-c\lambda_S^+(p-2)(1-X) - (p-1)f'(1)(1-X)] \\
& \sim - \frac{mf'(1)}{D(X)} \gg 0, \quad \text{for } X \sim 1
\end{aligned}
\]
where $D(X) = (p-1)(\lambda_S^+)^p(1-X)^{p-1}$. Hence, our trajectory is ``repulsive'' near the point $S$, i.e., there are not other \emph{locally linear} trajectories which ``enter'' into the point $S$ with slope $\lambda_S^+$.

\noindent Now, we claim that $T_c$ is the unique trajectory ``entering'' into the point $S$. Indeed, recalling the fact that $T_c$ has the same slope of the null isoclines near $S$ and using the fact that $\partial H/\partial Z > 0$ along $T_c$ for $X \sim 1$, it is evident that trajectories from $S$ above $T_c$ are not admitted. Furthermore, a simple calculation shows that for all $0 < X < 1$ fixed and $Z$ positive but small, the second derivative is negative:
\[
\frac{d^2Z}{dX^2} \sim - \frac{f_{m,p}(X)}{(p-1)X^2Z^{2(p-1)}}\bigg(\frac{c}{p-1} + \frac{f_{m,p}(X)}{Z} \bigg) < 0, \quad \text{for } Z \sim 0^+.
\]
Hence, we deduce our assertion observing that all trajectories ``entering'' into $S$ and lying below $T_c$ have to be convex, which is a contradiction with the previous formula.

\noindent\emph{Case $1 < p < 2$.} Proceeding as before, it is not difficult to prove that when $1 < p < 2$, there exists a trajectory ``entering'' in $S$ with local behaviour $Z \sim \lambda_S^-(1-X)^{2/p}$ for $X \sim 1$ and $\lambda_S^- := \{-pmf'(1)/[2(p-1)]\}^{1/p}$. Note that in this case the trajectory is not locally linear but presents a power-like local behaviour around $S$ with power greater than one. Moreover, exactly as in the case $p > 2$, it is simple to see that this trajectory is ``repulsive''. It remains to prove that our trajectory is the unique ``entering'' in $S$. Note that the computation of the second derivative is still valid. Moreover, in order to show that there are no orbits between our trajectory and the branch of the null isoclines, we consider the two-parameters family of curves:
\[
Z_a(x) = \lambda_S^-(1-X)^a, \qquad  1 < a < 2/p
\]
and we use an argument with invariant regions. We compute the derivative
\[
\frac{dZ_a/dX}{dZ/dX}(X,Z_a(X);c) = \frac{dZ_a/dX(X)}{H(X,Z_a(X);c)} \sim \frac{-a(p-1)(\lambda_S^-)^p}{mf'(1)}(1-X)^{ap-2} \sim +\infty \quad \text{as } X \sim 1
\]
for all $1 < a < 2/p$. This means that the ``flux'' derivative along the curve $Z_a = Z_a(X)$ is very small respect to the derivative of the curve when $X \sim 1$, i.e., the trajectories have horizontal slopes respect to the ones of $Z_a = Z_a(X)$. So, since $1 < a < 2/p$ is arbitrary, it is possible to conclude the uniqueness of our trajectory ``entering'' in $S$. We anticipate we will use a similar technique in \emph{Step5}.

\emph{Step3.} We are ready to prove that there exists a unique value $c = c_{\ast}$ and a unique trajectory $T_{c_{\ast}}$ linking $S$ and $R_{c_{\ast}}$.
\\
Before proceeding, we have to introduce some notations. Let $\Gamma_1$ be the left-branch of the isoclines in the case $0 < c < c_0$ (note that the study of the isoclines carried out in \emph{Step1} and the analysis in \emph{Step5} tell us that there are not trajectories linking the points $S$ and $R_c$ for $c \geq c_0$), define $\Gamma_2 := \{(X,C): 0 \leq X \leq 1\}$ and consider $\Gamma_c := \Gamma_1\cup R_c \cup \Gamma_2$. Finally, note that the trajectories are strictly monotone respect to the parameter $c$, i.e.,
\begin{equation}\label{eq:PARAMETERDERIVATIVEHTW}
\frac{\partial H}{\partial c} > 0 \quad \text{in } [0,1]\times \RR \times (0,\infty).
\end{equation}
Note that \eqref{eq:PARAMETERDERIVATIVEHTW} tells us that if $c_1 < c_2$ then $T_{c_1} > T_{c_2}$, i.e., we have strict monotonicity of the trajectories from $S$ respect to the parameter $c > 0$.
\\
Now, fix $0 < \overline{c} < c_0$ and let $T_{\overline{c}}$ be the linear trajectory corresponding to the value $\overline{c}$ from the point $S$. Since $T_{\overline{c}}$ come from the point $S$ and it lies in the region in which the slope is negative, it must join $S$ with $\Gamma_{\overline{c}}$, i.e., there exists a point $(\overline{X},\overline{Z}) \in \Gamma_{\overline{c}} \cap T_{\overline{c}}$. We have the following possibilities:

\noindent $\bullet$ If $(\overline{X},\overline{Z}) = R_{\overline{c}}$, then we have $c_{\ast} = \overline{c}$ and the uniqueness follows by \eqref{eq:PARAMETERDERIVATIVEHTW}.

\noindent $\bullet$ If $(\overline{X},\overline{Z}) \in \Gamma_1$, using \eqref{eq:PARAMETERDERIVATIVEHTW} it follows that there exists $0 < \overline{c}_2 < \overline{c}$ and a corresponding point $(\overline{X}_2,\overline{Z}_2) \in \Gamma_2 \cap T_{\overline{c}_2}$. Moreover, we have that for all $\overline{c} < c < c_0$ the trajectory $T_c$ links $S$ with $\Gamma_1$ and for all $0 < c < \overline{c}_2$ the trajectory $T_c$ links $S$ with $\Gamma_2$. Hence, from the continuity of the trajectories respect to the parameter $c$, there exists $\overline{c}_2 < c_{\ast} < \overline{c}$ such that $T_{c_{\ast}}$ joins $S$ with $R_{c_{\ast}}$ and the uniqueness follows from the strict monotonicity \eqref{eq:PARAMETERDERIVATIVEHTW}.

\noindent $\bullet$ If $(\overline{X},\overline{Z}) \in \Gamma_2$, using \eqref{eq:PARAMETERDERIVATIVEHTW} it follows that there exists $\overline{c} < \overline{c}_1 < c_0$ and a corresponding point $(\overline{X}_1,\overline{Z}_1) \in \Gamma_1 \cap T_{\overline{c}_1}$. Moreover, we have that for all $\overline{c}_1 < c < c_0$ the trajectory $T_c$ links $S$ with $\Gamma_1$ and for all $0 < c < \overline{c}$ the trajectory $T_c$ links $S$ with $\Gamma_2$. Hence, there exists a unique $\overline{c} < c_{\ast} < \overline{c}_1$ such that $T_{c_{\ast}}$ joins $S$ with $R_{c_{\ast}}$.

\noindent Thus, we can conclude that there exists exactly one value $c = c_{\ast}(m,p) < c_0(m,p)$ with corresponding trajectory $T_{c_{\ast}}$ joining the points $S$ and $R_{c_{\ast}}$. Moreover, since in  \emph{Step2} we have showed that the trajectory from the point $S$ is unique, it follows that the trajectory $T_{c_{\ast}}$ is unique too. We underline that our argument is completely qualitative and based on a topological observation: we proved the existence of two numbers $0 < \overline{c}_2 < \overline{c}_1 < c_0$ such that for all $0 < c < \overline{c}_2$, $T_c$ links $S$ with $\Gamma_2$, while for all $\overline{c}_1 < c < c_0$, $T_c$ links $S$ with $\Gamma_1$. Hence, since the trajectories are continuous respect with the parameter $c$ it follows the existence of a value $c_{\ast}$ such that $T_{c_{\ast}}$ links $S$ with $\Gamma_1\cap\Gamma_2 = R_{c_{\ast}}$.

\noindent Finally, it remains to prove that the TW corresponding to $c_{\ast}$ is finite i.e., it reaches the point $u = 0$ in finite ``time'' while the point $u = 1$ in infinite ``time'' (here the time is intended in the sense of the profile, i.e. the ``time'' is measured in terms of the variable $\xi$). In order to see this, it is sufficient to integrate the first equation in \eqref{eq:SYSTEMNONSINGULARTWSXI} by separation of variable between $X_0$ and $X_1$:
\begin{equation}\label{eq:DIFFERENTIALRELATIONTW}
\xi_1 - \xi_0 = m\int_{X_0}^{X_1}\frac{dX}{ZX^{(1-\mu)/(p-1)}}.
\end{equation}
Since $\gamma = \mu + p - 2$, it is clear that for $\gamma > 0$ we have $(1-\mu)/(p-1)<1$, so that for all $0 < X_1 < 1$, the integral of the right member is finite when $X_0 \to 0$ and so, our TW reaches the steady state $u=0$ in finite time. This is a crucial computation that shows that \emph{finite propagation} holds in this case. Conversely, recalling the behaviour of the trajectory $T_{c_{\ast}}$ near the point $S$ (see \emph{Step2}), it is not difficult to see that for all $0 < X_0 < 1$, the integral is infinite when $X_1 \to 1$ and so, the TW gets to the steady state $u = 1$ in infinite time.

We conclude this step with a brief analysis of the remaining trajectories. This is really simple once one note that the trajectory $T_{c_{\ast}}$ acts as barrier dividing the set $[0,1]\times\RR$ into subsets, one below and one above $T_{c_{\ast}}$. The trajectories in the subset below link the point $O$ with $R_{-\infty} := (0,-\infty)$ (the same methods we use in \emph{Step5} apply here), have slope zero on the left branch of the isoclines, are concave for $Z > 0$ and increasing for $Z < 0$. We name these trajectories CS-TWs (change sign TWs) of type 1. On the other hand, the trajectories above $T_{c_{\ast}}$ recall ``parabolas'' (see \emph{Step4}) connecting the points $R_{\infty} := (0,\infty)$ and $S_Z = (1,Z)$ for some $Z > 0$ and having slope zero on the right branch of the isoclines. Finally, note that, in this last case, the trajectories lying in the region $[0,1]\times[C_{\ast},\infty)$ are always decreasing.

\emph{Step4.} In the next paragraph, we show that there are no admissible TW solutions when $0 < c < c_{\ast}$. Suppose by contradiction that there exists $0 <c < c_{\ast}$ such that the corresponding trajectory $T_c$ joins $S$ and $R_{\overline{Z}}$ for some $\overline{Z} \geq 0$. Then, by \eqref{eq:PARAMETERDERIVATIVEHTW}, it must be $\overline{Z} > C_{\ast}$. Moreover, since $H(X,Z;c) < 0$ for all $Z > C$, there exists a ``right neighbourhood'' of $R_{\overline{Z}}$ in which the solution $Z = Z(X)$ corresponding to the trajectory $T_c$ is invertible and the function $X = X(Z)$ has derivative
\begin{equation}\label{eq:INVERSEEQUATIONTRAJECTORIESTWS}
\frac{dX}{dZ} = \frac{(p-1)XZ^{p-1}}{cZ - Z^p - f_{m,p}(X)} := K(X,Z;c).
\end{equation}
Choosing the neighbourhood $B_{\delta}(R_{\overline{Z}}) = \{(X,Z): \; X^2 + (Z - \overline{Z})^2 < \delta, \; X \geq 0 \}$, where $\delta > 0$ is small enough (for example, $\delta \leq \min\{1, \overline{Z} - C_{\ast}\}$), it is simple to obtain
\[
\Bigg|\frac{K(X,Z;c)}{X}\Bigg| \leq I \quad \text{in } B_{\delta}(R_{\overline{Z}}),
\]
where $I = I_{\overline{Z},\delta,m,p} > 0$ depends on $\overline{Z}$, $\delta$, $m$ and $p$. The previous inequality follows from the fact the quantity $|cZ - Z^p - f_{m,p}(X)|$ is greater than a positive constant in $B_{\delta}(R_{\overline{Z}})$ (see \emph{Step1}) and $Z$ is (of course) bounded in $B_{\delta}(R_{\overline{Z}})$. This means that the function $K(X,Z)$ is sub-linear respect with the variable $X$, uniformly in $Z$ in $B_{\delta}(R_{\overline{Z}})$ and this is sufficient to guarantee the local uniqueness of the solution near $R_{\overline{Z}}$. Since the null function solves \eqref{eq:INVERSEEQUATIONTRAJECTORIESTWS} with initial datum 0, from the uniqueness of this solution, it follows that $X = X(Z)$ is identically zero too and it cannot be invertible. This contradiction assures there are no TW solutions for $c < c_{\ast}$.

As we did in the previous step, we explain the qualitative properties of the trajectories. In the case $c < c_{\ast}$, the ``zoo'' of the trajectories is more various. First of all, the previous analysis shows that we have a connection between the points $S$ and $R_{\infty}$ always decreasing. Moreover, below this connection, there is a family of CS-TWs of type 1 linking the points $O$ and $R_{-\infty}$ and family of CS-TWs of type 2, i.e., trajectories linking $R_{\infty}$ with $R_{-\infty}$, decreasing for $Z>0$ and increasing for $Z<0$. Furthermore, for topological reasons, there exists a trajectory between the family of CS-TWs of type 1 and CS-TWs of type 2 which link the critical point $R_c$ and $R_{-\infty}$. Finally, we find again the ``parabolas'' described in \emph{Step3}.

\emph{Step5.} Finally, we focus on the case $c > c_{\ast}$. We follow the procedure used in \emph{Step4}, supposing by contradiction that there exists $c > c_{\ast}$ such that the corresponding trajectory $T_c$ joins $S$ and $R_{\overline{Z}}$ for some $\overline{Z} > 0$ (in this case it must be $\overline{Z} < C_{\ast}$). Note that, to be precise, we should treat separately the cases $c_{\ast} < c < c_0$ and $c \geq c_0$ since the phase plane changes markedly, but, however, our argument works independently of this distinction.
\\
Again, we want to prove the sub-linearity of the function $K(X,Z;c)$ respect with the variable $X$, uniformly in $Z$ in a ``right'' neighbourhood of $R_{\overline{Z}}$. Define
\[
\begin{aligned}
&\mathcal{R} := \{(X,Z): H(X,Z;c) > 0 \text{ for } 0 \leq X \leq 1, \; 0 \leq Z \leq C_{\ast}\} \\
&B_{\delta}(R_{\overline{Z}}) = \{(X,Z): \; X^2 + (Z - \overline{Z})^2 \leq \delta, \; X \geq 0 \},
\end{aligned}
\]
where $\delta > 0$ is taken small enough such that $B_{\delta}(R_{\overline{Z}}) \subset \mathcal{R}$. Hence, proceeding as in the previous step, we can state that, in $B_{\delta}(R_{\overline{Z}})$, the quantity $|cZ - Z^p - f_{m,p}(X)|$ is greater than a positive constant and so, it is simple to get
\[
\Bigg|\frac{K(X,Z;c)}{X}\Bigg| \leq I \quad \text{in } B_{\delta}(R_{\overline{Z}}),
\]
where $I = I_{\overline{Z},\delta,m,p} > 0$ depends only on $\overline{Z}$, $\delta$, $m$ and $p$.
Hence, reasoning as before, we obtain that the trajectories cannot ``touch'' the $Z$-axis. This means that the point $O$ ``attracts'' the trajectories and so for all $c > c_{\ast}$ there exists a connection between the points $S$ and $O$, i.e., a TW solution.

\noindent Now, in order to prove that these TWs are positive, we show that all trajectories approach the branch of the isoclines near $O$ given by equation \eqref{eq:ISOCLINESEQUATIONTWS}. First of all, from the equation of the null isoclines \eqref{eq:ISOCLINESEQUATIONTWS}, it is simple to see that the branch of the isoclines satisfies $\widetilde{Z}(X) \sim \widetilde{\lambda} X^{\gamma/(p-1)}$ as $X \sim 0$, where $\widetilde{\lambda} := mf'(0)/c$. Now, as we did in \emph{Step1}, we use an argument with invariant regions. Consider the one-parameter family
\[
Z_{\lambda}(X) = \lambda X^{\frac{\gamma}{p-1}}, \qquad \widetilde{\lambda} < \lambda < \infty.
\]
With straightforward calculations as in \emph{Step1}, it is simple to obtain
\[
\frac{dZ_{\lambda}/dX}{dZ/dX}(X,Z_{\lambda}(X);c) = \frac{dZ_{\lambda}/dX (X)}{H(X,Z_{\lambda}(X);c)} \sim \frac{\gamma \lambda^p}{c(\lambda - \widetilde{\lambda})}X^{\gamma} \sim 0 \quad \text{as } X \sim 0,
\]
i.e., for all $\lambda > \widetilde{\lambda}$ and for small values of the variable $X$, the derivative of the trajectories along the curve $Z_{\lambda} = Z_{\lambda}(X)$ is infinitely larger than the derivative of the same curve. This fact implies that for all $\lambda > \widetilde{\lambda}$, all trajectories $Z = Z(X)$ satisfy $\widetilde{Z}(X) \leq Z(X) \leq Z_{\lambda}(X)$ as $X \sim 0$ and so, since $\lambda$ can be chosen arbitrarily near to $\widetilde{\lambda}$ it follows $Z(X) \sim \widetilde{Z}(X)$ as $X \sim 0$.

\noindent We may now integrate the first differential equation in \eqref{eq:SYSTEMNONSINGULARTWSXI} by separation of variables. Since $\widetilde{Z}(X) \sim \widetilde{\lambda} X^{\gamma/(p-1)}$ and $\gamma + 1-\mu = p-1$, we have that $\widetilde{Z}(X)X^{(1-\mu)/(p-1)} = \widetilde{\lambda}X$. Hence, using the fact that $Z(X) \sim \widetilde{Z}(X)$ as $X \sim 0$, we get
\[
\xi - \xi_0 = m\int_{X_0}^{X} \frac{dX}{ZX^{(1-\mu)/(p-1)}} \sim m\int_{X_0}^{X} \frac{dX}{\widetilde{Z}X^{(1-\mu)/(p-1)}} \sim \frac{c}{f'(0)}\int_{X_0}^{X} \frac{dX}{X} \quad \text{for } X \sim 0.
\]
where $0 < X_0:=X(\xi_0) < X$ is fixed. Thus, since the last integral is divergent in $X = 0$, we deduce that the admissible TWs reach the level $u = 0$ at the time $\xi = -\infty$, i.e., they are positive.

\noindent Moreover, we can describe the exact shape of these TWs when $X \sim 0$. Indeed, from the previous formula we have that $\xi - \xi_0 \sim (c/f'(0))\ln(X/X_0)$ for $X \sim 0$, which can be easily re-written as
\begin{equation}\label{eq:ASYMBEHASUPERTWGGPOS}
X(\xi) \sim a_0\exp\Big(\frac{f'(0)}{c}\xi\Big) = a_0\exp\Big(\frac{f'(0)^{1/p}}{\nu}\xi\Big), \quad \text{for } \xi \sim -\infty
\end{equation}
where $a_0$ is a fixed positive constant (depending on $\xi_0$) and $\nu > 0$ is the speed with $f'(0)=1$. We recall that $X(\xi) = \varphi(x+ct)$ which is $u(x,t)$ for the TW.

We conclude this long analysis, describing the ``zoo'' of the trajectories obtained in the $(X,Z)$-plane. As we observed in the case $c = c_{\ast}$ the TW joining $S$ and $O$ represents a barrier and divide the region $[0,1]\times\RR$ in two subsets. We obtain again the CS-TWs of type 1 in the region below the positive TW and the ``parabolas'' in the region above it. Moreover, if $c_{\ast} < c < c_0$, we have a family of trajectories which ``come'' from $O$, are decreasing in the region between the branches of the isoclines and increasing in the regions between the $Z$-axis and the left branch and between the right branch and the line $X = 1$. Performing the complete phase plane analysis it is possible to see that these trajectories are increasing for $X \geq1 $ and go to infinity as $X \to \infty$. We name them infinite TWs of type 1. If $c \geq c_0$, the analysis is similar except for the fact that there exist a family (or exactly one trajectory in the case $c = c_0$) of infinite TWs of type 2, i.e., trajectories ``coming'' from $O$ and increasing for all $X \geq 0$. Finally, we have a trajectory from the point $R_c$ trapped between the family of infinite TWs and the family of ``parabolas'' which goes to infinity as $X \to \infty$. $\Box$

\begin{figure}[!ht]
  \centering
  \includegraphics[scale=0.3]{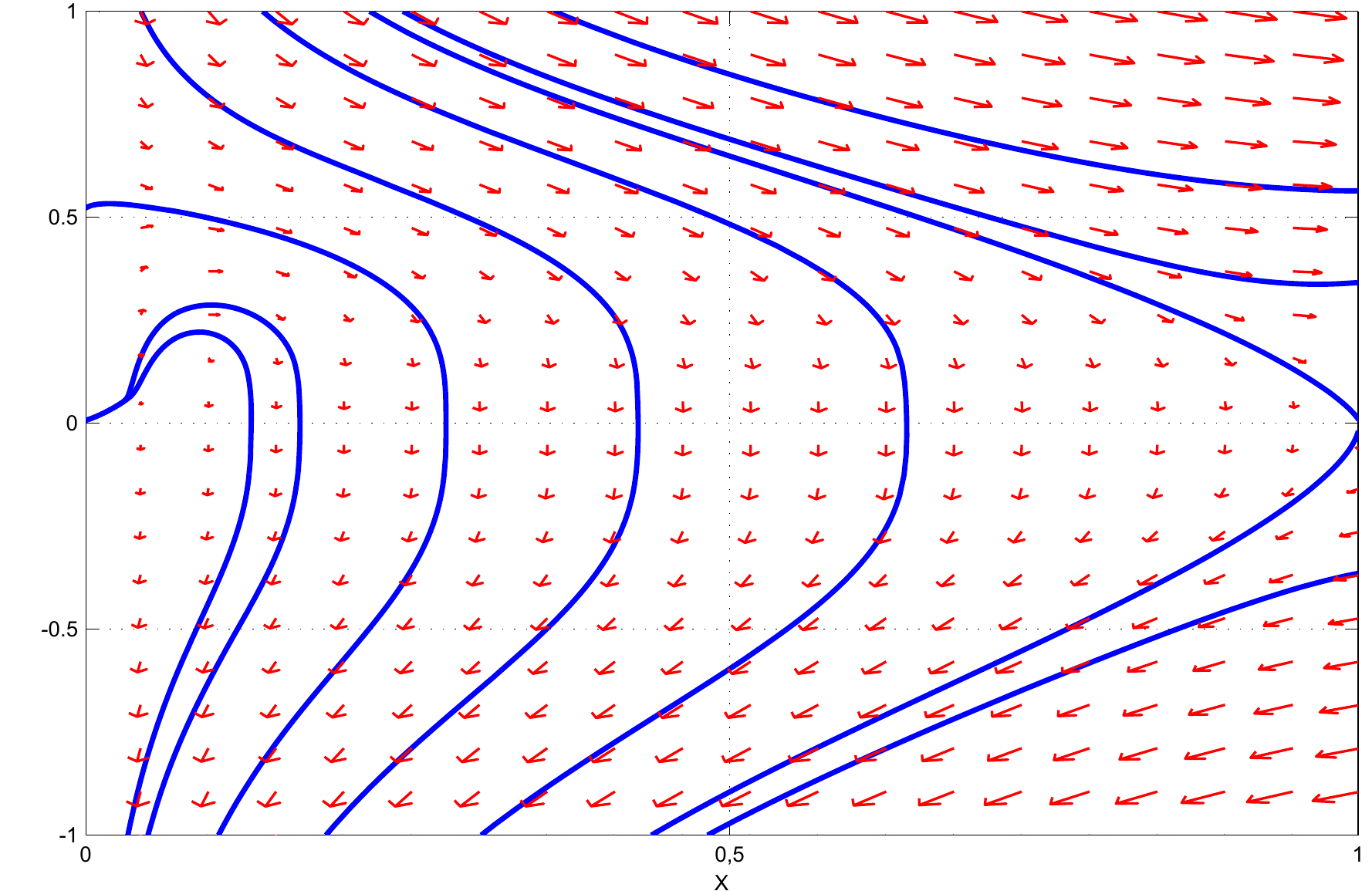}
  \includegraphics[scale=0.3]{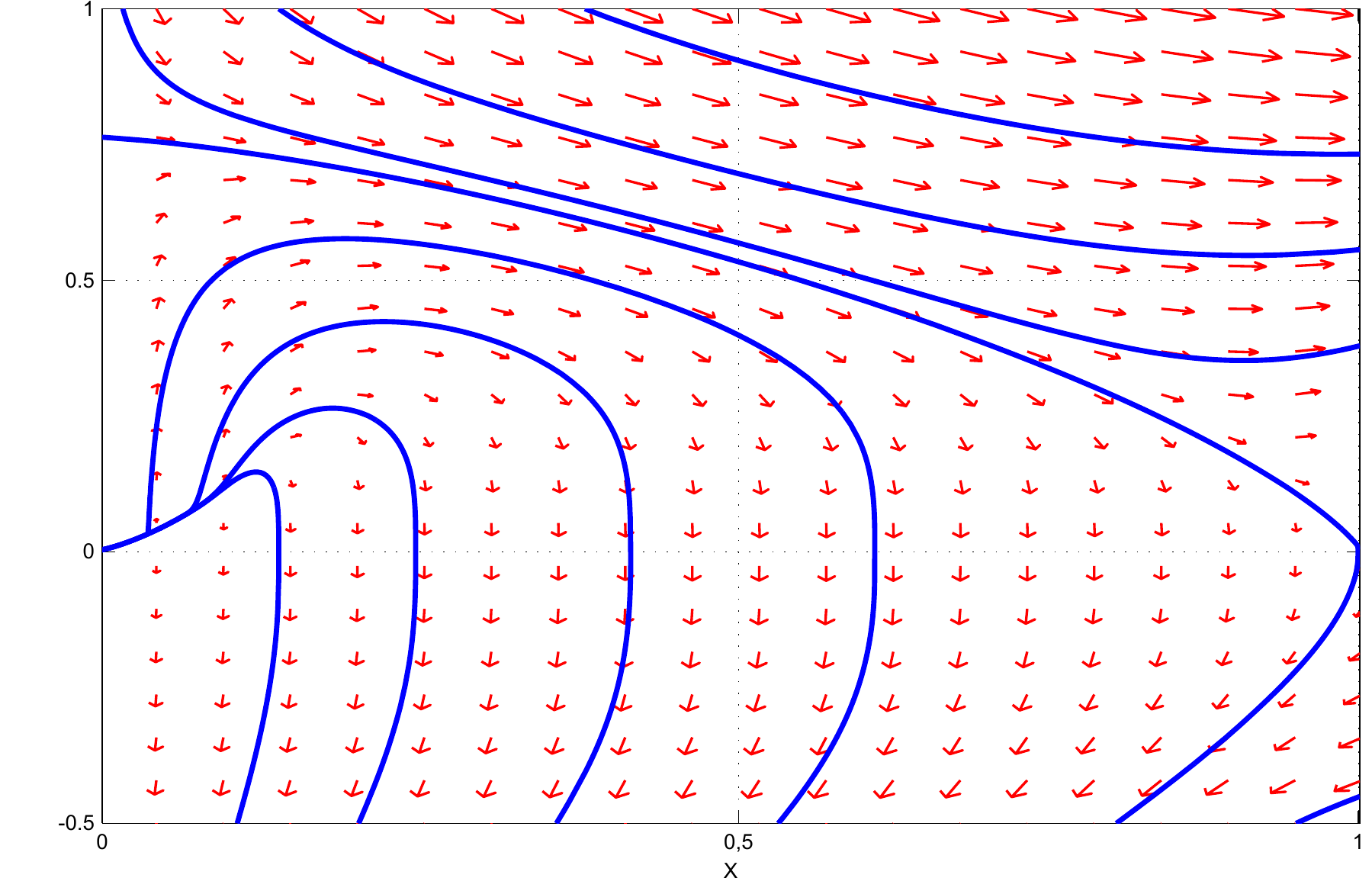}
  \includegraphics[scale=0.3]{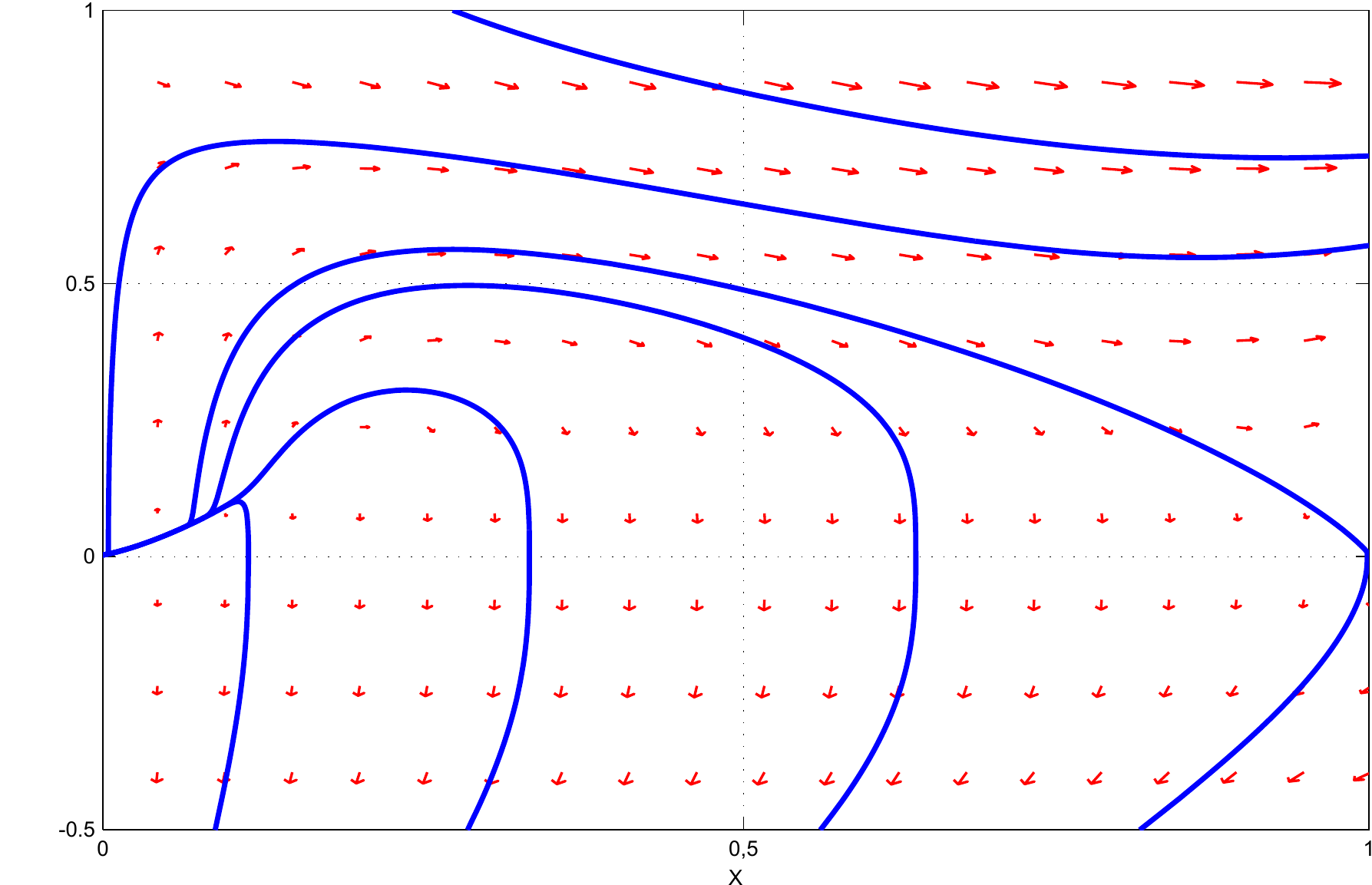}
  \caption{Phase plane analysis for $\gamma > 0$: Cases $c < c_{\ast}$, $c = c_{\ast}$ and $c > c_{\ast}$}\label{fig:TRAJECTORIESSLOWDIFFUSIONCASEGGPOSITIVE}
\end{figure}
\subsection{Analysis of the Change Sign TWs (CS-TWs) of type 2} In \emph{Step4} of the previous proof, we showed the existence of a family of  CS-TWs of type 2 which link the points $R_{\infty}$ and $R_{-\infty}$ in the $(X,Z)$-plane. These particular TWs will play a really important role in the study of more general solutions (see Section \ref{SECTIONASYMPTOTICBEHAVIOURDIMENSION1} and \ref{SECTIONASYMPTOTICBEHAVIOUR2}). In particular, we want to analyze their behaviour in the ``real'' plane $(\xi,\varphi(\xi)) = (\xi,X(\xi))$ near their global positive maximum (for example in $\xi = 0$) and near the two change sign points $\xi_0 < 0 < \xi_1$ where $\varphi(\xi_0) = 0 = \varphi(\xi_1)$.
\paragraph{Behaviour near the maximum point.} The existence of a positive global maximum follows from the analysis performed in the $(X,Z)$-plane and it corresponds to the value of $\xi$ such that $Z(\xi) = 0$ (we will always assume $\xi = 0$). However, in the proof of Theorem \ref{NTHEOREMCONVERGENCEINNEROUTERSETS} (see Section \ref{SECTIONASYMPTOTICBEHAVIOUR2}), we will need more information on the CS-TWs near their maximum point when $0 < m < 1$ and $p > 2$ (see relation \eqref{eq:ESTIMATEONSECONDDERIVATIVENEARMAXIMUMPOINT}). So, with this choice of parameters, we differentiate with respect to the variable $\xi$ the first equation in \eqref{eq:SYSTEMNONSINGULARTWSXI}:
\begin{equation}\label{eq:SECONDDERIVATIVEOFPROFILEFORMULA}
\frac{d^2X}{d\xi^2} = \frac{1}{m}\bigg[\Big(1-\frac{\gamma}{p-1}\Big)X^{-\gamma/(p-1)}\frac{dX}{d\xi}Z + X^{1-\gamma/(p-1)}\frac{dZ}{d\xi} \bigg].
\end{equation}
Then, substituting the expressions in \eqref{eq:SYSTEMNONSINGULARTWSXI} and taking $Z$ small, it is not difficult to see that
\[
\frac{d^2X}{d\xi^2} \sim -\frac{X^{-\frac{\gamma}{p-1}}f(X)}{m(p-1)} |Z|^{2-p}, \quad \text{for } Z \sim 0.
\]
In particular, for all $0 < X < 1$ fixed, it is straightforward to deduce the relation we are interested in:
\begin{equation}\label{eq:ESTIMATEONSECONDDERIVATIVENEARMAXIMUMPOINT}
\bigg|\frac{dX}{d\xi}\bigg|^{p-2}\bigg|\frac{d^2X}{d\xi^2}\bigg| \sim  \frac{m^{2-p}}{m(p-1)}X^{p-2-\gamma}f(X) > 0, \quad \text{for } Z \sim 0,
\end{equation}
near the maximum point of $\varphi = \varphi(\xi)$.
\paragraph{Behaviour near the change sign points.} For what concerns the profile's behaviour near zero, we can proceed formally considering the equation of the trajectories \eqref{eq:EQUATIONTRAJECTORIESTWS} and observe that if $X \sim 0$ and $Z \sim \infty$ we have $cZ - Z^p - f_{m,p}(X) \sim - Z^p$ and so
\[
\frac{dZ}{dX} \sim - \frac{1}{(p-1)}\frac{Z}{X} \qquad \text{i.e.} \qquad Z \sim a'X^{-\frac{1}{p-1}}, \quad \text{for } X \sim 0
\]
for some $a' > 0$. Now, from the first differential equation in \eqref{eq:SYSTEMNONSINGULARTWSXI}, we get
\begin{equation}\label{eq:ASYMPTOTICBEHAVIOURDERIVATIVEXZERO}
\frac{dX}{d\xi} \sim a X^{1-m}, \quad\text{for } X \sim 0
\end{equation}
where we set $a = a'/m$ and integrating it between $X = 0$ and $X = \varphi$ by separation of variables, we get that the profile $X = \varphi$ satisfies
\begin{equation}\label{eq:FREEBOUNDARYCONDITIONFAILED}
\varphi (\xi) \sim ma(\xi - \xi_0)^{1/m} \quad \text{for } \xi \sim \xi_0
\end{equation}
which not only explains us that the CS-TWs get to the level zero in finite time, but also tells us that they cannot be weak solutions of problem \eqref{eq:REACTIONDIFFUSIONEQUATIONPLAPLACIAN} since they violate the Darcy law: $\varphi(\xi) \sim (\xi - \xi_0)^{\gamma/(p-1)}$ near the \emph{free boundary} (see \cite{EstVaz:art}, or \cite{V2:book} for the Porous Medium equation). For this reason they cannot describe the asymptotic behaviour of more general solutions, but they turn out to be really useful when employed as ``barriers'' from below in the PDEs analysis. We ask the reader to note that, since the same procedure works when $X \sim 0$ and $Z \sim -\infty$, we obtain the existence of the second point $0 < \xi_1 < \infty$ such that $\varphi(\xi_1) = 0$ and, furthermore, the local analysis near $\xi_1$ is similar.

\noindent Finally, using that $Z \sim aX^{-1/(p-1)}$ for $X \sim 0$, we can again obtain information on the behaviour of the second derivative near the change sign point $\xi = \xi_0$. In particular, we get the relation
\begin{equation}\label{eq:SECONDESTIMATEONSECONDDERIVATIVENEARMAXIMUMPOINT}
X^{2m-1}\frac{d^2X}{d\xi^2} \sim \frac{a^2(p-2-\gamma)}{m^2(p-1)} \quad \text{for } X \sim 0 \text{ and } Z \sim \infty,
\end{equation}
near the ``free boundary point'' of $\varphi = \varphi(\xi)$. We will need this estimate in the proof of Theorem \ref{NTHEOREMCONVERGENCEINNEROUTERSETS} when in the case $0 < m < 1$ and $p > 2$.
\paragraph{Detailed derivation of relation (\ref{eq:ASYMPTOTICBEHAVIOURDERIVATIVEXZERO}).} Before, we deduced relation \eqref{eq:ASYMPTOTICBEHAVIOURDERIVATIVEXZERO} in a too formal way and we decided to end this section with a more complete proof. Let $Z = Z(X)$ be a branch of the trajectory of the CS-TW of type 2 and suppose $Z \geq 0$ (the case $Z \leq 0$ is similar).

\emph{Step1}. We start proving that $Z(X) \geq aX^{-1/(p-1)}$ for $X \sim 0$ and some $a > 0$. Since, $f_{m,p}(X) \geq 0$ we have
\[
\frac{dZ}{dX} \leq \frac{cZ - Z^p}{(p-1)XZ^{p-1}}
\]
which can be integrated by variable separation and gives $Z^{p-1}(X) \geq -X_0(c-Z_0)^{p-1}X^{-1}$, for all $0 < X \leq X_0 < 1$ ($X_0$ is the initial condition and $Z_0 = Z(X_0))$. Hence, taking $X_0 \sim 0$ and, consequently, $Z_0 \sim \infty$, we have that $-X_0(c-Z_0)^{p-1} \sim X_0Z_0^{p-1}$ and we deduce
\[
Z(X) \geq X_0^{\frac{1}{p-1}}Z_0X^{-\frac{1}{p-1}} \quad \text{for } X \sim 0.
\]

\emph{Step2}. Now, we show $Z(X) \leq aX^{-1/(p-1)}$ for $X \sim 0$. Using the fact that $f_{m,p}(X) \leq F_{m,p}$ and $Z \geq 0$, we get the differential inequality
\[
\frac{dZ}{dX} \geq -\frac{Z^p + F_{m,p}}{(p-1)XZ^{p-1}}.
\]
Proceeding as before, it is straightforward to deduce $Z^p \leq X_0^{p/(p-1)}(Z_0 + F) X^{-p/(p-1)}$, where $X_0$ and $Z_0$ are taken as before. Thus, since $X_0^{p/(p-1)}(Z_0 + F) \sim X_0^{p/(p-1)}Z_0$ for $X_0 \sim 0$, we obtain
\[
Z(X) \leq X_0^{\frac{1}{p-1}}Z_0X^{-\frac{1}{p-1}} \quad \text{for } X \sim 0,
\]
which allows us to conclude $Z(X) \sim aX^{-\frac{1}{p-1}}$ for $X \sim 0$ and $a = X_0^{\frac{1}{p-1}}Z_0$ and, consequently, \eqref{eq:ASYMPTOTICBEHAVIOURDERIVATIVEXZERO}.

\subsection{\texorpdfstring{\boldmath}{}Continuity of the function \texorpdfstring{$c_{\ast}$ for $\gamma > 0$}{gamma}}
We end this section studying the continuity of the function $c_{\ast} = c_{\ast}(m,p)$ with respect the parameters $m > 0$ and $p > 1$, when $\gamma > 0$ (we prove the first step of Theorem \ref{THMCONTINUITYOFCASTINCLOSUREOFH}). We will see that the continuity of the critical speed strongly depends on the stability of the orbit ``entering'' in the point $S = (0,1)$ (recall that we proved its existence and its uniqueness for all $c > 0$ in \emph{Step2} of Theorem \ref{THEOREMEXISTENCEOFTWS}). Before proceeding, we need to introduce the following notations:

\noindent $\bullet$ $Z_j = Z_j(X)$ stands for the analytic representation of the trajectory $T_{c_{\ast}(m_j,p_j)}$ (as a function of $X$) for the values $m_j > 0$ and $p_j > 1$ such that $\gamma_j = m_j(p_j-1)-1 > 0$ and $j=0,1$.

\noindent $\bullet$ $A_j = A_j(X)$ will indicate the analytic representation of the trajectories ``above'' $T_{c_{\ast}(m_j,p_j)}$, $j =0,1$.

\noindent $\bullet$ $B_j = B_j(X)$ will indicate the analytic representation of the trajectories ``below'' $T_{c_{\ast}(m_j,p_j)}$, $j =0,1$.

\noindent The following lemma proves that the orbit $T_{c_{\ast}}$ is continuous with respect to the parameters $m>0$ and $p>1$.
\begin{lem}\label{STABILITYOFTRAJECTORIESGGPOS}
Let $c = c_{\ast}$. Then the orbit $T_{c_{\ast}}$ linking $R_{c_{\ast}}= (0,C_{\ast})$ and $S=(1,0)$ is continuous with respect to the parameters $m>0$ and $p>1$ (with $\gamma > 0$) uniformly on $[0,1]$.

\noindent This means that for all $m_0 > 0$ and $p_0 > 1$ with $\gamma_0 > 0$, for all $\varepsilon > 0$ there exists $\delta > 0$ such that
\[
|Z_0(X) - Z_1(X)| \leq \varepsilon \quad \text{for all } |m_0 - m_1| + |p_0 - p_1| \leq \delta \text{ with } \gamma_1 > 0
\]
and for all $0 \leq X \leq 1$.
\end{lem}
\emph{Proof.} Fix $m_0 > 0$ and $p_0 > 1$ with $\gamma_0 > 0$ and note that the proof is trivial if $X = 1$.

\emph{Step1.} First of all, we show that for all $\varepsilon > 0$ and for all $0 < \overline{X} < 1$, there exists $\delta > 0$ such that
\[
|Z_0(\overline{X}) - Z_1(\overline{X})| \leq \varepsilon \quad \text{for all } |m_0 - m_1| + |p_0 - p_1| \leq \delta \text{ with } \gamma_1 > 0.
\]
So, fix $0 < \overline{X} < 1$ and $\varepsilon > 0$. We consider the trajectories $A_0 = A_0(X)$ and $B_0 = B_0(X)$ with
\begin{equation}\label{eq:SHOOTINGPOINTSCONTCAST}
A_0(\overline{X}) = Z_0(\overline{X}) + \varepsilon \quad \text{ and  } \quad B_0(\overline{X}) = Z_0(\overline{X})-\varepsilon,
\end{equation}
where $Z_0 = Z_0(X)$, as we explained before, is the analytic expression for the trajectory $T_{c_{\ast}}$ with parameters $m_0$ and $p_0$. Since we proved that $T_{c_{\ast}}$ is ``repulsive'' near $S$ (see \emph{Step2} of Theorem \ref{THEOREMEXISTENCEOFTWS}), we have that $A_0(\cdot)$ has to cross the line $Z = 1$ in some point with positive height while $B_0(\cdot)$ crosses the $X$-axis in a point with first coordinate in the interval $(\overline{X},1)$. Hence, we can apply the continuity of the trajectories with respect to the parameters $m$ and $p$ outside the critical points and deduce the existence of $\delta > 0$ such that for all $m_1$ and $p_1$ satisfying $|m_0 - m_1| + |p_0 - p_1| \leq \delta$, the trajectories $A_1 = A_1(X)$ and $B_1 = B_1(X)$ with
\begin{equation}\label{eq:SHOOTINGPOINTSCONTCAST1}
A_1(\overline{X}) = Z_0(\overline{X}) + \varepsilon \quad \text{ and  } \quad B_1(\overline{X}) = Z_0(\overline{X})-\varepsilon
\end{equation}
satisfy $|A_0(X) - A_1(X)| \leq \varepsilon$ and $|B_0(X) - B_1(X)| \leq \varepsilon$ for all $\overline{X} \leq X \leq 1$. In particular, $A_1(\cdot)$ crosses the line $Z = 1$ in a point with positive height and $B_1(\cdot)$ has to cross the $X$-axis in a point with first coordinate in the interval $(\overline{X},1)$. Consequently, since $B_1(X) \leq Z_1(X) \leq A_1(X)$ for all $\overline{X} \leq X \leq 1$, we deduce that $|Z_0(\overline{X}) - Z_1(\overline{X})| \leq \varepsilon$.

\emph{Step2.} Now, in order to show that
\[
|Z_0(X) - Z_1(X)| \leq \varepsilon \quad \text{ for all } \overline{X} \leq X \leq 1,
\]
we suppose by contradiction that there exists a point $\overline{X} < \overline{X}' < 1$ such that $|Z_0(\overline{X}') - Z_1(\overline{X}')| > \varepsilon$. Without loss generality, we can suppose $Z_0(\overline{X}') > Z_1(\overline{X}') + \varepsilon$. Then, we can repeat the procedure carried out before by taking $A_0(\cdot)$ and $B_0(\cdot)$ satisfying \eqref{eq:SHOOTINGPOINTSCONTCAST} with $\overline{X} = \overline{X}'$. Hence, the continuity of the trajectories with respect to the parameters $m$ and $p$ (outside the critical points) assures us the existence of $A_1(\cdot)$ and $B_1(\cdot)$ satisfying \eqref{eq:SHOOTINGPOINTSCONTCAST1} with $\overline{X} = \overline{X}'$. Hence, since $B_0(\cdot)$ crosses the $X$-axis in point with first coordinate in the interval $(\overline{X}',1)$ and $B_1(\cdot)$ has to behave similarly (by continuity), we have that the trajectory described by $B_1(\cdot)$ and $Z_1(\cdot)$ have to intersect, contradicting the uniqueness of the solutions.

\noindent We ask to the reader to note that, at this point, we have showed the continuity of the trajectory $T_{c_{\ast}}$ with respect to the parameters $m$ and $p$ uniformly in the interval $[\overline{X},1]$, where $0 < \overline{X} < 1$ is arbitrary.

\emph{Step3.} Finally, to conclude the proof, it is sufficient to check the continuity in $X = 0$, i.e., we have to prove that for all $m_0 > 0$ and $p_0 > 1$ such that $\gamma_0 > 0$, for all $\varepsilon > 0$ there exists $\delta > 0$ such that
\[
|Z_0(0) - Z_1(0)| \leq \varepsilon \quad \text{for all } |m_0 - m_1| + |p_0 - p_1| \leq \delta \text{ with } \gamma_1 > 0.
\]
Arguing by contradiction again, we suppose that there exist $\varepsilon >0$ such that for all $\delta > 0$, we can find $m_1$ and $p_1$ with $|m_0 - m_1| + |p_0 - p_1| \leq \delta$ and $\gamma_1 > 0$ such that it holds $|Z_0(0) - Z_1(0)| > \varepsilon$.

\noindent Hence, since the trajectories are continuous with respect to the variable $X$, we deduce the existence of a small $0 < \overline{X} < 1$ such that $|Z_0(\overline{X}) - Z_1(\overline{X})| > \varepsilon/2$ and so, thanks to the result from the \emph{Step1} and \emph{Step2}, we can take $\delta > 0$ small so that $|Z_0(\overline{X}) - Z_1(\overline{X})| \leq \varepsilon/2$, obtaining the desired contradiction. $\Box$

\bigskip

A direct consequence of the previous lemma is the continuity of the function $c_{\ast} = c_{\ast}(m,p)$ in the region $\mathcal{H}:=\{(m,p): \gamma = m(p-1)-1 > 0$, that we enunciate in the following corollary.
\begin{cor}\label{COROLLARYCONTINUITYCASTINNERSET}
The function $c_{\ast} = c_{\ast}(m,p)$ is continuous in the region $\mathcal{H}$, i.e., for all $m_0$ and $p_0$ with $\gamma_0 > 0$ and for all $\varepsilon > 0$, there exists $\delta > 0$ such that it holds $|c_{\ast}(m_0,p_0) - c_{\ast}(m,p)| \leq \varepsilon$, for all $m$ and $p$ with $\gamma >0$ and satisfying $|m-m_0|+|p-p_0| \leq \delta$.
\end{cor}
%
%
%
%
%
%
%
%
%
%
%
%

\section{\texorpdfstring{\boldmath}{}Existence of Travelling Wave solutions: case \texorpdfstring{$\gamma = 0$}{gamma}}\label{CLASSIFICATIONEXISTENCETW}
In this section, we address to problem of the existence of TW solutions for equation \eqref{eq:REACDIFFPLAPLACIANTWS1} in the ``pseudo-linear'' (limit) case $\gamma = 0$, i.e., $m(p-1) = 1$. Note that the relation $m(p-1) = 1$ generalizes the linear setting $m=1$ and $p=2$ but the diffusion operator is still nonlinear. For this reason, the Liapunov's linearization method cannot be applied and so we are forced to apply the same techniques used in the proof of Theorem \ref{THEOREMEXISTENCEOFTWS}. However, as we will see in a moment, there are important differences respect to the case $\gamma > 0$ that have to be underlined: in particular, it turns out that the admissible TWs are always positive while the TWs with \emph{free boundary} disappear (see Theorem \ref{THEOREMEXISTENCEOFTWS1}).
\paragraph{Proof of Theorem \ref{THEOREMEXISTENCEOFTWS1}.} We proceed as we did at the beginning of the previous proof, considering the equation of the profile \eqref{eq:PROFILEEQUATIONTWS} and performing the change of variable \eqref{eq:NONSTANDARDCHANGEOFVARIABLESTWS} which, in this case, can be written as
\begin{equation}\label{eq:NONSTANDARDCHANGEOFVARIABLE2}
X = \varphi \qquad \text{ and } \qquad Z = m\varphi^{-1}\varphi' = mX^{-1}X'.
\end{equation}
Hence, we get the system
\begin{equation}\label{eq:SYSTEMNONSINGULARTWSXI1}
\frac{dX}{d\xi} = (1/m)XZ, \quad\quad |Z|^{p-2}\, \frac{dZ}{d\xi} = cZ - |Z|^p - mX^{-1}f(X),
\end{equation}
which can be re-written (after the re-parametrization $d\xi = |Z|^{p-2}d\tau$) as the non-singular system
\begin{equation}\label{eq:SYSTEMNONSINGULARTWSXI2}
\frac{dX}{d\tau} = (p-1)X|Z|^{p-2}Z, \quad\quad \frac{dZ}{d\tau} = cZ - |Z|^p - F(X),
\end{equation}
where $F(X) := mX^{-1}f(X)$ (note that $F(0) = mf'(0)$, $F(1) = 0$, $F(X) \geq 0$ and $F'(X) \leq 0$ for all $0 \leq X \leq 1$).

\noindent System \eqref{eq:SYSTEMNONSINGULARTWSXI2} possesses the critical point $S = (1,0)$ for all $c \geq 0$. Now, define
\[
c_{0\ast}(m,p) := p(m^2f'(0))^{\frac{1}{mp}},
\]
as in \eqref{eq:DEFINITIONOFCASTPSEUDO}. Then it is no difficult to prove:

\noindent $\bullet$ If $c < c_{0\ast}(m,p)$, then there are no other critical points for system \eqref{eq:SYSTEMNONSINGULARTWSXI2}.

\noindent $\bullet$ If $c = c_{0\ast}(m,p)$, then system \eqref{eq:SYSTEMNONSINGULARTWSXI2} has another critical point $R_{\lambda_{\ast}} := (0,\lambda_{\ast})$, where we define for simplicity
\[
\lambda_{\ast} := (c_{0\ast}(m,p)/p)^m = (m^2f'(0))^{1/p},
\]
according to the definition of $c_{0\ast}(m,p)$.

\noindent $\bullet$ If $c > c_{0\ast}$, then system \eqref{eq:SYSTEMNONSINGULARTWSXI2} has also two critical points $R_{\lambda_i} = (0,\lambda_i)$, $i = 1,2$ where $0 < \lambda_1 < \lambda_{\ast} < \lambda_2 < c^m$.

\noindent This follows from the fact that if $X = 0$ then $dZ/d\tau = 0$ if and only if $cZ - |Z|^p - mf'(0) = 0$ and the number of solutions of this equation depends on the parameter $c > 0$ following our classification. Equivalently, one can write the equation of the trajectories
\begin{equation}\label{eq:EQUATIONOFTHETRAJECTORIESPSEUDO}
\frac{dZ}{dX} = \frac{cZ - |Z|^p - F(X)}{(p-1)X|Z|^{p-2}Z}
\end{equation}
and study the null isoclines imposing (exactly as we did at the beginning of the proof of Theorem \ref{THEOREMEXISTENCEOFTWS}):
\[
\max_{\widetilde{Z} \in [0,C]} \{c\widetilde{Z} - |\widetilde{Z}|^p\}= mf'(0).
\]
Solving the previous equation, we get the same value for $c_{0\ast}(m,p)$. In particular, for $c = c_{0\ast}(m,p)$ we have $c_{0\ast}\lambda_{\ast} - \lambda_{\ast}^p - mf'(0) = 0$ and $c_{0\ast} - p\lambda_{\ast}^{p-1} = 0$. When $c > c_{0\ast}(m,p)$, we have $c\lambda_i - \lambda_i^p - mf'(0) = 0$, for $i = 1,2$, while $c - p\lambda_1^{p-1} > 0$ and $c - p\lambda_2^{p-1} < 0$. Moreover, we can classify the null isoclines according to the ranges of the parameter $c > 0$: remembering the properties of the function $F(\cdot)$, it is not difficult to show that for the value $c = c_{0\ast}(m,p)$, we have a continuous isocline curve recalling a ``horizontal parabola'' with vertex in $R_{\lambda_{\ast}}$ and linking this point with $S$ and $(1,c_{0\ast}^m)$. Moreover, the trajectories are increasing in the area between the isocline and the line $X = 1$ whilst decreasing in the remaining part of the rectangle $[0,1]\times[0,c_{0\ast}^m]$. If $c < c_{0\ast}(m,p)$, the branch of the null isocline is again a ``horizontal parabola'', but, in this case, it does not have intersections with the axis $X = 0$. So, the trajectories have negative slope at the left of this curve while positive at the right. Finally, in the case $c > c_{0\ast}(m,p)$, the isoclines are composed by two branches: one in lower position linking the points $R_{\lambda_1}$ and $S$ and one in higher position linking $R_{\lambda_2}$ and $(1,c^m)$. In the area between these two branches the trajectories are increasing while in the rest of the rectangle $[0,1]\times[0,c^m]$ are decreasing. Note that for all $c \geq 0$, the slopes of the trajectories are negative when $Z \geq c^m$ while positive $Z \leq 0$.
\begin{figure}[!ht]
  \centering
  \includegraphics[scale=0.3]{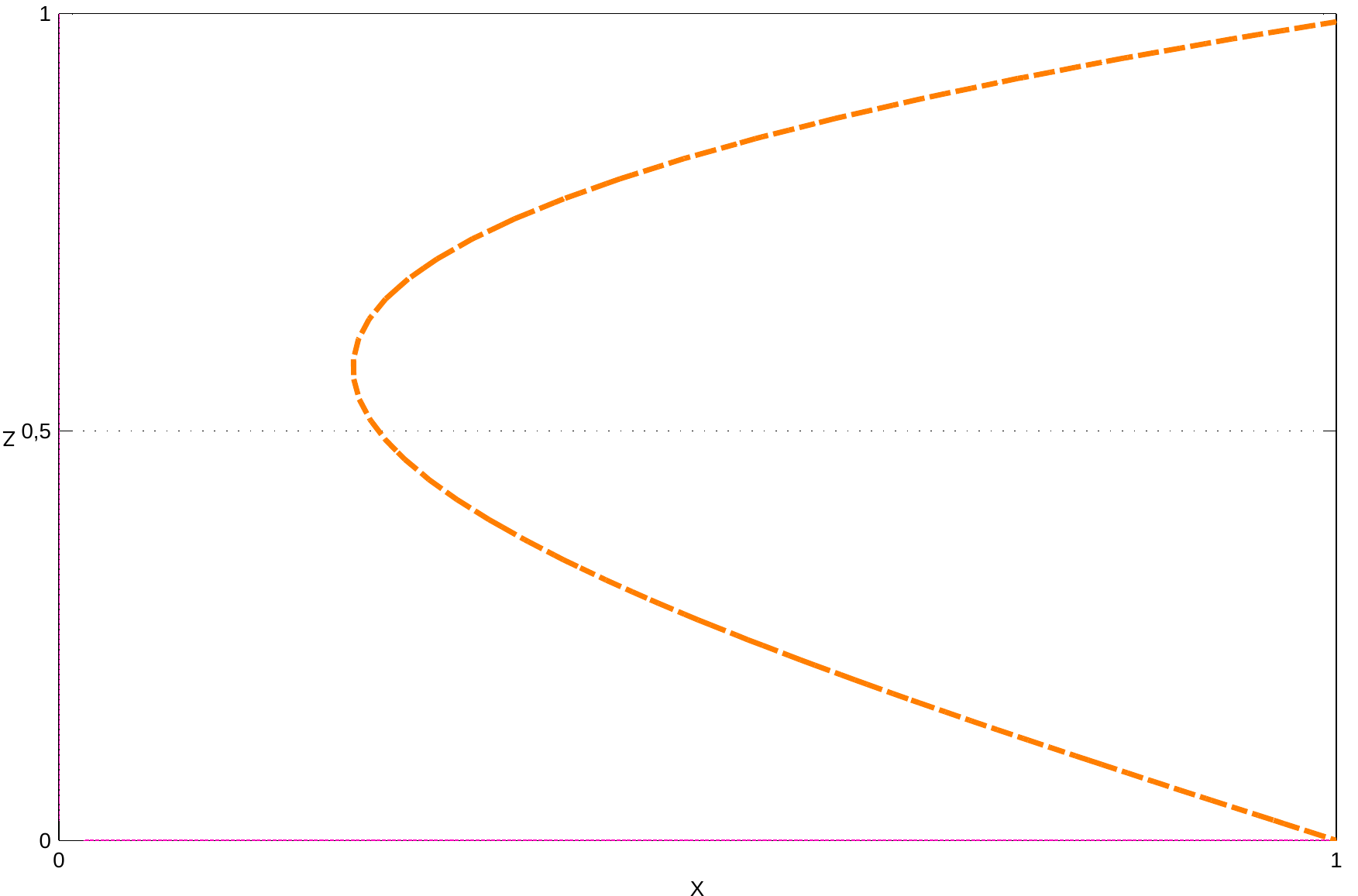}
  \includegraphics[scale=0.3]{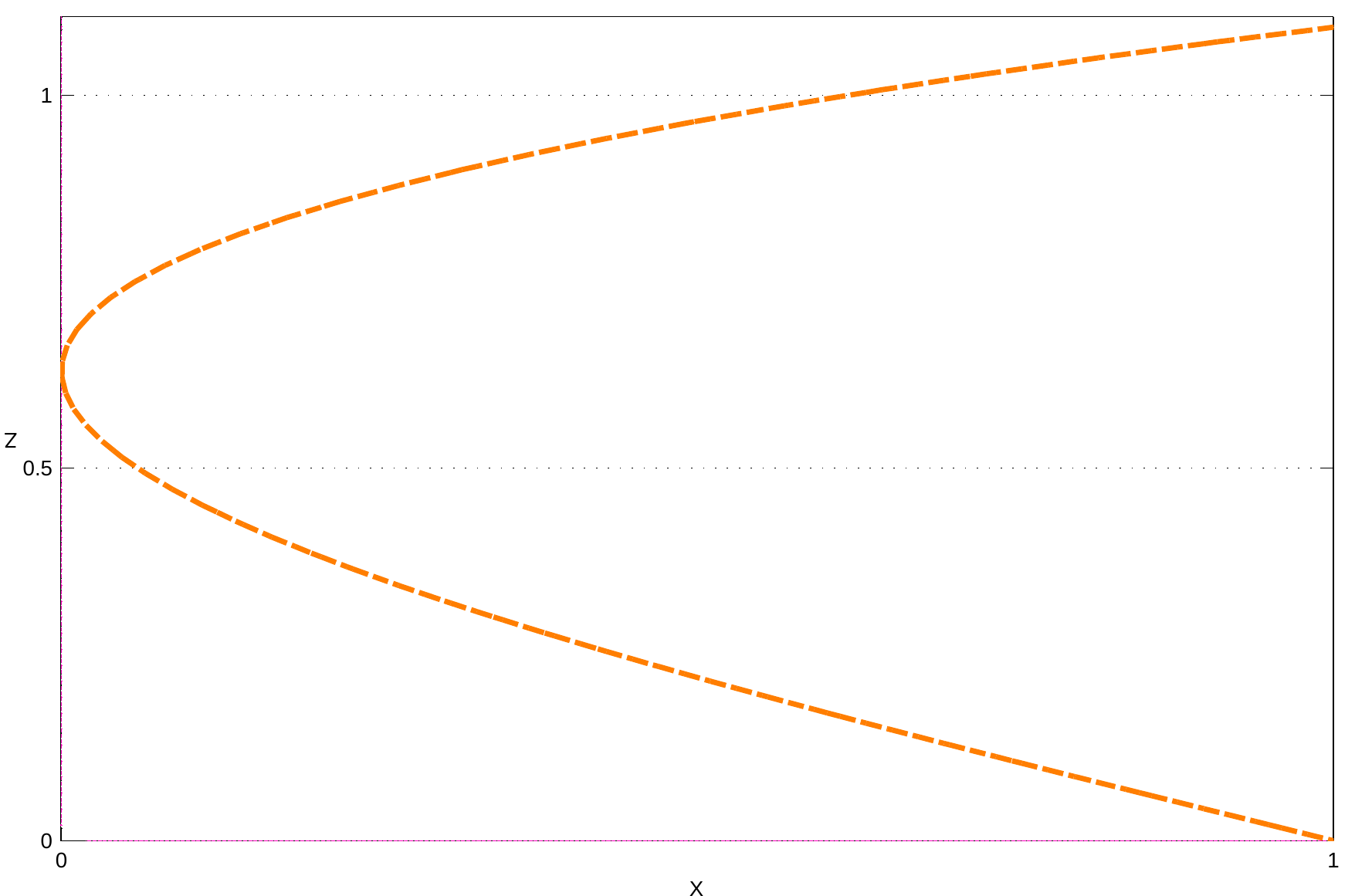}
  \includegraphics[scale=0.3]{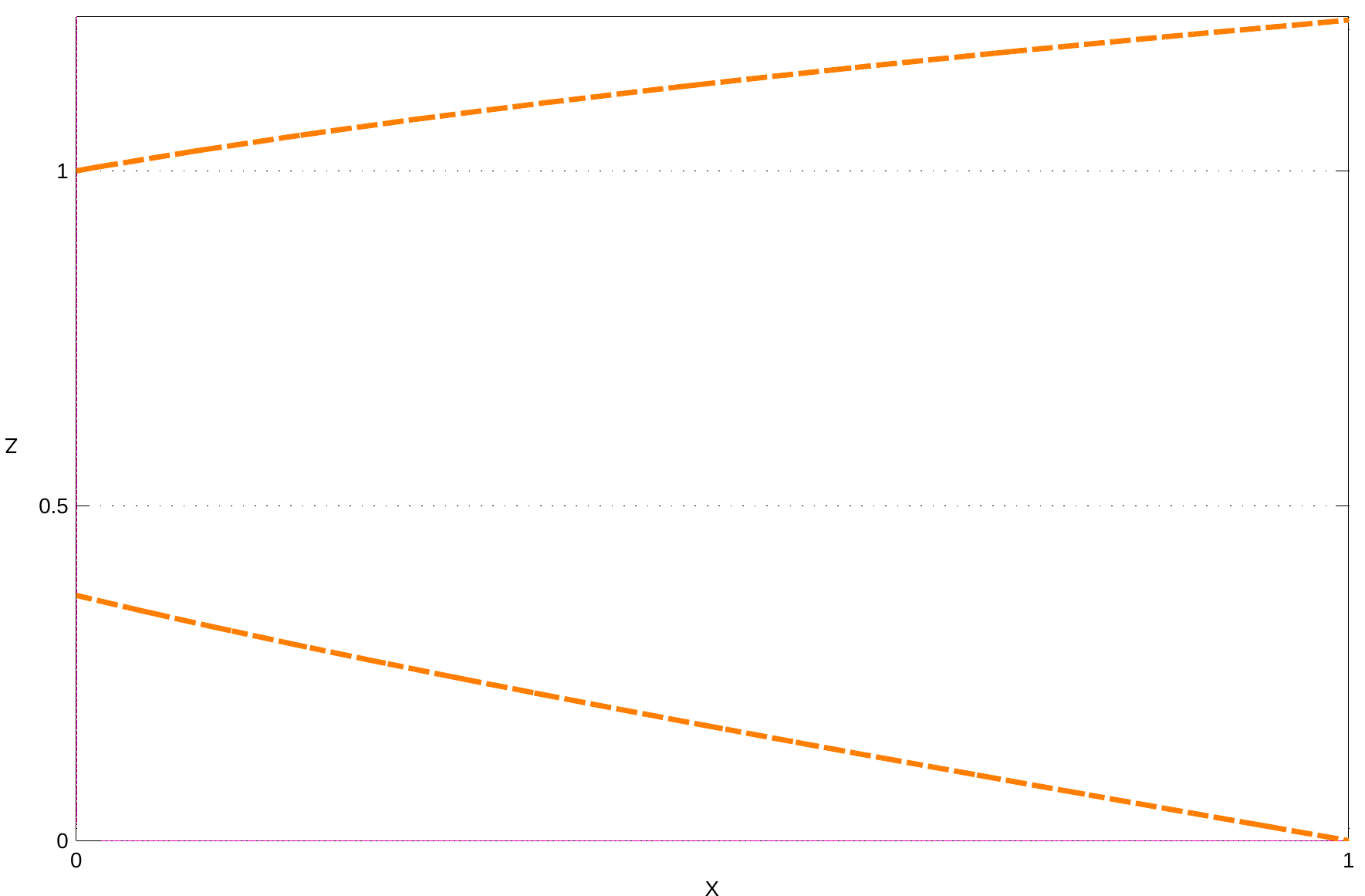}
  \caption{Null isoclines of system \eqref{eq:SYSTEMNONSINGULARTWSXI2}: Case $c < c_{0\ast}$, $c = c_{0\ast}$ and $c > c_{0\ast}$.}\label{fig:NULLISOCLINESPSEUDO}
\end{figure}

\noindent So it is clear, that the null isoclines in the case $\gamma = 0$ are quite different from the case $\gamma > 0$. However, it is no difficult to prove the existence and the uniqueness of a trajectory ``entering'' in the critical point $S$. This trajectory satisfies the same properties of the case $\gamma > 0$ (the fact that the same techniques used when $\gamma > 0$ apply here is not surprising since $f_{m,p}(X) \sim F(X) \sim mf(X)$ as $X \sim 1$). Now we are ready to prove our assertions.

\emph{Case} $c = c_{0\ast}(m,p)$. In this case, it sufficient to observe that the trajectory ``entering'' in $S$ has to cross the axis $X = 0$ in a point $(0,\overline{Z})$ with $0 < \overline{Z} \leq \lambda_{\ast}$. Nevertheless, proceeding as in \emph{Step4} and \emph{Step5} of the proof of Theorem \ref{THEOREMEXISTENCEOFTWS} it is simple to see that the only possibility is $\overline{Z} = \lambda_{\ast}$. Hence we have proved the existence of a connection $R_{\lambda_{\ast}} \leftrightarrow S$. To show that this admissible TW is positive we consider the first differential equation in \eqref{eq:SYSTEMNONSINGULARTWSXI1} and we integrate between $X = 0$ and $X = 1$ the equivalent differential relation $d\xi = (XZ)^{-1}dX$. It is straightforward to see that the correspondent integral is divergent both near $X = 0$ and $X = 1$, which means that the TW is positive.

\noindent Now we would like to describe the exact shape of this TW (with speed $c_{0\ast}(m,p) := p(m^2f'(0))^{\frac{1}{mp}}$) as $\xi \sim -\infty$. In \emph{Step5} of Theorem \ref{THEOREMEXISTENCEOFTWS}, in the case $\gamma > 0$, we have given an analytic representation of the (positive) TWs corresponding to the value $c > c_{\ast}(m,p)$. We found that the profile is a simple exponential for $\xi \sim -\infty$, see formula \eqref{eq:ASYMBEHASUPERTWGGPOS}. This has been possible since we have been able to describe the asymptotic behaviour of the trajectories in the $(X,Z)$-plane near $X=0$. The case $\gamma = 0$ is more complicated and we devote Section \ref{Appendix1} to the detailed analysis. Here, we simply report the asymptotic behaviour of our TW $\varphi = \varphi(\xi)$ with critical speed $c_{0\ast}(m,p)$:
\begin{equation}\label{eq:ASYMPTOTICTALESTWPSEUDOCAST}
\varphi(\xi) \sim a_0 |\xi|^{\frac{2}{p}}e^{\frac{\lambda_{\ast}}{m}\xi} = a_0 |\xi|^{\frac{2}{p}}\exp\Big(m^{\frac{2-p}{p}}f'(0)^{\frac{1}{p}}\xi\Big) \quad \text{for }\; \xi \sim -\infty,
\end{equation}
where $\lambda_{\ast} := (c_{0\ast}(m,p)/p)^m = (m^2f'(0))^{1/p}$ and $a_0 > 0$.

\noindent We conclude this paragraph with a brief description of the remaining trajectories. Below the positive TW, we have a family of trajectories linking $R_{\lambda_{\ast}}$ with $R_{-\infty} = (0,-\infty)$ while above it, there are a family of ``parabolas'' (exactly as in the case $\gamma > 0$). Between these two families, there are trajectories from the point $R_{\lambda_{\ast}}$ which cross the line $X = 1$.

\emph{Case} $c < c_{0\ast}(m,p)$. In this case there are not admissible TWs. Indeed, using again the same methods of \emph{Step4} and \emph{Step5} it is simple to show that the orbit ``entering'' in $S$ cannot touch the axis $X = 0$, i.e., it links $R_{\infty} = (0,\infty)$ and $S$ and so it is not an admissible TW. Below this trajectory there are a family of CS-TWs of type 2 and above it we have a family of ``parabolas''. We stress that, with the same techniques used for the case $\gamma >0$, the CS-TWs of type 2 satisfy $Z(X) \sim a' X^{-1/(p-1)} = a'X^{-m}$ for $X \sim 0$ and a suited positive constant $a'$. Hence, exactly as we did in Section \ref{SECTIONEXISTENCEOFTWS}, it is possible to show that the profile $X = \varphi$ reaches the level zero in to finite points $-\infty < \xi_0 < \xi_1 < \infty$. Indeed, exactly as we did before, using the first differential equation in \eqref{eq:SYSTEMNONSINGULARTWSXI1}, we get the estimate
\[
\frac{dX}{d\xi} \sim a X^{1-m} \quad \text{for } X \sim 0,
\]
from which we can estimate the (finite) times $\xi_0$ and $\xi_1$ (note that $a = a'/m$).
\begin{figure}[!ht]
  \centering
  \includegraphics[scale=0.3]{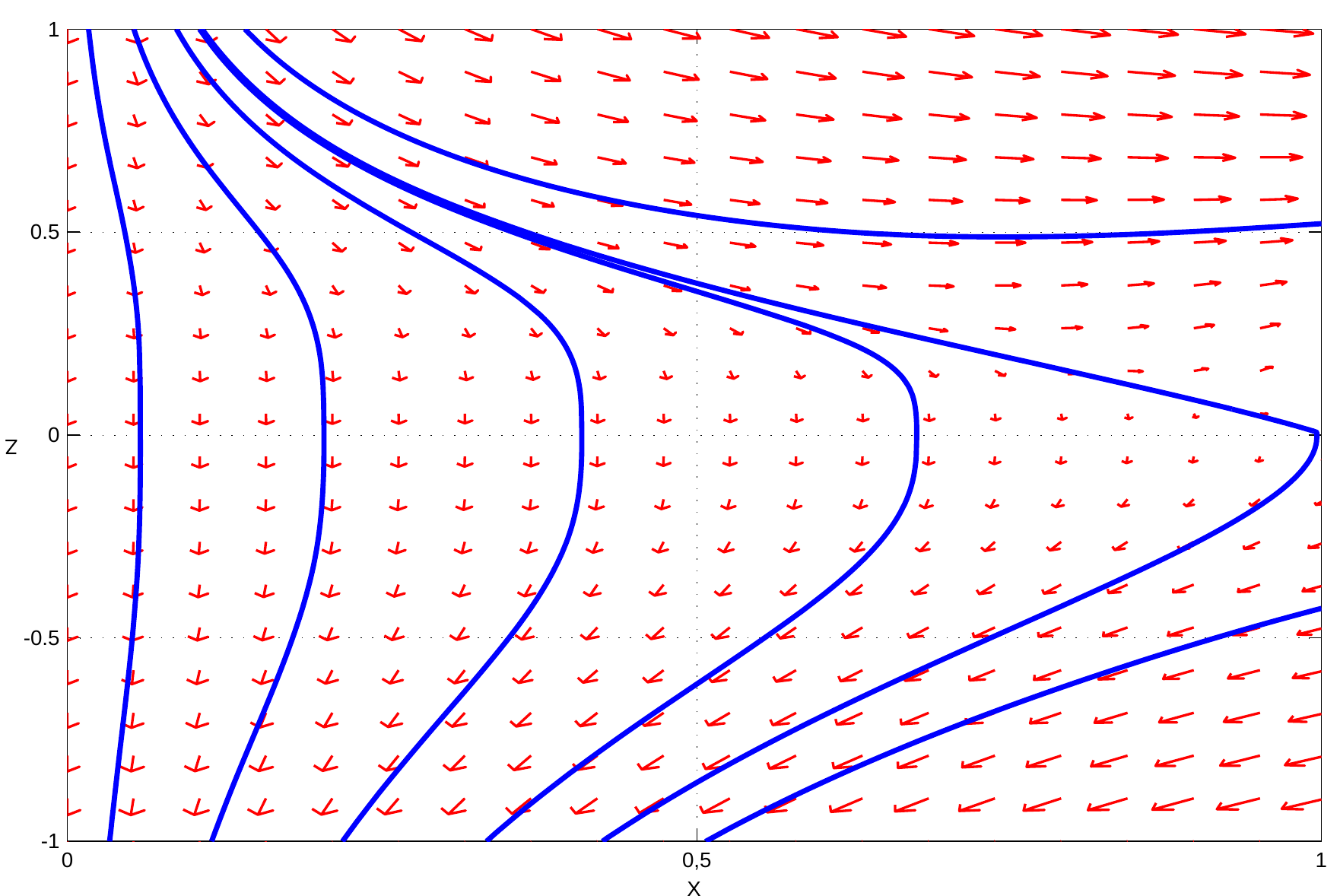}
  \includegraphics[scale=0.3]{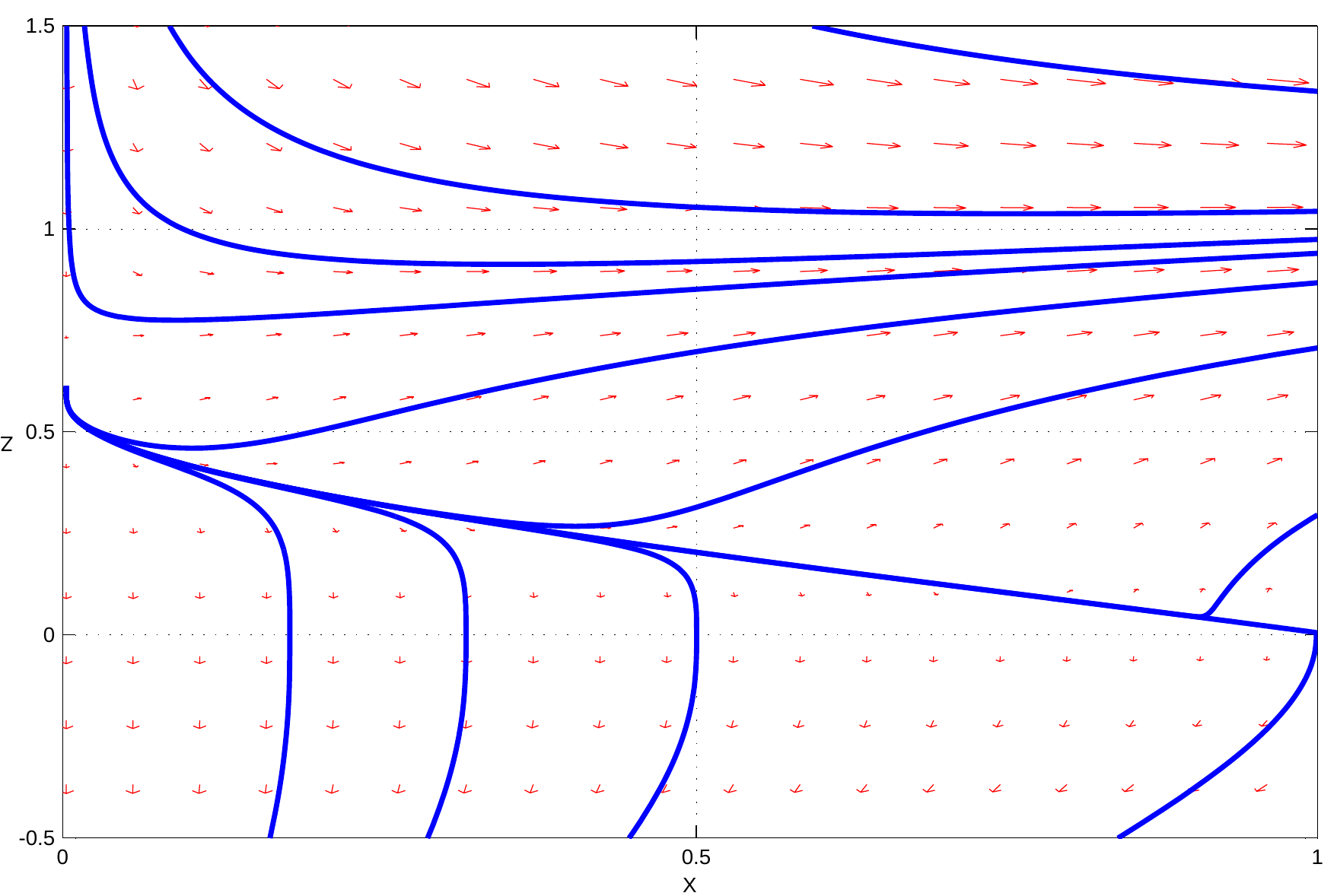}
  \includegraphics[scale=0.3]{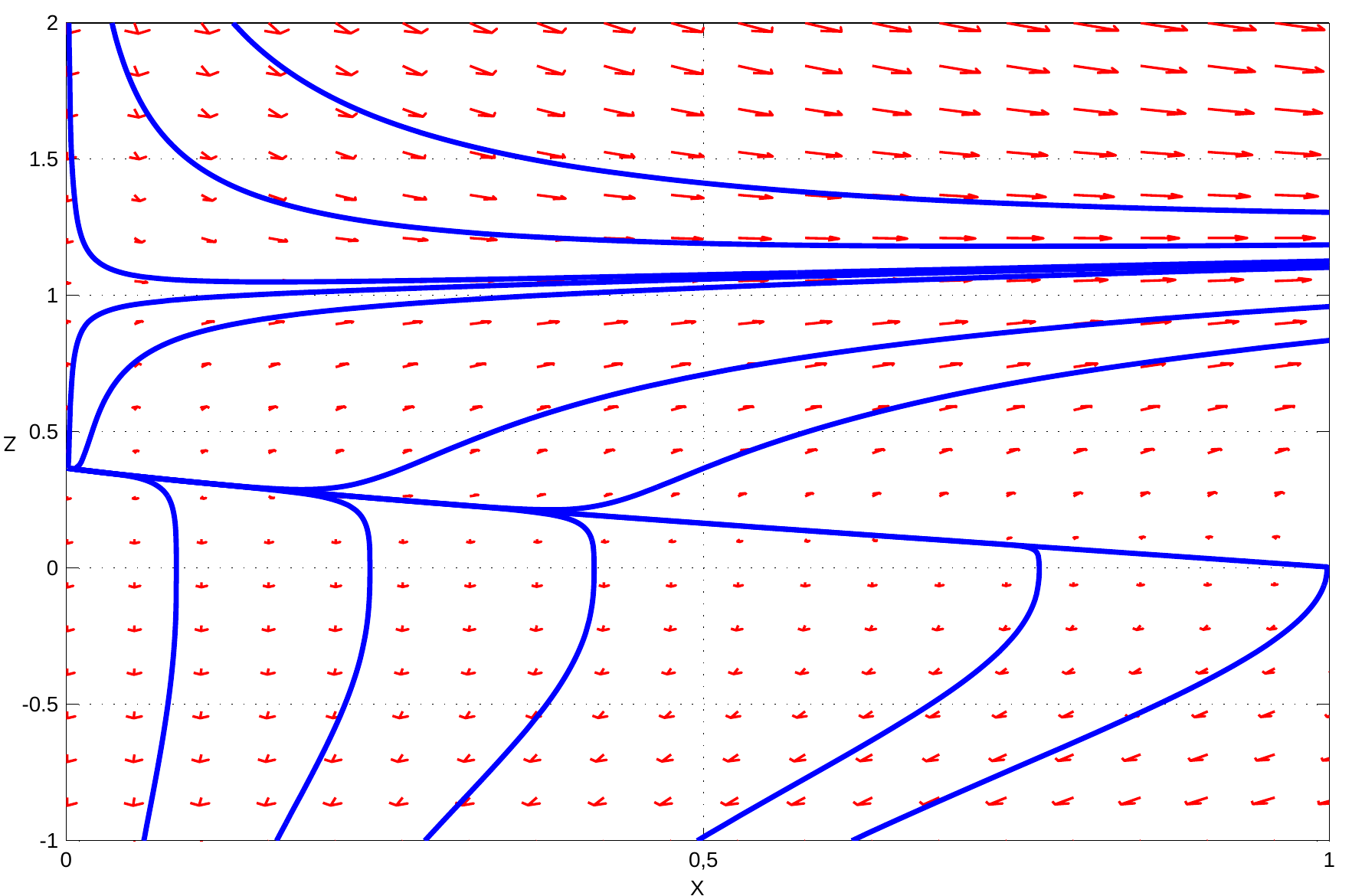}
  \caption{Phase plane analysis for $\gamma = 0$: Cases $c < c_{\ast}$, $c = c_{\ast}$ and $c > c_{\ast}$}\label{fig:TRAJECTORIESSLOWDIFFUSIONCASEGGZERO}
\end{figure}

\emph{Case} $c > c_{0\ast}(m,p)$. To show the existence of a connection $R_{\lambda_1} \leftrightarrow S$ and to prove that this TW is positive, it is sufficient to apply the methods used in the case $c = c_{0\ast}(m,p)$. We recall that for all $c > c_{0\ast}(m,p)$, $\lambda_1 = \lambda_1(c) < \lambda_{\ast}$, it solves $c\lambda_1 - \lambda_1^p - mf'(0) = 0$, and $c - p\lambda_1^{p-1} > 0$.

\noindent Now, by using the Lyapunov method, it is possible to linearize system \eqref{eq:SYSTEMNONSINGULARTWSXI2} around the critical point $R_{\lambda_1} = (0,\lambda_1)$, and it is straightforward to see that its Jacobian matrix (calculated in $(0,\lambda_1)$) has two positive eigenvalues $\nu_1 = (p-1)\lambda_1^{p-1}$ and $\nu_2 = c-p\lambda_1^{p-1}$. This means that $R_{\lambda_1}$ is a \emph{node unstable}. Hence, we deduce that the trajectory $Z = Z(X)$ from $R_{\lambda_1} = (0,\lambda_1)$ satisfies $Z(X) \sim \lambda_1 - \nu_1 X$, for $X \sim 0$ and we can re-write the first equation in \eqref{eq:SYSTEMNONSINGULARTWSXI1} as
\[
\frac{dX}{d\xi} \sim (1/m)X(\lambda_1 - \nu_1 X) \quad \text{for } X \sim 0.
\]
This is a first order \emph{logistic} type ODE, and so we easily obtain that the profile $X(\xi) = \varphi(x+ct)$ has the exponential decay
\begin{equation}\label{eq:ASYMBEHGGGEQ0SUPERTW}
X(\xi) \sim a_0e^{\frac{\lambda_1}{m}\xi} , \quad \text{for }\; \xi \sim -\infty,
\end{equation}
where $a_0 > 0$ is a fixed constant.

\noindent Below this connection we have a family of trajectories joining $R_{\lambda_1}$ and $R_{-\infty}$. Above it, there are trajectories from $R_{\lambda_2}$ and crossing the line $X=1$, one trajectory from $R_{\lambda_2}$ crossing the line $X=1$ and, finally, a family of ``parabolas'' as in the other cases. $\Box$
\paragraph{Asymptotic behaviour of $\boldsymbol{c_{0\ast}(m,p)}$.} We want to comment the asymptotic behaviour of the function $c_{0\ast}(m,p)$. Since we assume $m(p-1) = 1$, it is simple to re-write the critical speed as a function of $m > 0$ or $p > 1$:
\[
c_{0\ast}(m) = (1+m)m^{\frac{1-m}{1+m}}f'(0)^{\frac{1}{m+1}} \quad\text{ or } \quad c_{0\ast}(p) = p (p-1)^{-\frac{2(p-1)}{p}}f'(0)^{\frac{p-1}{p}}.
\]
Then we have:
\[
\lim_{m\to0}c_{0\ast}(m) = \lim_{p\to\infty}c_{0\ast}(p) = 0 \quad\text{ and } \quad \lim_{m\to\infty}c_{0\ast}(m) = \lim_{p\to1}c_{0\ast}(p) = 1.
\]
These limits allow us to guess if in the cases $m \to 0$, $p \to \infty$ and $m \to \infty$, $p \to 1$ such that $m(p-1) = 1$, there exist TWs or not. It seems natural to conjecture that the answer is negative in the first case, while admissible TWs could exist in the second case.
\begin{figure}[!ht]
  \centering
  \includegraphics[scale=0.4]{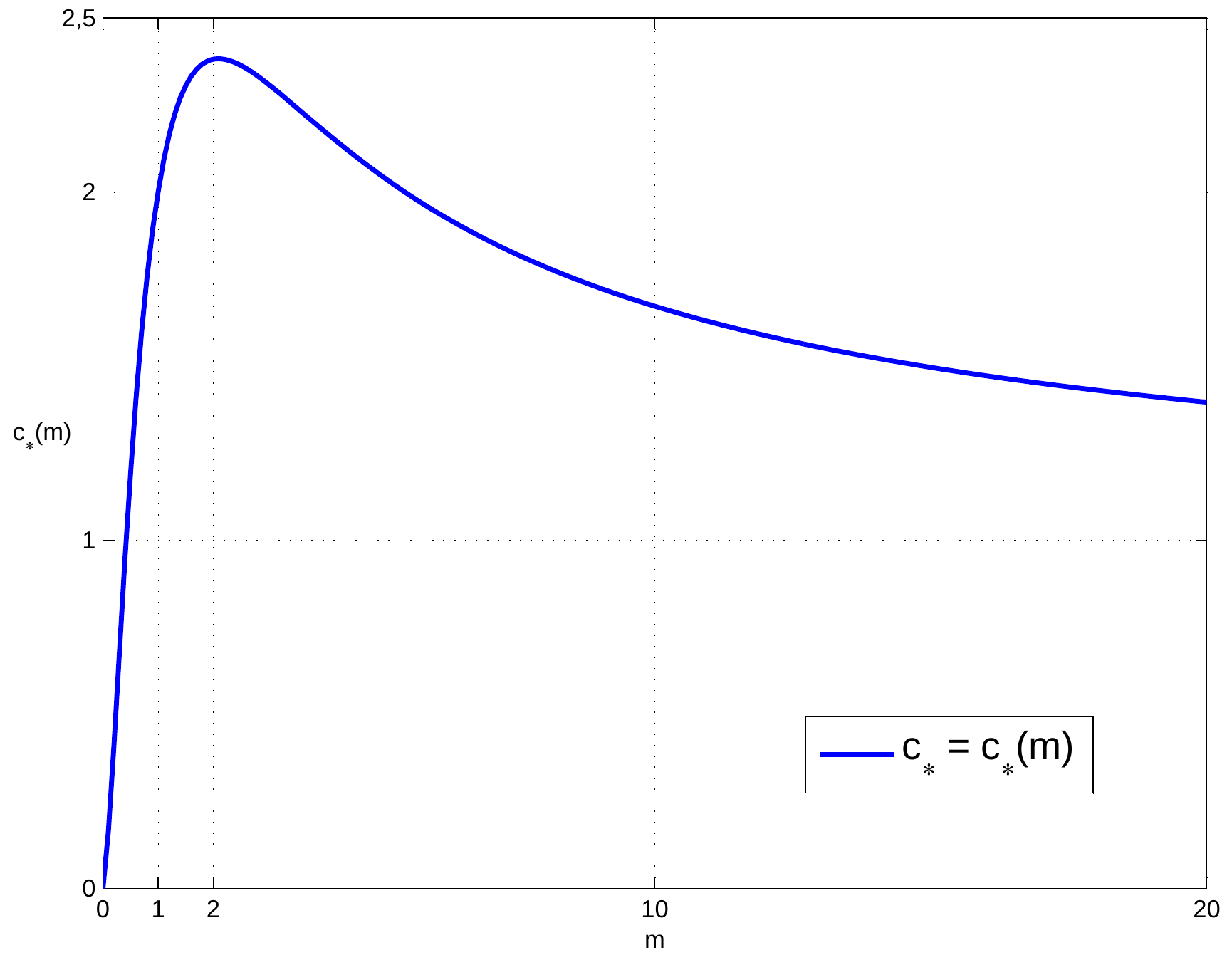} \quad
  \includegraphics[scale=0.4]{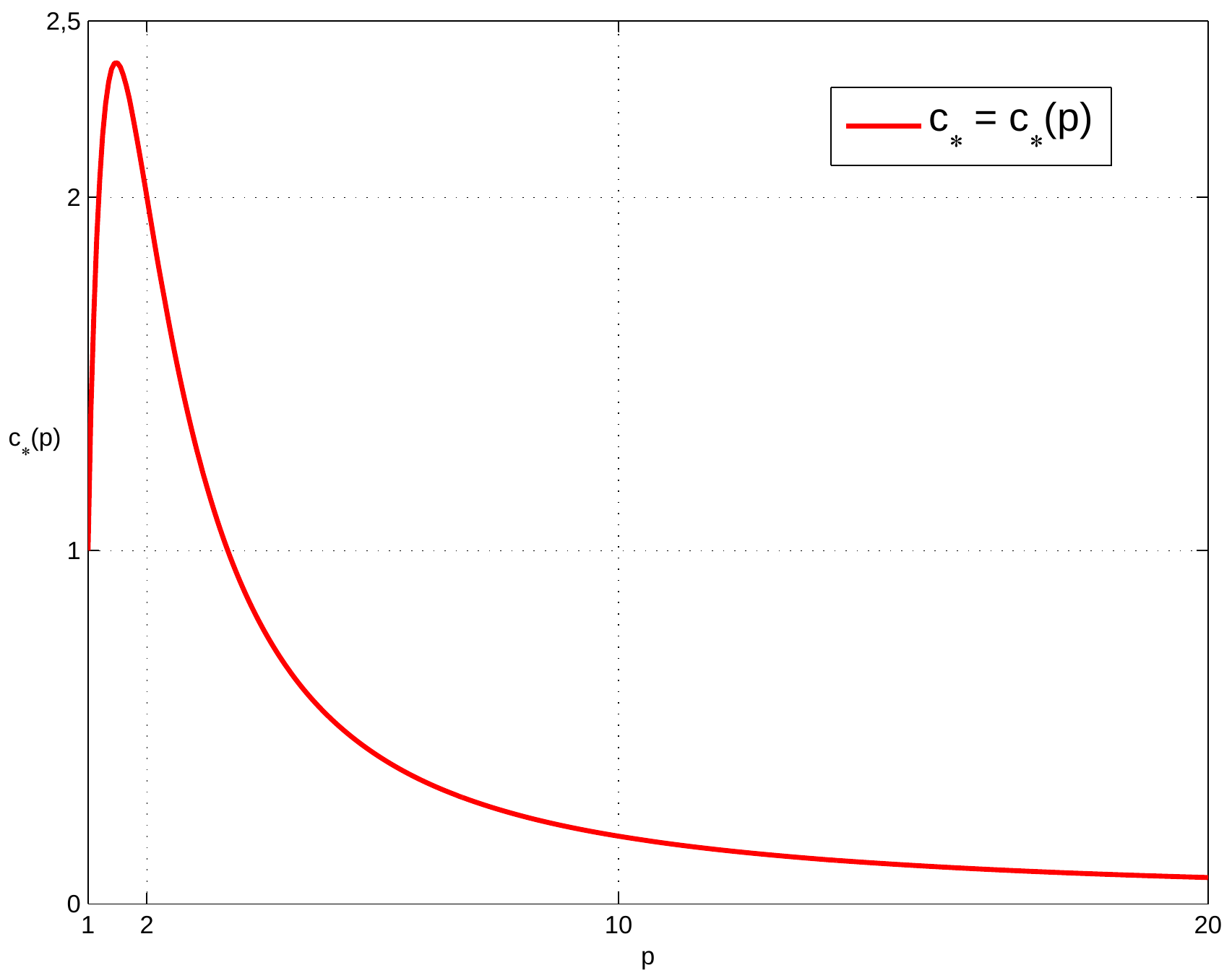}
  \caption{Pseudo-linear case: the graphics of $c_{0\ast}$ in dependence of $m > 0$ and $p > 1$.}\label{fig:CASTASFUNCTIONOFPMPSEUDO}
\end{figure}

\noindent Finally, it is not difficult to calculate the derivative of $c_{0\ast}(m)$ (or $c_{0\ast}(p)$) in the case $f(u) = u(1-u)$:
\[
c_{0\ast}'(m) = (m+1)^{-1}m^{-\frac{2m}{m+1}}(m + 1 - 2m\log m)
\]
and conclude that the maximum of the critical speed is assumed for a critical value of $m = m_{\ast}$, with $2 < m_{\ast} < 3$ (of course, we can repeat this procedure with the function $c_{0\ast}(p)$ finding a critical value $4/3 < p_{\ast} < 3/2$). This simple calculation assures the the linear case $m = 1$ and $p = 2$ is not critical for the function $c_{0\ast}(m,p)$ while the maximal speed of propagation is found choosing $m = m_{\ast}$ and/or $p = p_{\ast}$.
\paragraph{Remarks.} We end the discussion on the value $c_{0\ast}$ in the ``pseudo-linear'' case with two observations. First of all, since the formulas for $c_{\ast}(m)$ and $c_{\ast}(p)$ depend on the value of $f'(0)$, it is always possible to modify the reaction term $f(\cdot)$ with $\lambda f(\cdot)$ for some $\lambda > 0$ (depending on $m$ or $p$) and re-scale the space variable to have $c_{0\ast}(m) = c_{0\ast}(p) = 1$.

\noindent Secondly, we remind the reader to note that we obtained a value of $c_{0\ast}(m,p) = p(m^2f'(0))^{1/(mp)}$ which generalizes the ``classical one'' (with $m=1$ and $p=2$) found in \cite{Aro-Wein1:art, Aro-Wein2:art} and \cite{K-P-P:art}, $c_{\ast} = 2\sqrt{f'(0)}$. Moreover, it is interesting to see that our methods give an alternative and independent proof of the classical result stated in Theorem \ref{THMEXISTENCETWLINEARINTRO}.

\subsection{\texorpdfstring{\boldmath}{}Continuity of the function \texorpdfstring{$c_{\ast}$ for $\gamma \geq 0$}{gamma}} In these last paragraphs, we complete the proof of Theorem \ref{THMCONTINUITYOFCASTINCLOSUREOFH} showing the following lemma.
\begin{lem}\label{CONTINUITYOFCASTCOMPLETECLOSURELEM}
The function defined in \eqref{eq:CASTEXTENDEDDEFINITION} is continuous in $\overline{\mathcal{H}}$, i.e., for all $m_0$ and $p_0$ with $\gamma_0 \geq 0$ and for all $\varepsilon > 0$, there exists $\delta > 0$ such that it holds $|c_{0\ast}(m_0,p_0) - c_{\ast}(m,p)| \leq \varepsilon$, for all $m$ and $p$ with $\gamma \geq 0$ and satisfying $|m-m_0|+|p-p_0| \leq \delta$.
\end{lem}
\emph{Proof.} We divide the proof in some steps. In the first ones, we simply recall what we have showed early.

\emph{Step1.} First of all, we observe that, thanks to Corollary \ref{COROLLARYCONTINUITYCASTINNERSET} and since the function $c_{0\ast}(m,p)$ is continuous in the set $\{\gamma = 0\}$, it is sufficient to check the continuity along the boundary of the region $\mathcal{R}$. More precisely, we will show that for all $m_0$ and $p_0$ with $\gamma_0 = 0$, for all $\varepsilon > 0$, there exists $\delta > 0$ such that it holds $|c_{0\ast}(m_0,p_0) - c_{\ast}(m,p)| \leq \varepsilon$, for all $m$ and $p$ with $\gamma > 0$ and satisfying $|m-m_0|+|p-p_0| \leq \delta$.

\emph{Step2.} Now, we ask the reader to note that, with the same notations and techniques used in the proof of Lemma \ref{STABILITYOFTRAJECTORIESGGPOS}, it is possible to prove the continuity of the trajectory ``entering'' in $S = (1,0)$ with respect to the parameters $m$ and $p$ with $\gamma \geq 0$ uniformly in the variable $X$ in all sets $[\overline{X},1]$, where $\overline{X}$ is fixed in $(0,1]$. This fact can be easily checked since, as we explained in the proof of Theorem \ref{THEOREMEXISTENCEOFTWS1}, the local behaviour of the trajectories near the point $S$ is the same for $\gamma > 0$ and $\gamma = 0$.

\emph{Step3.} Before, proceeding with the proof we need to recall a last property. In the proof of Theorem \ref{THEOREMEXISTENCEOFTWS}, we showed that $c_{\ast}(m,p) < c_0(m,p)$ for all $m$ and $p$ with $\gamma > 0$, where
\[
c_0(m,p) := p\Bigg(\frac{F_{m,p}}{p-1}\Bigg)^{(p-1)/p},
\]
and $F_{m,p}$ is the maximum of the function $f_{m,p}(X) = mX^{\frac{\gamma}{p-1}-1}f(X)$. Since $f_{m,p}(X) \to mX^{-1}f(X)$ as $\gamma \to 0$ for all $0 < X \leq 1$ and the limit function is decreasing, we have that the maximum point of $f_{m,p}(\cdot)$ converges to zero as $\gamma \to 0$. Consequently, we obtain that $F_{m,p} \to mf'(0)$ as $\gamma \to 0$ and we deduce
\[
c_0(m,p) \to p(m^2f'(0))^{1/(m+1)} = c_{0\ast}(m,p)  \quad \text{as } \gamma \to 0.
\]
Now we have all the elements needed for completing the proof by proving the assertion stated in \emph{Step1}.

\emph{Step4.} Suppose by contradiction that there exists $m_0$ and $p_0$ with $\gamma_0 = 0$, $\varepsilon > 0$ and a sequence $(m_j,p_j) \to (m_0,p_0)$ as $j \to \infty$ with $\gamma_j > 0$, such that $|c_{0\ast}(m_0,p_0) - c_{\ast}(m_j,p_j)| \geq \varepsilon$ for all $j \in \NN$.

\noindent Hence, there exists a subsequence that we call again $c_{\ast}(m_j,p_j)$ converging to a value $c \not= c_{0\ast}(m_0,p_0)$. Note that, since we proved $c_{\ast}(m_j,p_j) < c_0(m_j,p_j) \to c_{0\ast}(m_0,p_0)$, it has to be $c < c_{0\ast}(m_0,p_0)$ (see \emph{Step3}).

\noindent Now, since $\gamma_j > 0$, we have that the trajectory $T_{c_{\ast}(m_j,p_j)}$ joins the point $S = (1,0)$ with the point $R_{c_{\ast}(m_j,p_j)} = (0,c_{\ast}(m_j,p_j)^{1/(p_j-1)})$ for all $j \in \NN$ and $R_{c_{\ast}(m_j,p_j)} \to R_c = (0, c^{1/(p_0-1)})$ as $j \to \infty$. In particular, it follows that the sequence of trajectories $T_{c_{\ast}(m_j,p_j)}$ has to be bounded. However, thanks to the analysis done in the proof of Theorem \ref{THEOREMEXISTENCEOFTWS1}, we know that if $c < c_{0\ast}(m_0,p_0)$, the trajectory ``entering'' in the point $S$ is unbounded and joins the previous point with $R_{\infty} = (0,\infty)$. Consequently, applying the continuity of the trajectory ``entering'' in $S$ with respect to the parameters $m$ and $p$ (stated in \emph{Step2}) we obtain the desired contradiction and we conclude the proof. $\Box$
%
%
%
%
%
%
%
%
%
%
%
%
%
%
\section{\texorpdfstring{\boldmath}{}Case \texorpdfstring{$\gamma > 0$}{gamma}. Existence of minimal expanding super-level sets}\label{SECTIONLINEAREXPANSIONSUPERLEVELSETS}
From this moment on we address the question of the asymptotic analysis of general solutions. More precisely, we devote five sections to the issue. In this section we prove Proposition \ref{LEMMANONDISCRETEVERSIONLEVELSETS} by appropriately combining a series of lemmas. As mentioned in the introduction, we study problem \eqref{eq:REACTIONDIFFUSIONEQUATIONPLAPLACIAN}, with reaction term $f(\cdot)$ satisfying \eqref{eq:ASSUMPTIONSONTHEREACTIONTERM} and with the particular choice of the initial datum:
\begin{equation}\label{eq:INITIALDATUMTESTLEVELSETS}
\widetilde{u}_0(x) :=
\begin{cases}
\begin{aligned}
\widetilde{\varepsilon} \qquad &\text{if } |x| \leq \widetilde{\varrho}_0 \\
0 \qquad &\text{if } |x| > \widetilde{\varrho}_0,
\end{aligned}
\end{cases}
\end{equation}
where $\widetilde{\varepsilon}$ and $\widetilde{\varrho}_0$ are positive real numbers. The choice of $\widetilde{u}_0(\cdot)$ in \eqref{eq:INITIALDATUMTESTLEVELSETS} is related to the finite propagation of the Barenblatt solutions in the case $\gamma > 0$ (see Section \ref{PRELIMINARIESINTRO}) and its usefulness will be clear in the next sections (see Section \ref{SECTIONASYMPTOTICBEHAVIOUR}).

\noindent We begin our study with an elementary lemma. Even though it has a quite simple proof, its meaning is important: it assures the existence of a (large) time such that the solution $u = u(x,t)$ does not extinguish in time, but, on the contrary, remains larger than a (small) positive level, on all compact sets of $\RR^N$.
%
%
%
%
\begin{lem}\label{EXPANDINGLEVELSETS2}
Let $m > 0$ and $p > 1$ such that $\gamma > 0$ and let $N \geq 1$. Then for all $0 < \widetilde{\varepsilon} < 1$, for all $\widetilde{\varrho}_0 > 0$ and for all $\widetilde{\varrho}_1 \geq \widetilde{\varrho}_0$, there exists $t_1 > 0$ and $n_1 \in \NN$ such that the solution $u(x,t)$ of problem \eqref{eq:REACTIONDIFFUSIONEQUATIONPLAPLACIAN} with initial datum \eqref{eq:INITIALDATUMTESTLEVELSETS} satisfies
\[
u(x,t_1) \geq \widetilde{\varepsilon}/n_1 \quad \text{in } \{|x| \leq \widetilde{\varrho}_1\}.
\]
\end{lem}
\emph{Proof.} Fix $0 < \widetilde{\varepsilon} < 1$ and $0 < \widetilde{\varrho}_0 \leq \widetilde{\varrho}_1$.

\noindent We start with constructing a Barenblatt solution with positive parameters $M_1$ and $\theta_1$ such that $B_{M_1}(x,\theta_1) \leq \widetilde{u}_0(x)$ in $\RR^N$. Since the profile of the Barenblatt solution is decreasing, we can choose $M_1, \theta_1 > 0$ such that $B_{M_1}(0,\theta_1) = \widetilde{\varepsilon}$ and $B_{M_1}(x,\theta_1)|_{|x| = \widetilde{\varrho}_0} = 0$. Thus, it is simple to obtain the relation
\begin{equation}\label{eq:RELATIONBETWEENPARAMETERSBARINITIAL}
(M_1^{\gamma}\theta_1)^{\alpha} = (k/C_1)^{\frac{(p-1)N}{p}}\widetilde{\varrho}_0^N
\end{equation}
and the constants
\begin{equation}\label{eq:CHOICEOFBARENBLATTPARAMETERSLEVELSETS}
\theta_1 = k^{p-1}\widetilde{\varrho}_0^{\,p} \widetilde{\varepsilon}^{\,-\gamma} \qquad \text{ and } \qquad M_1 = C_1^{-\frac{p-1}{\gamma}}(k/C_1)^{\frac{(p-1)N}{p}}\widetilde{\varrho}_0^N \widetilde{\varepsilon},
\end{equation}
where $C_1$ is the constant corresponding to the profile $F_1(\cdot)$ and
\[
\alpha = \frac{1}{\gamma + p/N}, \quad k = \frac{\gamma}{p}\Big(\frac{\alpha}{N}\Big)^{\frac{1}{p-1}}
\]
are defined in Section \ref{PRELIMINARIESINTRO}. Then, we consider the solution of the problem
\[
\begin{cases}
\begin{aligned}
\partial_tw = \Delta_p w^m  \qquad\qquad\; &\text{in } \RR^N\times(0,\infty)\\
w(x,0) = B_{M_1}(x,\theta_1) \quad &\text{in } \RR^N,
\end{aligned}
\end{cases}
\]
i.e., $w(x,t) = B_{M_1}(x,\theta_1 + t)$, which satisfies $u(x,t) \geq w(x,t)$ in $\RR^N\times(0,\infty)$ thanks to the Maximum Principle (recall that by construction we have $w(x,0) \leq \widetilde{u}_0(x)$).

\noindent Now, we take $t_1 > 0$ and $n_1 \in \NN$ satisfying
\begin{equation}\label{eq:CONDITIONSONT1ANDN1}
t_1 \geq 2^{\frac{N(p-1)}{\alpha p}}\theta_1 \bigg(\frac{\widetilde{\varrho}_1}{\widetilde{\varrho}_0} \bigg)^{\frac{N}{\alpha}} \qquad \text{and} \qquad
n_1 \geq 2^{\frac{p-1}{\gamma}} \Bigg(1 + \frac{t_1}{\theta_1}\Bigg)^{\alpha}.
\end{equation}
Thus, since the profile of the Barenblatt solutions is decreasing, in order to have $u(x,t_1) \geq \widetilde{\varepsilon}/n_1$ in $\{|x| \leq \widetilde{\varrho}_1\}$, it sufficient to impose
\begin{equation}\label{eq:SUFFICIENTCONDITIONLEMMA1EXPANDINGSMALLERLEVELSETS}
w(x,t_1)|_{|x| = \widetilde{\varrho}_1} = B_{M_1}(x,\theta_1 + t_1)|_{|x| = \widetilde{\varrho}_1} \geq \widetilde{\varepsilon}/n_1.
\end{equation}
Now, using the relations in \eqref{eq:RELATIONBETWEENPARAMETERSBARINITIAL} and \eqref{eq:CHOICEOFBARENBLATTPARAMETERSLEVELSETS}, it is not difficult to compute
\[
\begin{aligned}
& B_{M_1}(x,\theta_1 + t_1)|_{|x| = \widetilde{\varrho}_1} = M_1 B_1(x,M_1^{\gamma}(\theta_1 + t_1))|_{|x| = \widetilde{\varrho}_1}  \\
= \;& \frac{M_1}{(M_1^{\gamma}\theta_1)^{\alpha}}\Bigg(\frac{\theta_1}{\theta_1 + t_1}\Bigg)^{\alpha} \Bigg\{ C_1 - k\bigg(\frac{\widetilde{\varrho}_1}{\widetilde{\varrho}_0} \bigg)^{\frac{p}{p-1}}\Bigg[ (M_1^{\gamma}\theta_1)^{-\frac{\alpha}{N}}\Bigg( \frac{\theta_1}{\theta_1 + t_1} \Bigg)^{\frac{\alpha}{N}}\Bigg]^{\frac{p}{p-1}}\Bigg\}_+^{\frac{p-1}{\gamma}} \\
= \;&\widetilde{\varepsilon} \,\Bigg(\frac{\theta_1}{\theta_1 + t_1}\Bigg)^{\alpha}\Bigg[1 - \bigg(\frac{\widetilde{\varrho}_1}{\widetilde{\varrho}_0} \bigg)^{\frac{p}{p-1}}\Bigg( \frac{\theta_1}{\theta_1 + t_1} \Bigg)^{\frac{\alpha p}{N(p-1)}} \Bigg]_+^{\frac{p-1}{\gamma}}.
\end{aligned}
\]
So, requiring \eqref{eq:SUFFICIENTCONDITIONLEMMA1EXPANDINGSMALLERLEVELSETS} is equivalent to
\begin{equation}\label{eq:CONDITIONONT1ANDN1GREATERTHANVEP}
\Bigg[1 - \Bigg( \frac{\theta_1}{\theta_1 + t_1} \Bigg)^{\frac{\alpha p}{N(p-1)}} \bigg(\frac{\widetilde{\varrho}_1}{\widetilde{\varrho}_0} \bigg)^{\frac{p}{p-1}}\Bigg]_+^{\frac{p-1}{\gamma}} \geq \;\; \frac{1}{\,n_1}\Bigg(1 + \frac{t_1}{\theta_1}\Bigg)^{\alpha}.
\end{equation}
Since the first condition in \eqref{eq:CONDITIONSONT1ANDN1} assures that the term on the left side of the previous inequality is larger than $2^{-(p-1)/\gamma}$, we have that a sufficient condition so that \eqref{eq:CONDITIONONT1ANDN1GREATERTHANVEP} is satisfied is
\[
n_1 \geq 2^{\frac{p-1}{\gamma}} \Bigg(1 + \frac{t_1}{\theta_1}\Bigg)^{\alpha},
\]
which is our second assumption in \eqref{eq:CONDITIONSONT1ANDN1}, and so our proof is complete. $\Box$

\bigskip

We proceed in our work, proving that, for all $\widetilde{\varepsilon} > 0$ small enough, the super-level sets $E_{\widetilde{\varepsilon}}^+(t) := \{x \in \RR^N: \; u(x,t) \geq \widetilde{\varepsilon}\}$ of the solution $u = u(x,t)$ of problem \eqref{eq:REACTIONDIFFUSIONEQUATIONPLAPLACIAN} with initial datum \eqref{eq:INITIALDATUMTESTLEVELSETS} do not contract in time for $t$ large enough and, in particular, we will show that for all $m > 0$ and $p > 1$ such that $\gamma > 0$ and for all $\widetilde{\varrho}_0 > 0$, it holds
\[
\{|x| \leq \widetilde{\varrho}_0/2\} \subset E_{\widetilde{\varepsilon}}^+(t), \quad \text{for large times}.
\]
This result highlights the role of the reaction term $f(\cdot)$. Indeed, the solution of the ``pure diffusive'' equation converges to zero as $t \to \infty$, whilst the presence of the function $f(\cdot)$ is sufficient to guarantee the strict positivity of the solution in a (small) compact set for large times.
%
%
%
%
%
%
\begin{lem}\label{LEMMAEXPANDINGLEVELSETS}
Let $m > 0$ and $p > 1$ such that $\gamma > 0$ and let $N \geq 1$. Then, for all $\widetilde{\varrho}_0 > 0$, there exist $t_2 > 0$ and $0 < \widetilde{\varepsilon}_0 < 1$ which depend only on $m$, $p$, $N$, $f$ and $\widetilde{\varrho}_0$, such that for all $0 < \widetilde{\varepsilon} \leq \widetilde{\varepsilon}_0$, the solution $u(x,t)$ of problem \eqref{eq:REACTIONDIFFUSIONEQUATIONPLAPLACIAN} with initial datum \eqref{eq:INITIALDATUMTESTLEVELSETS} satisfies
\[
u(x,jt_2) \geq \widetilde{\varepsilon} \quad \text{in } \{|x| \leq \widetilde{\varrho}_0/2\}, \text{ for all } \; j \in \NN_+ = \{1,2,\ldots\}.
\]
\end{lem}
\emph{Proof.} We prove the assertion by induction on $j \in \NN_+$. We ask the reader to note that we follow the ideas previously used by Cabr\'e and Roquejoffre \cite{C-R2:art} and, later, in \cite{S-V:art} adapting them to our (quite different) setting.

\emph{Step0.} In this step, we introduce some basic definitions and quantities we will use during the proof.

\noindent We fix $j = 1$ and $\widetilde{\varrho}_0 > 0$. Moreover, let $0 < \delta < 1$ and set $\lambda := f(\delta)/\delta$. Then take $t_2$ large enough such that
\begin{equation}\label{eq:PRELIMINARYCONDITIONSONT0}
e^{\lambda t_2} \geq 2^{\frac{p-1}{\gamma}}\Bigg(1 + \frac{\tau(t_2)}{\widetilde{C}_1}\Bigg)^{\alpha}
\end{equation}
where $\alpha > 0$ and $k > 0$ are defined in Section \ref{PRELIMINARIESINTRO} (see also the beginning of the proof of Lemma \ref{EXPANDINGLEVELSETS2}). The constant $\widetilde{C}_1$ is defined by the formula
\[
\widetilde{C}_1:= k^{p-1}\widetilde{\varrho}_0^{\,p}
\]
and the function $\tau(\cdot)$ is defined as follows:
\begin{equation}\label{eq:FIRSTCHANGEOFVARIABLELEVELSETS}
\tau(t) = \frac{1}{\gamma\lambda}\Big[ e^{\gamma\lambda t} - 1\Big], \quad \text{for } t \geq 0.
\end{equation}
Then, we set $\widetilde{\varepsilon}_0 := \delta e^{-f'(0)t_2}$ and, finally, we fix $0 < \widetilde{\varepsilon} \leq \widetilde{\varepsilon}_0$. Note that the choices in \eqref{eq:PRELIMINARYCONDITIONSONT0} are admissible since $\gamma\alpha < 1$ and $t_2 > 0$ does not depend on $\widetilde{\varepsilon}$.

\emph{Step1.} Construction of a sub-solution of problem \eqref{eq:REACTIONDIFFUSIONEQUATIONPLAPLACIAN}, \eqref{eq:INITIALDATUMTESTLEVELSETS} in $\RR^N\times[0,t_2]$.
\\
First of all, as we did at the beginning of the proof of Lemma \ref{EXPANDINGLEVELSETS2}, we construct a Barenblatt solution of the form $B_{M_1}(x,\theta_1)$ such that $B_{M_1}(x,\theta_1) \leq \widetilde{u}_0(x)$ for all $x \in \RR^N$. Evidently, we obtain the same formulas for $M_1 > 0$ and $\theta_1 >0$ (see \eqref{eq:RELATIONBETWEENPARAMETERSBARINITIAL} and \eqref{eq:CHOICEOFBARENBLATTPARAMETERSLEVELSETS}). Before proceeding, we note that, using \eqref{eq:CHOICEOFBARENBLATTPARAMETERSLEVELSETS} and the fact that $\widetilde{\varepsilon} < 1$, it is simple to get
\begin{equation}\label{eq:BOUNDFORTHETA1LEVELSETS}
\theta_1 \geq \widetilde{C}_1 > 0.
\end{equation}
Now, consider the change of time variable $\tau = \tau(t)$ defined in \eqref{eq:FIRSTCHANGEOFVARIABLELEVELSETS} and the ``linearized'' problem
\begin{equation}\label{eq:FIRSTLINEARIZEDPROBLEMLEVELSETS}
\begin{cases}
\begin{aligned}
\partial_tw = \Delta_p w^m + \lambda w \quad &\text{in } \RR^N\times(0,\infty)\\
w(x,0) = \widetilde{u}_0(x) \qquad\;\; &\text{in } \RR^N.
\end{aligned}
\end{cases}
\end{equation}
Then, the function $\widetilde{w}(x,\tau) = e^{-\lambda t}w(x,t)$ solves the problem
\[
\begin{cases}
\begin{aligned}
\partial_{\tau}\widetilde{w} = \Delta_p \widetilde{w}^m  \;\,\qquad &\text{in } \RR^N\times(0,\infty)\\
\widetilde{w}(x,0) = \widetilde{u}_0(x) \;\;\quad &\text{in } \RR^N.
\end{aligned}
\end{cases}
\]
Since $B_{M_1}(x,\theta_1) \leq \widetilde{u}_0(x) \leq \widetilde{\varepsilon}$ for all $x \in \RR^N$, from the Maximum Principle we get
\begin{equation}\label{eq:INEQUALITYSUBSOLUTIONLEVELSETS}
B_{M_1}(x,\theta_1 + \tau) \leq \widetilde{w}(x,\tau) \leq \widetilde{\varepsilon} \quad \text{in } \RR^N\times(0,\infty).
\end{equation}
Hence, using the concavity of $f$ and the second inequality in \eqref{eq:INEQUALITYSUBSOLUTIONLEVELSETS} we get
\[
w(x,t) = e^{\lambda t}\widetilde{w}(x,\tau) \leq \widetilde{\varepsilon} e^{f'(0)t} \leq \widetilde{\varepsilon}_0 e^{f'(0)t_2} = \delta, \quad \text{in } \RR^N\times[0,t_2]
\]
and so, since $w \leq \delta$ implies $f(\delta)/\delta \leq f(w)/w$, we have that $w$ is a sub-solution of problem \eqref{eq:REACTIONDIFFUSIONEQUATIONPLAPLACIAN}, \eqref{eq:INITIALDATUMTESTLEVELSETS} in $\RR^N\times[0,t_2]$. Finally, using the first inequality in \eqref{eq:INEQUALITYSUBSOLUTIONLEVELSETS}, we obtain
\begin{equation}\label{eq:FIRSTINEQUALITYFORCHOOSINGEPSLEVELSETS}
u(x,t) \geq e^{\lambda t} \widetilde{w}(x,\tau) \geq e^{\lambda t} B_{M_1}(x,\theta_1 + \tau), \quad \text{in } \RR^N\times[0,t_2].
\end{equation}

\emph{Step2.} In this step, we verify that the assumptions \eqref{eq:PRELIMINARYCONDITIONSONT0} on $t_2 > 0$ (depending only on $m$, $p$, $N$, $f$ and $\widetilde{\varrho}_0$) are sufficient to prove $u(x,t_2) \geq \widetilde{\varepsilon}$ in the set  $\{|x| \leq \widetilde{\varrho}_0/2\}$.
\\
First of all, we note that, from the second inequality in \eqref{eq:FIRSTINEQUALITYFORCHOOSINGEPSLEVELSETS} and since the profile of the Barenblatt solution is decreasing, it is clear that it is sufficient to have $t_2$ such that
\begin{equation}\label{eq:FIRSTSUFFICIENTINEQUALITYT0LEVELSETS}
e^{\lambda t_2} B_{M_1}(x,\theta_1 + \tau_2)|_{|x| = \widetilde{\varrho}_0/2} \geq \widetilde{\varepsilon},
\end{equation}
where $\tau_2 := \tau(t_2)$. Using the relations in \eqref{eq:RELATIONBETWEENPARAMETERSBARINITIAL} and \eqref{eq:CHOICEOFBARENBLATTPARAMETERSLEVELSETS} and proceeding as in the proof of Lemma \ref{EXPANDINGLEVELSETS2}, we get
\[
e^{\lambda t_2} B_{M_1}(x,\theta_1 + \tau_2)|_{|x| = \widetilde{\varrho}_0/2} = \widetilde{\varepsilon}\, e^{\lambda t_2}\Bigg(\frac{\theta_1}{\theta_1 + \tau_2}\Bigg)^{\alpha}\Bigg[1 - 2^{-\frac{p}{p-1}}\Bigg( \frac{\theta_1}{\theta_1 + \tau_2} \Bigg)^{\frac{\alpha p}{N(p-1)}} \Bigg]_+^{\frac{p-1}{\gamma}}
\]
Hence, we have that \eqref{eq:FIRSTSUFFICIENTINEQUALITYT0LEVELSETS} is equivalent to
\begin{equation}\label{eq:SECONDSUFFICIENTINEQUALITYT0LEVELSETS}
e^{\lambda t_2} \Bigg[1 - 2^{-\frac{p}{p-1}}\Bigg( \frac{\theta_1}{\theta_1 + \tau_2} \Bigg)^{\frac{\alpha p}{N(p-1)}} \Bigg]_+^{\frac{p-1}{\gamma}} \geq \Bigg(1 + \frac{\tau_2}{\theta_1}\Bigg)^{\alpha}.
\end{equation}
Since for all fixed $\tau > 0$, the function $\theta/(\theta + \tau)$ satisfies
\[
\frac{\theta}{\theta + \tau} \leq 1 \leq 2^{\frac{N}{\alpha p}}, \quad \text{ for all } \theta \geq 0
\]
and since $\theta_1 \geq \widetilde{C}_1$ (see \eqref{eq:BOUNDFORTHETA1LEVELSETS}), it is simple to deduce that a sufficient condition so that \eqref{eq:SECONDSUFFICIENTINEQUALITYT0LEVELSETS} is satisfied is
\[
e^{\lambda t_2} \geq 2^{\frac{p-1}{\gamma}}\Bigg(1 + \frac{\tau_2}{\widetilde{C}_1}\Bigg)^{\alpha}
\]
Note that it does not depend on $0 < \widetilde{\varepsilon} \leq \widetilde{\varepsilon}_0$ and it is exactly the assumption \eqref{eq:PRELIMINARYCONDITIONSONT0} on $t_2 > 0$. The proof of the case $j = 1$ is completed.

However, before studying the iteration of this process, we need to do a last effort. Let $\widetilde{\varrho}_1 \geq \widetilde{\varrho}_0/2$ be such that
\[
e^{\lambda t_2} B_{M_1}(x,\theta_1 + \tau_2)|_{|x| = \widetilde{\varrho}_1 } = \widetilde{\varepsilon}
\]
and introduce the function
\[
v_0(x) :=
\begin{cases}
\begin{aligned}
\widetilde{\varepsilon} \;\;\;\quad\qquad\qquad\qquad\qquad &\text{if } |x| \leq \widetilde{\varrho}_1 \\
e^{\lambda t_2} B_{M_1}(x,\theta_1 + \tau_2) \qquad &\text{if } |x| > \widetilde{\varrho}_1.
\end{aligned}
\end{cases}
\]
A direct computation (which we leave as an exercise for the interested reader) shows that condition \eqref{eq:PRELIMINARYCONDITIONSONT0} is sufficient to prove
\begin{equation}\label{eq:IMPORTANTESTIMATEONINITIALDATA}
u(x,t_2) \geq v_0(x) \geq B_{M_1}(x, \theta_1) \quad \text{in } \RR^N.
\end{equation}
\paragraph{Iteration.} We suppose to have proved that the solution of problem \eqref{eq:REACTIONDIFFUSIONEQUATIONPLAPLACIAN}, \eqref{eq:INITIALDATUMTESTLEVELSETS} satisfies
\[
u(x,jt_2) \geq \widetilde{\varepsilon} \quad \text{in } \{|x| \leq \widetilde{\varrho}_0/2\}, \text{ for some }  j \in \NN_+
\]
with the property
\begin{equation}\label{eq:IMPORTANTESTIMATEONINITIALDATA1}
u(x,jt_2) \geq v_0(x) \geq B_{M_1}(x, \theta_1) \quad \text{in } \RR^N
\end{equation}
which we can assume since it holds in the case $j = 1$ (see \eqref{eq:IMPORTANTESTIMATEONINITIALDATA}) and we prove
\begin{equation}\label{eq:FINALTHESISITERATION}
u(x,(j+1)t_2) \geq \widetilde{\varepsilon} \quad \text{in } \{|x| \leq \widetilde{\varrho}_0/2\}.
\end{equation}
Thanks to \eqref{eq:IMPORTANTESTIMATEONINITIALDATA1}, it follows that the solution $v(x,t)$ of the problem
\begin{equation}\label{eq:KCAUCHYPROBLEMLEVELSETS}
\begin{cases}
\begin{aligned}
\partial_tv = \Delta_pv^m + f(v) \quad &\text{in } \RR^N\times(0,\infty)\\
v(x,0) = v_0(x) \quad\qquad &\text{in } \RR^N
\end{aligned}
\end{cases}
\end{equation}
satisfies $u(x,t + jt_2) \geq v(x,t)$ in $\RR^N\times[0,\infty)$ by the Maximum Principle. Consequently, we can work with the function $v(x,t)$ and proceeding similarly as before.

\emph{Step1'.} Construction of a sub-solution of problem \eqref{eq:KCAUCHYPROBLEMLEVELSETS} in $\RR^N\times[0,t_2]$.
\\
This step is identical to \emph{Step1} (case $j = 1$). However, we ask the reader to keep in mind that now we are building a sub-solution of the function $v = v(x,t)$. This means that we consider the ``linearized'' problem
\[
\begin{cases}
\begin{aligned}
\partial_tw = \Delta_p w^m + \lambda w \quad\, &\text{in } \RR^N\times(0,\infty)\\
w(x,0) = v_0(x) \;\;\;\qquad &\text{in } \RR^N
\end{aligned}
\end{cases}
\]
where $\lambda :=f(\delta)/\delta$ and, using again the change of variable \eqref{eq:FIRSTCHANGEOFVARIABLELEVELSETS}, we deduce that the function $\widetilde{w}(x,\tau) = e^{-\lambda t}w(x,t)$ solves the problem
\[
\begin{cases}
\begin{aligned}
\partial_{\tau}\widetilde{w} = \Delta_p \widetilde{w}^m  \;\,\qquad &\text{in } \RR^N\times(0,\infty)\\
\widetilde{w}(x,0) = v_0(x) \;\;\quad &\text{in } \RR^N.
\end{aligned}
\end{cases}
\]
Since $B_{M_1}(x,\theta_1) \leq v_0(x) \leq \widetilde{\varepsilon}$ for all $x \in \RR^N$, from the Maximum Principle we get again
\[
B_{M_1}(x,\theta_1 + \tau) \leq \widetilde{w}(x,\tau) \leq \widetilde{\varepsilon} \quad \text{in } \RR^N\times(0,\infty)
\]
and, moreover, $w(x,t) \leq \delta$ in $\RR^N\times[0,t_2]$ which allows us to conclude that $w = w(x,t)$ is a sub-solution of problem \eqref{eq:KCAUCHYPROBLEMLEVELSETS} in $\RR^N\times[0,t_2]$. In particular, it holds
\[
u(x,(j+1)t_2) \geq v(x,t_2) \geq w(x,t_2) \geq e^{\lambda t_2} B_{M_1}(x,\theta_1 + \tau_2).
\]

\emph{Step2'.} This step is identical to \emph{Step2}, since we have to verify that
\[
e^{\lambda t_2} B_{M_1}(x,\theta_1 + \tau_2)|_{|x| = \widetilde{\varrho}_0/2} \geq \widetilde{\varepsilon}.
\]
Since we showed in \emph{Step2} that \eqref{eq:PRELIMINARYCONDITIONSONT0} is as sufficient condition so that the previous inequality is satisfied, we obtain \eqref{eq:FINALTHESISITERATION} and we conclude the proof.

\noindent However, in order to be precise and be sure that our iteration actually works, we have to prove that
\[
u(x,(j+1)t_2) \geq B_{M_1}(x,\theta_1) \quad \text{in } \RR^N,
\]
but this follows from the fact that $w(x,t_2) \geq v_0(x) \geq B_{M_1}(x,\theta_1)$ in $\RR^N$. $\Box$
%
%
%
%
%
%
\begin{cor}\label{LEMMANONDISCRETEVERSIONLEVELSETS1}
Let $m > 0$ and $p > 1$ such that $\gamma > 0$ and let $N \geq 1$. Then, there exist $t_2 > 0$ and $0 < \widetilde{\varepsilon}_0 < 1$ which depend only on $m$, $p$, $N$, $f$ and $\widetilde{\varrho}_0$ such that, for all $0 < \widetilde{\varepsilon} \leq \widetilde{\varepsilon}_0$, the solution $u(x,t)$ of problem \eqref{eq:REACTIONDIFFUSIONEQUATIONPLAPLACIAN} with initial datum \eqref{eq:INITIALDATUMTESTLEVELSETS1} satisfies
\[
u(x,t) \geq \widetilde{\varepsilon} \quad \text{in }  \{|x| \leq \widetilde{\varrho}_0/2\} \;\text{ for all } t \geq t_2.
\]
\end{cor}
\emph{Proof}. The previous lemma states that for all $\widetilde{\varrho}_0 > 0$ and for the sequence of times $t_j = (jt_2)_{j \in \NN_+}$, the solution of problem \eqref{eq:REACTIONDIFFUSIONEQUATIONPLAPLACIAN}, \eqref{eq:INITIALDATUMTESTLEVELSETS} reaches a positive value $\widetilde{\varepsilon}$ in the set $ \{|x| \leq \widetilde{\varrho}_0/2\}$, i.e.
\[
u(x,jt_2) \geq \widetilde{\varepsilon} \quad \text{in } \{|x| \leq \widetilde{\varrho}_0/2\}, \text{ for all } \; j \in \NN_+,
\]
for all $0 < \widetilde{\varepsilon} \leq \widetilde{\varepsilon}_0 = \delta e^{-f'(0)t_2}$. We can improve this result choosing a smaller $\widetilde{\varepsilon}_0 > 0$.

\noindent Indeed, since conditions \eqref{eq:PRELIMINARYCONDITIONSONT0}:
\[
e^{\lambda t_2} \geq 2^{\frac{(p-1)}{\gamma}}\Bigg(1 + \frac{\tau(t_2)}{\widetilde{C}_1}\Bigg)^{\alpha}
\]
are satisfied for all $t_2 \leq t \leq 2t_2$, we can repeat the same proof of Lemma \ref{LEMMAEXPANDINGLEVELSETS}, modifying the value of $\widetilde{\varepsilon}_0$ and choosing a different value $\underline{\widetilde{\varepsilon}}_0 = \delta e^{-2f'(0)t_2} > 0$, which is smaller but strictly positive for all $t_2 \leq t \leq 2t_2$. Hence, it turns out that for all $0 < \widetilde{\varepsilon} \leq \underline{\widetilde{\varepsilon}}_0$, it holds
\[
u(x,t) \geq \widetilde{\varepsilon} \quad \text{in } \{ |x| \leq \widetilde{\varrho}_0/2\},  \text{ for all } \; t_2 \leq t \leq 2t_2.
\]
Now, iterating this procedure as in the proof of Lemma \ref{LEMMAEXPANDINGLEVELSETS}, we do not have to change the value of $\underline{\widetilde{\varepsilon}}_0$ when $j \in \NN_+$ grows and so, for all $0 < \widetilde{\varepsilon} \leq \underline{\widetilde{\varepsilon}}_0$, we obtain
\[
u(x,t) \geq \widetilde{\varepsilon} \quad \text{in } \{|x| \leq \widetilde{\varrho}_0/2\}, \text{ for all } \; j \in \NN_+ \;\text{ and for all }\;\; j t_2 \leq t \leq (j+1)t_2.
\]
Hence, using the arbitrariness of $j \in \NN_+$, we end the proof of this corollary. $\Box$

\bigskip

Now, in order to prove our main instrument (Proposition \ref{LEMMANONDISCRETEVERSIONLEVELSETS}), needed in the study of the asymptotic behaviour of the general solutions of problem \eqref{eq:REACTIONDIFFUSIONEQUATIONPLAPLACIAN} (see Lemma \ref{LEMMACOMPACTCONVERGENCETO1ASYMPTOTICBEHAVIOURPGREATER2}), we combine Lemma \ref{EXPANDINGLEVELSETS2} and Corollary \ref{LEMMANONDISCRETEVERSIONLEVELSETS1}.
%
%
%
%
%
%
%
%
\paragraph{Proof of Proposition \ref{LEMMANONDISCRETEVERSIONLEVELSETS} (case $\boldsymbol{ \gamma > 0}$).} Fix $\widetilde{\varrho}_0 > 0$, $\widetilde{\varrho}_1 \geq \widetilde{\varrho}_0$ and consider the solution $u = u(x,t)$ of problem \eqref{eq:REACTIONDIFFUSIONEQUATIONPLAPLACIAN} with initial datum \eqref{eq:INITIALDATUMTESTLEVELSETS}. So, thanks to Lemma \ref{EXPANDINGLEVELSETS2}, for all $0 < \widetilde{\varepsilon}_1 < 1$, there exist $t_1 > 0$ and $n_1 \in \NN$ such that
\[
u(x,t_1) \geq \widetilde{\varepsilon} := \widetilde{\varepsilon}_1/n_1 \quad \text{in } \{|x| \leq \widetilde{\varrho}_1\}.
\]
Now, define the function
\[
\widetilde{v}_0(x) :=
\begin{cases}
\begin{aligned}
\widetilde{\varepsilon} \qquad &\text{if } |x| \leq \widetilde{\varrho}_1 \\
0  \qquad &\text{if } |x| > \widetilde{\varrho}_1
\end{aligned}
\end{cases}
\]
and note that $u(x,t_1) \geq \widetilde{v}_0(x)$ in $\RR^N$. Hence, the solution of the problem
\[
\begin{cases}
\begin{aligned}
v_t = \Delta_pv^m + f(v) \quad &\text{in } \RR^N\times(0,\infty)\\
v(x,0) = \widetilde{v}_0(x) \;\qquad &\text{in } \RR^N
\end{aligned}
\end{cases}
\]
satisfies $u(x,t_1+t) \geq v(x,t)$ in $\RR^N\times(0,\infty)$. Finally, applying Corollary \ref{LEMMANONDISCRETEVERSIONLEVELSETS1}, we get
\[
u(x,t) \geq \widetilde{\varepsilon} \quad \text{in }  \{|x| \leq \widetilde{\varrho}_1/2\} \;\text{ for all } t \geq t_0
\]
with $t_0 := t_1 + t_2$ (note that, since $n_1 \in \NN$ can be chosen larger,  $\widetilde{\varepsilon} > 0$ is arbitrarily small and we ``end up'' in the hypotheses of Corollary \ref{LEMMANONDISCRETEVERSIONLEVELSETS1}). $\Box$
%
%
%
%
%
%
%
%
%
%
%
%
%
%
\section{\texorpdfstring{\boldmath}{}Case \texorpdfstring{$\gamma = 0$}{gamma}. Existence of minimal expanding super-level sets} \label{SECTIONLINEAREXPANSIONSUPERLEVELSETSPSEUDO}
In this section we study problem \eqref{eq:REACTIONDIFFUSIONEQUATIONPLAPLACIAN}, with reaction term $f(\cdot)$ satisfying \eqref{eq:ASSUMPTIONSONTHEREACTIONTERM} and with initial datum:
\begin{equation}\label{eq:INITIALDATUMTESTLEVELSETSPSEUDO}
\widetilde{u}_0(x) :=
\begin{cases}
\begin{aligned}
\widetilde{\varepsilon} \;\;\qquad\qquad\qquad &\text{if } |x| \leq \widetilde{\varrho}_0 \\
a_0e^{-b_0|x|^{\frac{p}{p-1}}} \qquad &\text{if } |x| > \widetilde{\varrho}_0,
\end{aligned}
\end{cases}
\qquad a_0 := \widetilde{\varepsilon} e^{b_0\widetilde{\varrho}_0^{\,\frac{p}{p-1}}}
\end{equation}
where $\widetilde{\varepsilon}$, $\widetilde{\varrho}_0$ and $b_0$ are positive numbers (note that $\widetilde{u}_0 \in L^1(\RR^N)$). We ask the reader to note that, alternatively, we can fix the constants $\widetilde{\varrho}_0$, $a_0$ and $b_0$ and obtaining $\widetilde{\varepsilon} > 0$ by the ``inverse'' of the second formula in \eqref{eq:INITIALDATUMTESTLEVELSETSPSEUDO}.

\noindent Our goal is to prove Proposition \ref{LEMMANONDISCRETEVERSIONLEVELSETS}. Anyway, in order to motivate our choice from the beginning, we start proving the following lemma which will be useful in the proof of Lemma \ref{LEMMACOMPACTCONVERGENCETO1ASYMPTOTICBEHAVIOURPGREATER2}.
%
%
%
%
%
%
\begin{lem}\label{LEMMAPLACINGBARENBLATTUNDERSOLUTIONPSEUDO}
Let $m > 0$ and $p > 1$ such that $\gamma = 0$ and let $N \geq 1$. Then for all  $\widetilde{\varrho}_0 > 0$, there exist $t_1 > 0$, $a_0 > 0$ and $b_0 > 0$ such that the solution $u(x,t)$ of problem \eqref{eq:REACTIONDIFFUSIONEQUATIONPLAPLACIAN} with initial datum \eqref{eq:ASSUMPTIONSONTHEINITIALDATUM} satisfies
\[
u(x,t_1) \geq \widetilde{u}_0(x) \quad \text{in } \RR^N
\]
where $\widetilde{u}_0(\cdot)$ is defined in \eqref{eq:INITIALDATUMTESTLEVELSETSPSEUDO}.
\end{lem}
\emph{Proof.} Let $u = u(x,t)$ the solution of problem \eqref{eq:REACTIONDIFFUSIONEQUATIONPLAPLACIAN}, \eqref{eq:ASSUMPTIONSONTHEINITIALDATUM} and consider the solution $v = v(x,t)$ of the ``pure diffusive'' Cauchy problem:
\[
\begin{cases}
\begin{aligned}
\partial_t v = \Delta_pv^m \;\quad\quad &\text{in } \RR^N\times(0,\infty) \\
v(x,0) = u_0(x) \quad &\text{in } \RR^N,
\end{aligned}
\end{cases}
\]
which satisfies $v(x,t) \leq u(x,t)$ in $\RR^N\times[0,\infty)$ thanks to the Maximum Principle.

\emph{Step1.} We begin by proving that for all $\tau > 0$, there exist $\varepsilon > 0$, $M > 0$ such that for all $\tau < t_1 < \tau + \varepsilon$ it holds
\[
u(x,t_1) \geq B_M(x,t_1 - \tau) \quad \text{in } \RR^N
\]
where $B_M(x,t)$ is the Barenblatt solution in the ``pseudo-linear'' case (presented in Section \ref{PRELIMINARIESINTRO}) with exponential form
\[
B_M(x,t) = C_M t^{-\frac{N}{p}} \exp \Big(-k|xt^{-\frac{1}{p}}|^{\frac{p}{p-1}}\Big),
\]
where $C_M>0$ is chosen depending on the mass $M$ and $k = (p-1)p^{-p/(p-1)}$.

\noindent Fix $\tau > 0$, $\varrho > 0$ and let $\varepsilon > 0$ (for the moment arbitrary). Furthermore, take a mass $M > 0$ such that $C_M \leq 1$. We want to compare the general solution $v = v(x,t)$ with the ``delayed'' Barenblatt solution $B_M = B_M(x,t-\tau)$ in the strip
\[
S = [\tau,\tau+\varepsilon]\times\{x\in\RR^N: |x| \geq \varrho\},
\]
in order to deduce $v(x,t) \geq B_M(x,t-\tau)$ in $S$. Hence, from the Maximum Principle, is sufficient to check this inequality on the parabolic boundary of $S$. Note that we have $v(x,\tau) \geq B_M(x,0) = 0$ for all $|x| \geq \varrho$. Now, in order to check that $v(\varrho,t) \geq B_M(\varrho,t-\tau)$ for all $\tau \leq t \leq \tau + \varepsilon$, we compute the derivative of the function $b(t) := B_M(\varrho,t-\tau)$ deducing that
\[
b'(t) \geq 0 \quad \text{if and only if} \quad t \leq \tau + \bigg[\frac{pk}{N(p-1)}\bigg]^{p-1}\varrho^{\,p}.
\]
Thus, choosing $\varepsilon \leq \{(pk)/[N(p-1)]\}^{p-1}\varrho^{\,p}$ and noting that $v(\varrho,t) \geq \underline{v} > 0$ for all $\tau \leq t \leq \tau+\varepsilon$ and some constant $\underline{v} > 0$ (the strict positivity follows from the Harnack Inequality proved in \cite{Kin-Kuusi:art}), we have that it is sufficient to prove $\underline{v} \geq B_M(\varrho,\varepsilon)$. Since, $B_M(\varrho,\varepsilon) \to 0$ as $\varepsilon \to 0$, we obtain the required inequality by taking eventually $\varepsilon > 0$ smaller. We ask to the reader to note that the assumption $C_M \leq 1$ guarantees that the choice of $\varepsilon > 0$ does not depend on $M > 0$.

\noindent Now, take $\tau < t_1 < \tau + \varepsilon$. We know $v(x,t_1) \geq B_M(x,t_1-\tau)$ if $|x|\geq \varrho$ and so, by taking $C_M > 0$ smaller (depending on $\varepsilon > 0$) and using again the positivity of the solution $v = v(x,t)$, it is straightforward to check that the same inequality holds for all $|x|\leq\varrho$, completing the proof of the first step.

\emph{Step2.} Fix $\widetilde{\varrho}_0 > 0$ and let $t_1>0$ be fixed as in the previous step. To end the proof, it is sufficient to prove
\[
B_M(x,t_1 - \tau) \geq \widetilde{u}_0(x) \quad \text{in } \{|x|\geq \widetilde{\varrho}_0\},
\]
since the profile of the Barenblatt solution is radially decreasing. A direct computation shows that the previous inequality is satisfied by taking
\[
a_0 \leq C_M(t_1 - \tau)^{-N/p} \quad \text{ and } \quad b_0 \geq k(t_1 - \tau)^{-1/(p-1)}
\]
and so, the proof is complete. $\Box$

\bigskip

Following the ideas of the previous section (see Lemma \ref{LEMMAEXPANDINGLEVELSETS}), we prove that for all $\widetilde{\varepsilon} > 0$ small enough, the super-level sets $E_{\widetilde{\varepsilon}}^+(t)$ of the solution $u = u(x,t)$ of problem \eqref{eq:REACTIONDIFFUSIONEQUATIONPLAPLACIAN} with initial datum \eqref{eq:INITIALDATUMTESTLEVELSETSPSEUDO} expand in time for $t$ large enough. More precisely, we will show that for all $\widetilde{\varrho}_0 > 0$ and for all $\widetilde{\varrho}_1 \geq \widetilde{\varrho}_0$, it holds
\[
\{|x| \leq \widetilde{\varrho}_1\} \subset E_{\widetilde{\varepsilon}}^+(t), \quad \text{for large times}.
\]
This will result a slightly simpler respect to the case $\gamma > 0$, since now we work with strictly positive Barenblatt solutions (see Section \ref{PRELIMINARIESINTRO}). For instance, we do not need a version of Lemma \ref{EXPANDINGLEVELSETS2}, which seems to be essential in the case $\gamma > 0$.
%
%
%
%
%
%
\begin{lem}\label{LEMMAEXPANDINGLEVELSETSPSEUDO}
Let $m > 0$ and $p > 1$ such that $\gamma = 0$ and let $N \geq 1$. Then, for all $\widetilde{\varrho}_0 > 0$ and for all $\widetilde{\varrho}_1 \geq \widetilde{\varrho}_0$, there exist $t_0 > 0$ and $0 < \widetilde{\varepsilon}_0 < 1$ which depend only on $m$, $p$, $N$, $f$, $\widetilde{u}_0$ and $\widetilde{\varrho}_1$, such that for all $0 < \widetilde{\varepsilon} \leq \widetilde{\varepsilon}_0$, the solution $u(x,t)$ of problem \eqref{eq:REACTIONDIFFUSIONEQUATIONPLAPLACIAN} with initial datum \eqref{eq:INITIALDATUMTESTLEVELSETSPSEUDO} satisfies
\[
u(x,jt_0) \geq \widetilde{\varepsilon} \quad \text{in } \{|x| \leq \widetilde{\varrho}_1\}, \text{ for all } \; j \in \NN_+ = \{1,2,\ldots\}.
\]
\end{lem}
\emph{Proof.} We prove the assertion by induction on $j = 1,2,\dots$.

\emph{Step0.} As in the proof of Lemma \ref{LEMMAEXPANDINGLEVELSETS}, this ``first'' step is devoted to the introduction of basic definitions and quantities we need during the proof.

\noindent We fix $j = 1$, $\widetilde{\varrho}_0 > 0$ and $\widetilde{\varrho}_1 \geq \widetilde{\varrho}_0$. Moreover, let $0 < \delta < 1$, set $\lambda := f(\delta)/\delta$ and fix $0 < \lambda_0 < \lambda$. Then take $t_0$ large enough such that
\begin{equation}\label{eq:PRELIMINARYCONDITIONSONT0PSEUDO}
t_0 \geq k^{p-1}\widetilde{\varrho}^{\,p}_1, \qquad e^{\lambda_0 t_0} \geq \bigg(\frac{\theta_1 + t_0}{\theta_1} \bigg)^{N/p} \quad \text{and} \quad t_0 \geq \frac{1}{\lambda-\lambda_0}
\end{equation}
and, furthermore:
\begin{equation}\label{eq:PRELIMINARYCONDITIONSONT0PSEUDO1}
t_0 \geq \bigg(\frac{2k}{b_0}\bigg)^{p-1} \quad \text{ and } \quad t_0 \geq \frac{b_0}{2(\lambda-\lambda_0)}\widetilde{\varrho}_1^{\,\frac{p}{p-1}}
\end{equation}
where (we ask the reader to note that we defined $k$ in Section \ref{PRELIMINARIESINTRO}, case $\gamma = 0$):
\[
k := (p-1)p^{-p/(p-1)} \quad \text{ and } \quad \theta_1 := \bigg( \frac{k}{b_0}\bigg)^{p-1}.
\]
Then, we set $\widetilde{\varepsilon}_0 := \delta e^{-f'(0)t_0}$ and, finally, we fix $0 < \widetilde{\varepsilon} \leq \widetilde{\varepsilon}_0$.

\emph{Step1.} Construction of a sub-solution of problem \eqref{eq:REACTIONDIFFUSIONEQUATIONPLAPLACIAN}, \eqref{eq:INITIALDATUMTESTLEVELSETSPSEUDO} in $\RR^N\times[0,t_0]$.
\\
First of all, as we did at the beginning of the proof of Lemma \ref{EXPANDINGLEVELSETS2} and Lemma \ref{LEMMAEXPANDINGLEVELSETS}, we construct a Barenblatt solution of the form $B_{M_1}(x,\theta_1)$ such that $B_{M_1}(x,\theta_1) \leq \widetilde{u}_0(x)$ for all $x \in \RR^N$. Note that in the case $\gamma = 0$, the Barenblatt solutions have exponential profile (see Section \ref{PRELIMINARIESINTRO}). Imposing again $B_{M_1}(0, \theta_1) = \widetilde{\varepsilon}$ we get the first relation $C_{M_1}\theta_1^{-N/p} = \widetilde{\varepsilon}$ which is sufficient to have $B_{M_1}(x,\theta_1) \leq \widetilde{u}_0(x)$ in $\{|x| \leq \widetilde{\varrho}_0\}$. On the other hand, for all $|x| \geq \widetilde{\varrho}_0$, we need to have:
\[
B_{M_1}(x,\theta_1) \leq a_0e^{-b_0|x|^{\frac{p}{p-1}}}, \qquad |x| \geq \widetilde{\varrho}_0.
\]
Using the relation $C_{M_1}\theta_1^{-N/p} = \widetilde{\varepsilon}$ and the fact that $\widetilde{\varepsilon} \leq a_0$, it is not difficult to see that a possible choice of parameters such that the previous inequality is satisfied too is
\begin{equation}\label{eq:CHOICEOFBARENBLATTPARAMETERSLEVELSETSPSEUDO}
\theta_1 = \bigg( \frac{k}{b_0}\bigg)^{p-1} \qquad \text{ and } \qquad C_{M_1} = \widetilde{\varepsilon}\, \bigg( \frac{k}{b_0}\bigg)^{\frac{N(p-1)}{p}}.
\end{equation}
Note that, respect to the case $\gamma > 0$, $\theta_1$ does not depend on $\widetilde{\varepsilon}>0$.
Now, as in the proof of Lemma \ref{LEMMAEXPANDINGLEVELSETS}, we consider the ``linearized'' problem
\begin{equation}\label{eq:FIRSTLINEARIZEDPROBLEMLEVELSETSPSEUDO}
\begin{cases}
\begin{aligned}
\partial_tw = \Delta_p w^m + \lambda w \quad &\text{in } \RR^N\times(0,\infty)\\
w(x,0) = \widetilde{u}_0(x) \;\;\qquad &\text{in } \RR^N
\end{aligned}
\end{cases}
\end{equation}
and we deduce that the function $\widetilde{w}(x,t) = e^{-\lambda t}w(x,t)$ solves the problem (note that when when $\gamma = 0$ the doubly nonlinear operator has homogeneity 1)
\[
\begin{cases}
\begin{aligned}
\partial_{t}\widetilde{w} = \Delta_p \widetilde{w}^m  \;\,\qquad &\text{in } \RR^N\times(0,\infty)\\
\widetilde{w}(x,0) = \widetilde{u}_0(x) \;\,\quad &\text{in } \RR^N.
\end{aligned}
\end{cases}
\]
Again we have that $B_{M_1}(x,\theta_1) \leq \widetilde{u}_0(x) \leq \widetilde{\varepsilon}$ for all $x \in \RR^N$ and so, by comparison we deduce
\begin{equation}\label{eq:INEQUALITYSUBSOLUTIONLEVELSETSPSEUDO}
B_{M_1}(x,\theta_1 + t) \leq \widetilde{w}(x,t) \leq \widetilde{\varepsilon} \quad \text{in } \RR^N\times(0,\infty).
\end{equation}
Hence, using the concavity of $f$ and the second inequality in \eqref{eq:INEQUALITYSUBSOLUTIONLEVELSETSPSEUDO} we get
\[
w(x,t) = e^{\lambda t}\widetilde{w}(x,t) \leq \widetilde{\varepsilon} e^{f'(0)t} \leq \widetilde{\varepsilon}_0 e^{f'(0)t_0} = \delta, \quad \text{in } \RR^N\times[0,t_0]
\]
and, consequently, $w = w(x,t)$ is a sub-solution of problem \eqref{eq:REACTIONDIFFUSIONEQUATIONPLAPLACIAN}, \eqref{eq:INITIALDATUMTESTLEVELSETSPSEUDO} in $\RR^N\times[0,t_0]$. Finally, using the first inequality in \eqref{eq:INEQUALITYSUBSOLUTIONLEVELSETSPSEUDO}, we obtain
\begin{equation}\label{eq:FIRSTINEQUALITYFORCHOOSINGEPSLEVELSETSPSEUDO}
u(x,t) \geq e^{\lambda t} B_{M_1}(x,\theta_1 + t), \quad \text{in } \RR^N\times[0,t_0].
\end{equation}

\emph{Step2.} In this step, we verify that the assumptions \eqref{eq:PRELIMINARYCONDITIONSONT0PSEUDO} on $t_0 > 0$ are sufficient to prove $u(x,t_0) \geq \widetilde{\varepsilon}$ in the set  $\{|x| \leq \widetilde{\varrho}_1\}$.
\\
First of all, we note that, from the second inequality in \eqref{eq:FIRSTINEQUALITYFORCHOOSINGEPSLEVELSETSPSEUDO} and since the profile of the Barenblatt solution is decreasing, it is clear that it is sufficient to have $t_0$ such that
\begin{equation}\label{eq:FIRSTSUFFICIENTINEQUALITYT0LEVELSETSPSEUDO}
e^{\lambda t_0} B_{M_1}(x,\theta_1 + t_0)|_{|x| = \widetilde{\varrho}_1} \geq \widetilde{\varepsilon}.
\end{equation}
Using the relations, we compute
\[
e^{\lambda t_0} B_{M_1}(x,\theta_1 + t_0)|_{|x| = \widetilde{\varrho}_1} = \widetilde{\varepsilon}\, \Bigg(\frac{\theta_1}{\theta_1 + t_0}\Bigg)^{N/p} \exp \Bigg[\lambda t_0 - k \bigg(\frac{\widetilde{\varrho}_1^{\,p}}{\theta_1 + t_0}\bigg)^{1/(p-1)} \Bigg].
\]
Hence, we have that \eqref{eq:FIRSTSUFFICIENTINEQUALITYT0LEVELSETSPSEUDO} is equivalent to
\begin{equation}\label{eq:SECONDSUFFICIENTINEQUALITYT0LEVELSETSPSEUDO}
\Bigg(\frac{\theta_1}{\theta_1 + t_0}\Bigg)^{N/p} \exp \Bigg[\lambda t_0 - k \bigg(\frac{\widetilde{\varrho}_1^{\,p}}{\theta_1 + t_0}\bigg)^{1/(p-1)} \Bigg] \geq 1.
\end{equation}
Now, using the first and the second relation in \eqref{eq:PRELIMINARYCONDITIONSONT0PSEUDO} it is not difficult to see that a sufficient condition so that \eqref{eq:SECONDSUFFICIENTINEQUALITYT0LEVELSETSPSEUDO} is satisfied is
\[
e^{(\lambda-\lambda_0)t_0 - 1} \geq 1,
\]
which is equivalent to the third assumption in \eqref{eq:PRELIMINARYCONDITIONSONT0PSEUDO} and so, we have $u(x,t_0) \geq \widetilde{\varepsilon}$ in $\{|x| \leq \widetilde{\varrho}_1\}$, i.e. the thesis for $j = 1$.

\noindent As we did in Lemma \ref{LEMMAEXPANDINGLEVELSETS}, we have to do a last step in order to conclude the case $j=1$. Let $\widetilde{\varrho}_2 \geq \widetilde{\varrho}_1$ such that
\[
e^{\lambda t_0} B_{M_1}(x,\theta_1 + t_0)|_{|x| = \widetilde{\varrho}_2} = \widetilde{\varepsilon}.
\]
Evidently, defining the function:
\[
v_0(x) :=
\begin{cases}
\begin{aligned}
\widetilde{\varepsilon} \;\;\,\quad\qquad\qquad\qquad\qquad &\text{if } |x| \leq \widetilde{\varrho}_2 \\
e^{\lambda t_0} B_{M_1}(x,\theta_1 + t_0) \qquad &\text{if } |x| > \widetilde{\varrho}_2,
\end{aligned}
\end{cases}
\]
it is sufficient to prove $v_0(x) \geq B_{M_1}(x,\theta_2)$ in $\{|x| > \widetilde{\varrho}_1\}$ and, using the second relation in \eqref{eq:PRELIMINARYCONDITIONSONT0PSEUDO}, it is simple to see that this is automatically satisfied if
\[
\exp\big[(\lambda - \lambda_0)t_0 \big] \geq \exp\big[b_0\widetilde{\varrho}_0^{\,\frac{p}{p-1}} + (kt_0^{-\frac{1}{p-1}} - b_0)|x|^{\frac{p}{p-1}}\big].
\]
Now, keeping in mind that $\widetilde{\varrho}_1 \geq \widetilde{\varrho}_0$ and using the first inequality in \eqref{eq:PRELIMINARYCONDITIONSONT0PSEUDO1}, we deduce that  a sufficient condition so that the previous inequality is satisfied is
\[
e^{(\lambda - \lambda_0)t_0} \geq e^{(b_0/2)\widetilde{\varrho}_1^{\,\frac{p}{p-1}}},
\]
which is equivalent to the second condition in \eqref{eq:PRELIMINARYCONDITIONSONT0PSEUDO1}. Hence, we have
\begin{equation}\label{eq:IMPORTANTESTIMATEONINITIALDATAPSEUDO}
u(x,t_0) \geq v_0(x) \geq B_{M_1}(x, \theta_1) \quad \text{in } \RR^N.
\end{equation}
\paragraph{Iteration.} We suppose to have proved that the solution of problem \eqref{eq:REACTIONDIFFUSIONEQUATIONPLAPLACIAN}, \eqref{eq:INITIALDATUMTESTLEVELSETSPSEUDO} satisfies
\[
u(x,jt_0) \geq \widetilde{\varepsilon} \quad \text{in } \{|x| \leq \widetilde{\varrho}_1\}, \text{ for some }  j \in \NN_+
\]
with the property
\begin{equation}\label{eq:IMPORTANTESTIMATEONINITIALDATA1PSEUDO1}
u(x,jt_0) \geq v_0(x) \geq B_{M_1}(x, \theta_1) \quad \text{in } \RR^N
\end{equation}
and we prove
\[
u(x,(j+1)t_0) \geq \widetilde{\varepsilon} \quad \text{in } \{|x| \leq \widetilde{\varrho}_1\}.
\]
As in the proof of Lemma \ref{LEMMAEXPANDINGLEVELSETS}, we have that \eqref{eq:IMPORTANTESTIMATEONINITIALDATA1PSEUDO1} implies that the solution $v = v(x,t)$ of the problem
\begin{equation}\label{eq:KCAUCHYPROBLEMLEVELSETSPSEUDO}
\begin{cases}
\begin{aligned}
v_t = \Delta_pv^m + f(v) \quad &\text{in } \RR^N\times(0,\infty)\\
v(x,0) = v_0(x) \;\qquad &\text{in } \RR^N
\end{aligned}
\end{cases}
\end{equation}
satisfies $u(x,t + jt_0) \geq v(x,t)$ in $\RR^N\times[0,\infty)$ by the Maximum Principle and this allows us to work with the function $v = v(x,t)$. The rest of the proof is almost identical to the case $\gamma > 0$ (see Lemma  \ref{LEMMAEXPANDINGLEVELSETS}) and we leave the details to the interested reader. $\Box$
\paragraph{Proof of Proposition \ref{LEMMANONDISCRETEVERSIONLEVELSETSPSEUDO} (case $\boldsymbol{\gamma = 0}$).} The proof is identical to the proof of Corollary \ref{LEMMANONDISCRETEVERSIONLEVELSETS1}, substituting condition \eqref{eq:PRELIMINARYCONDITIONSONT0} with conditions \eqref{eq:PRELIMINARYCONDITIONSONT0PSEUDO} and \eqref{eq:PRELIMINARYCONDITIONSONT0PSEUDO1}, $\widetilde{\varrho}_0/2$ with $\widetilde{\varrho}_1$ and $t_2$ with $t_0$. $\Box$

%
%
%
%
%
%
%
%
%
%
%
%
%
%
\section{Asymptotic behaviour. Convergence to 1 on compact sets}\label{SECTIONASYMPTOTICBEHAVIOUR}
In this section, we prove a technical lemma (see Lemma \ref{LEMMACOMPACTCONVERGENCETO1ASYMPTOTICBEHAVIOURPGREATER2}) which will be crucial in the study of the asymptotic behaviour of the solution of the Cauchy problem \eqref{eq:REACTIONDIFFUSIONEQUATIONPLAPLACIAN}, \eqref{eq:ASSUMPTIONSONTHEINITIALDATUM} with reaction term $f(\cdot)$ satisfying \eqref{eq:ASSUMPTIONSONTHEREACTIONTERM}. We make the additional assumption $u_0(0) = \max_{x \in \RR}u_0(x)$ (note that this choice is clearly admissible by performing a translation of the $x$-axis and it will avoid tedious technical work in the rest of the section).

\noindent In what follows, we show that a general solution (for all $\gamma \geq 0$) converges to one on every compact set of $\RR^N$ (``large enough'') for large times.
%
%
%
%
%
%
\begin{lem}\label{LEMMACOMPACTCONVERGENCETO1ASYMPTOTICBEHAVIOURPGREATER2}
Let $m>1$ and $p>1$ such that $\gamma \geq 0$ and let $N \geq 1$. Then, for all $\varepsilon > 0$, there exist $0 < \widetilde{a}_{\varepsilon} < 1$ and $\widetilde{\varrho}_{\varepsilon} > 0$ (which depend on $\varepsilon > 0$) such that for all $\widetilde{\varrho}_1 \geq \widetilde{\varrho}_{\varepsilon}$, there exists $t_1 > 0$ (depending on $\widetilde{\varrho}_1$, $\varepsilon$) such that it holds
\[
u(x,t) \geq 1 - \varepsilon \quad \text{in } \{|x| \leq \widetilde{a}_{\varepsilon}\widetilde{\varrho}_1 \}, \text{ for all } t \geq t_1.
\]
\end{lem}
\emph{Proof.} We begin supposing $\gamma > 0$. Fix $\varepsilon > 0$.

\emph{Step1.} In this step, we construct a sub-solution of problem \eqref{eq:REACTIONDIFFUSIONEQUATIONPLAPLACIAN} on a domain of the form $\Omega\times(t_1,\infty)$, where $\Omega$ is a ball in $\RR^N$ and $t_1 > 0$.
\\
First of all, we fix $\widetilde{\varrho}_2 = 2\widetilde{\varrho}_1 > 0$ arbitrarily large and we apply Proposition \ref{LEMMANONDISCRETEVERSIONLEVELSETS} deducing the existence of a value $\widetilde{\varepsilon} > 0$ and a time $t_0 > 0$ such that
\[
u(x,t) \geq \widetilde{\varepsilon} \quad \text{in } \Omega\times[t_0,\infty)
\]
where $\Omega := \{ |x| \leq \widetilde{\varrho}_1 \}$.
Hence, since $\widetilde{\varepsilon} \leq u \leq 1$ and remembering that the function $f(\cdot)$ is concave, we can deduce a linear bound from below for the reaction term
\[
f(u) \geq q(1-u) \quad  \text{in } \Omega\times[t_0,\infty), \quad \text{where } q := f(\widetilde{\varepsilon})/(1 - \widetilde{\varepsilon}).
\]
Now, for all $\widetilde{\varrho}_1 > 0$, we consider a time $t_1 \geq t_0$ (note that $t_1$ is now ``almost'' arbitrary and its precise value depending on $\varepsilon > 0$ will be specified later). Thus, the solution of the problem
\begin{equation}\label{eq:PARABOLICCOMPARISONSUBSOLUTIONASYMPTOTICBEHAVIOURPGREATER2}
\begin{cases}
\begin{aligned}
\partial_tv = \Delta_p v^m + q(1-v) \quad &\text{in } \Omega \times (t_1,\infty)\\
v = \widetilde{\varepsilon} \;\;\qquad\qquad\qquad\qquad &\text{in } \partial\Omega \times (t_1,\infty) \\
v(x,t_1) = \widetilde{\varepsilon} \qquad\qquad\qquad &\text{in } \Omega
\end{aligned}
\end{cases}
\end{equation}
satisfies $\widetilde{\varepsilon} \leq v(x,t) \leq u(x,t) \leq 1$ in $\Omega \times [t_1,\infty)$ by the Maximum Principle. Furthermore, since for all fixed $\tau > 0$, the function $w(x,t) = v(x,\tau + t)$ satisfies the equation in \eqref{eq:PARABOLICCOMPARISONSUBSOLUTIONASYMPTOTICBEHAVIOURPGREATER2} and $w(x,t_1) \geq v(x,t_1)$, it follows $v(x,\tau + t) \geq v(x,t)$ in $\Omega \times [t_1,\infty)$, i.e., for all $x \in \Omega$, the function $v(x,\cdot)$ is non-decreasing. Consequently, since $v$ is uniformly bounded, there exists the uniform limit $v_{\infty}(x) := \lim_{t\to\infty} v(x,t)$ and it solves the elliptic problem
\begin{equation}\label{eq:ELLIPTICCOMPARISONSUBSOLUTIONASYMPTOTICBEHAVIOURPGREATER2}
\begin{cases}
\begin{aligned}
- \Delta_p v_{\infty}^m = q(1-v_{\infty}) \quad &\text{in } \Omega \\
v_{\infty} = \widetilde{\varepsilon} \;\;\qquad\qquad\qquad &\text{in } \partial\Omega
\end{aligned}
\end{cases}
\end{equation}
in the weak sense.

\emph{Step2.} In this step, we define the constants $a_{\varepsilon}$ and $\widetilde{\varrho}_{\varepsilon}$ and we complete the proof of the lemma. The value of these constants comes from the construction of a particular sub-solution of the elliptic problem \eqref{eq:ELLIPTICCOMPARISONSUBSOLUTIONASYMPTOTICBEHAVIOURPGREATER2}. Since our argument is quite technical, we try to sketch it skipping some computations that can be verified directly by the reader.

We look for a sub-solution of the elliptic problem \eqref{eq:ELLIPTICCOMPARISONSUBSOLUTIONASYMPTOTICBEHAVIOURPGREATER2} in the form
\[
w^m(r) = a\big[e^{g(r)} - 1\big],
\]
where $r = |x|$, $x \in \RR^N$ while the function $g(\cdot)$ and the constant $a > 0$ are taken as follows:
\[
g(r) := 1 - \bigg(\frac{r}{\widetilde{\varrho}_1}\bigg)^{\lambda}, \quad \text{with } \lambda:= \frac{p}{p-1} \quad \text{and} \quad \frac{(1-\varepsilon/2)^m}{e - 1} < a < \frac{1}{e - 1}.
\]
Note that the radially decreasing function $w(\cdot)$ is well defined in $[0,\widetilde{\varrho}_1]$ and, moreover, $w(r) = 0$ on the boundary $\partial\Omega = \{|x| = \widetilde{\varrho}_1\}$.

\noindent Now, we define the value $\widetilde{\varrho}_{\varepsilon}$ (note it is well defined and positive thanks to the assumption on $a$) with the formula
\[
\widetilde{\varrho}_{\varepsilon}^{\,p} := \frac{N(ae\lambda)^{p-1}}{q\{1 - [a(e-1)]^{1/m} \}}
\]
and we show that for all $\widetilde{\varrho}_1 \geq \widetilde{\varrho}_{\varepsilon}$, $w(\cdot)$ is a sub-solution of the equation in \eqref{eq:ELLIPTICCOMPARISONSUBSOLUTIONASYMPTOTICBEHAVIOURPGREATER2}. Note that, since $w(\cdot)$ is radially decreasing, it is sufficient to consider our operator (the $p$-Laplacian) in radial coordinates and verify that for all $\widetilde{\varrho}_1 \geq \widetilde{\varrho}_{\varepsilon}$ and $0 \leq r \leq \widetilde{\varrho}_1$, it holds
\begin{equation}\label{eq:CONDITIONOFSUNSOLUTIONELLIPTIC}
-\Delta_{p,r} w^m := -r^{1-N}\partial_r(r^{N-1}|\partial_r w^m|^{p-2}\partial_r w^m) \leq q(1-w).
\end{equation}
A direct computation shows that
\[
-\Delta_{p,r} w^m = \bigg(\frac{a\lambda}{\widetilde{\varrho}_1^{\,\lambda}}\bigg)^{p-1}\big[N - p(r/\widetilde{\varrho}_1)^{\lambda}\big]e^{(p-1)g(r)}
\]
and that sufficient condition so that inequality \eqref{eq:CONDITIONOFSUNSOLUTIONELLIPTIC} is satisfied is
\[
N(ae\lambda)^{p-1}\widetilde{\varrho}_1^{\,-p} \leq q\{1 - [a(e-1)]^{1/m} \},
\]
which is equivalent to say $\widetilde{\varrho}_1 \geq \widetilde{\varrho}_{\varepsilon}$. Finally, we define $\widetilde{a}_{\varepsilon}$ with the formula:
\[
\widetilde{a}_{\varepsilon}^{\,\lambda} := 1 - \log\Bigg[\frac{a + (1 - \varepsilon/2)^m}{a} \Bigg].
\]
It is not difficult to see that $w(r) \geq 1 - \varepsilon/2$ in $\{r \leq \widetilde{a}_{\varepsilon}\widetilde{\varrho}_1\}$ and that again our assumption on $a$ guarantee the well definition of $\widetilde{a}_{\varepsilon}$.

Hence, if we suppose $\widetilde{\varrho}_1 \geq \widetilde{\varrho}_{\varepsilon}$, we can apply the elliptic Maximum Principle (recall that $w = 0$ on the boundary $\partial\Omega$) deducing
\[
v_{\infty}(x) \geq 1 - \varepsilon/2 \quad \text{in } \{|x| \leq \widetilde{a}_{\varepsilon}\widetilde{\varrho}_1\}.
\]
So, since $v(x,t) \to v_{\infty}(x)$ as $t \to \infty$, a similar inequality holds for the function $v = v(x,t)$ and large times:
\[
v(x,t) \geq 1 - \varepsilon \quad \text{in } \{|x| \leq \widetilde{a}_{\varepsilon}\widetilde{\varrho}_1\}, \text{ for all } t \geq t_1
\]
where $t_1 > 0$ is chosen large enough (depending on $\varepsilon > 0$) and the same conclusion is true for the solution $u = u(x,t)$ of problem \eqref{eq:REACTIONDIFFUSIONEQUATIONPLAPLACIAN} since it holds $v(x,t) \leq u(x,t)$ in $\Omega\times[t_1,\infty)$.

\noindent As the reader can easily check, the proof when $\gamma = 0$ is identical to the case $\gamma > 0$ except for the fact the we apply Proposition \ref{LEMMANONDISCRETEVERSIONLEVELSETSPSEUDO} instead of Proposition \ref{LEMMANONDISCRETEVERSIONLEVELSETS} (see also the following remark).   $\Box$
%
%
%
%
%
%
\paragraph{Remark.} At the beginning of \emph{Step1} (case $\gamma > 0$), we have applied Proposition \ref{LEMMANONDISCRETEVERSIONLEVELSETS}, even though the assumptions \eqref{eq:ASSUMPTIONSONTHEINITIALDATUM} on $u_0$ are not sufficient to guarantee its hypotheses. However, it is really simple to see that for an initial datum $u_0(\cdot)$ satisfying \eqref{eq:ASSUMPTIONSONTHEINITIALDATUM} and $u_0(0) = \max_{x \in\RR} u_0(x)$ there exist $\widetilde{\varrho}_0 > 0$ and $\widetilde{\varepsilon} > 0$ such that the function $\widetilde{u}_0 = \widetilde{u}_0 (x)$ defined in \eqref{eq:INITIALDATUMTESTLEVELSETS} satisfies $\widetilde{u}_0(x) \leq u_0(x)$ in $\RR^N$. Hence, using the Maximum Principle, it follows that the solution $u(x,t)$ of the problem \eqref{eq:REACTIONDIFFUSIONEQUATIONPLAPLACIAN} with initial datum \eqref{eq:ASSUMPTIONSONTHEINITIALDATUM} is greater than the solution of the problem \eqref{eq:REACTIONDIFFUSIONEQUATIONPLAPLACIAN} with initial datum \eqref{eq:INITIALDATUMTESTLEVELSETS}. Consequently, we deduce that $u(x,t)$ satisfies the assertion of Proposition \ref{LEMMANONDISCRETEVERSIONLEVELSETS} applying again the Maximum Principle.

\noindent If $\gamma = 0$, we can proceed similarly. However, this time, we have to start applying Lemma \ref{LEMMAPLACINGBARENBLATTUNDERSOLUTIONPSEUDO} which guarantees the existence of a time $t_1' > 0$ large enough so that we can place an initial datum with form \eqref{eq:INITIALDATUMTESTLEVELSETSPSEUDO} under the solution $u = u(x,t)$ at the time $t = t_1'$. Then, we can apply Proposition \ref{LEMMANONDISCRETEVERSIONLEVELSETSPSEUDO} and proceed with the proof of Lemma \ref{LEMMACOMPACTCONVERGENCETO1ASYMPTOTICBEHAVIOURPGREATER2}.

\noindent We ask the reader to note that, in both cases, the key point consists in deducing that for any ball $\Omega \subset \RR^N$ there exists a time $t_0 > 0$ such that
\[
u(x,t) \geq \widetilde{\varepsilon} \quad \text{in } \Omega\times[t_0,\infty).
\]
If $\gamma > 0$, we get the previous relation noting that we can always find a function satisfying  \eqref{eq:INITIALDATUMTESTLEVELSETS} which can be placed under an initial datum satisfying  \eqref{eq:ASSUMPTIONSONTHEINITIALDATUM} and applying Proposition \ref{LEMMANONDISCRETEVERSIONLEVELSETS}, while, when $\gamma = 0$, we can repeat this procedure but we need both Lemma \ref{LEMMAPLACINGBARENBLATTUNDERSOLUTIONPSEUDO} and Proposition \ref{LEMMANONDISCRETEVERSIONLEVELSETSPSEUDO}.
%
%
%
%
%
%
%
%
%
%
\section{\texorpdfstring{\boldmath}{} Asymptotic behaviour. The case N = 1}\label{SECTIONASYMPTOTICBEHAVIOURDIMENSION1}
Now we address our attention to the study of the asymptotic behaviour of the solutions of problem \eqref{eq:REACTIONDIFFUSIONEQUATIONPLAPLACIAN}, \eqref{eq:ASSUMPTIONSONTHEINITIALDATUM} when the dimension $N = 1$. This represents an important step for two reasons. First of all, we are going to understand the importance of the TWs we found in Theorem \ref{THEOREMEXISTENCEOFTWS} and Theorem \ref{THEOREMEXISTENCEOFTWS1}, since we employ them as sub-solutions and super-solutions of the general solution of problem \eqref{eq:REACTIONDIFFUSIONEQUATIONPLAPLACIAN}, \eqref{eq:ASSUMPTIONSONTHEINITIALDATUM}. Secondly, we will see that the one-dimensional solutions plays an important role in the study of higher-dimension solutions (see Theorem \ref{NTHEOREMCONVERGENCEINNEROUTERSETS}).
\begin{prop}\label{THEOREMCONVERGENCEINNEROUTERSETS}
Let $m > 0$ and $p>1$ such that $\gamma \geq 0$ and fix $N=1$. Then:

(i) For all $0 < c < c_{\ast}$, the solution $u = u(x,t)$ of problem \eqref{eq:REACTIONDIFFUSIONEQUATIONPLAPLACIAN} with initial datum \eqref{eq:ASSUMPTIONSONTHEINITIALDATUM} satisfies
\[
u(x,t) \to 1 \text{ uniformly in } \{|x| \leq ct\} \;\text{ as } t \to \infty,
\]
(ii) Moreover, for all $c > c_{\ast}$ it satisfies,
\[
u(x,t) \to 0 \text{ uniformly in } \{|x| \geq ct\} \;\text{ as } t \to \infty.
\]
\end{prop}
\emph{Proof.} (i) Let $\gamma \geq 0$. We prove that for all $\varepsilon > 0$ and for all $0 <c < c_{\ast}$, there exists $t_2 > 0$ such that
\[
u(x,t) \geq 1 - \varepsilon \quad \text{in } \{|x| \leq ct\} \;\text{ for all } t \geq t_2.
\]
First of all, fix $0 < c < c_{\ast}$ arbitrarily and $\varepsilon > 0$ such that $c+\varepsilon < c_{\ast}$. Now, consider a \emph{CS-TW solution of type 2}: $\varphi = \varphi(\xi)$, where $\xi = x - (c+\varepsilon)t$ (here we are considering the ``reflected'' TWs). As we explained at the end of Section \ref{SECTIONEXISTENCEOFTWS} and Section \ref{CLASSIFICATIONEXISTENCETW}, up to a translation of the $\xi$-axis, we can assume that  $\max_{\xi \in \RR}\varphi(\xi) = \varphi(0) = 1 - \varepsilon$ and $\varphi(\xi_1) = 0 = \varphi(\xi_0)$ for some $\xi_1 < 0 < \xi_0$ (we change notation for the zeroes of $\underline{\varphi}(\cdot)$ since we are working with the ``reflected'' TWs).

\noindent We define the function
\[
\underline{\varphi}(\xi) =
\begin{cases}
\begin{aligned}
1-\varepsilon \quad &\text{if }  \xi \leq 0 \\
\varphi(\xi) \quad &\text{if } 0 \leq \xi \leq \xi_0 \\
0         \;\qquad &\text{otherwise}.
\end{aligned}
\end{cases}
\]
Computing the derivative with respect to the variable $\xi = x -(c+\varepsilon)t$, it not difficult to verify that
\begin{equation}\label{eq:1PROFILEEQUATIONPSIASYMPTOTICBEHAVIOURPGREATER2}
-c\underline{\varphi}' - \varepsilon\underline{\varphi}' \leq [|(\underline{\varphi}^m)'|^{p-2}(\underline{\varphi}^m)']' + f(\underline{\varphi}), \quad \text{for all } \xi \in \RR.
\end{equation}
We proceed proving that $\underline{u}(x,t) = \underline{\varphi}(x -ct)$ is a sub-solution of $u = u(x,t)$ in $\RR_+\times[0,\infty)$, i.e., we verify that
\[
\partial_t\underline{u} \leq \partial_x(|\partial_x\underline{u}^m|^{p-2}\partial_x\underline{u}^m) + f(\underline{u}) \quad \text{in } \RR_+\times[0,\infty),
\]
where $\RR_+$ is the set of the positive real numbers. Imposing this condition and using \eqref{eq:PROFILEEQUATIONPSIASYMPTOTICBEHAVIOURPGREATER2}, it is simple to obtain that a sufficient condition so that the previous inequality is satisfied is $\varphi' \leq 0$ for all $\xi \in \RR$, which is true by construction (note that in this case $\xi = x - ct$, but we do not introduce other variables to avoid weighting down our presentation).

\noindent Now, we fix $\widetilde{\varrho}_1 \geq \widetilde{\varrho}_{\varepsilon}$ large enough such that $\xi_0 \leq \widetilde{a}_{\varepsilon}\widetilde{\varrho}_1$ where $\widetilde{a}_{\varepsilon}$ and $\widetilde{\varrho}_{\varepsilon}$ are the values found in Lemma \ref{LEMMACOMPACTCONVERGENCETO1ASYMPTOTICBEHAVIOURPGREATER2} (of course they refer to a general solution $u = u(x,t)$ of problem \eqref{eq:REACTIONDIFFUSIONEQUATIONPLAPLACIAN} with initial datum \eqref{eq:ASSUMPTIONSONTHEINITIALDATUM}). The function $\underline{u} = \underline{u}(x,t)$ satisfies $u(x,t_1) \geq \underline{u}(x,0)$ in $\RR_+$ thanks to Lemma \ref{LEMMACOMPACTCONVERGENCETO1ASYMPTOTICBEHAVIOURPGREATER2} and $u(0,t_1+t) \geq 1-\varepsilon \geq \underline{u}(0,t)$ for all $t\geq0$ (this follows from the construction of $\underline{u}$). Hence, we obtain
\[
u(x,t_1 + t) \geq \underline{u}(x,t) \quad \text{in } \RR_+\times[0,\infty).
\]
In particular, we deduce $u(x,t_1 +t) \geq \underline{u}(x,t) \geq 1 - \varepsilon$ for all $0 < x \leq ct$ and for all $t \geq 0$. We ask the reader to note that we can conclude that $u(x,t_1 +t) \geq 1-\varepsilon$ for all $-ct < x < 0$ and for all $t \geq0$, simply constructing a sub-solution in the set $\RR_{-}\times [0,\infty)$  ($\RR_{-}$ denotes the set of non-positive real numbers), considering from the beginning the usual ``non-reflected'' CS-TWs of type 2 and proceeding similarly as we did previously. Hence, we can assume
\[
u(x,t_1 + t) \geq 1-\varepsilon \quad \text{in } \{|x| \leq ct\} \quad \text{for all } t \geq 0.
\]
Finally, fix $0 < \widetilde{c} < c$ and $t_2:= ct_1/(c - \widetilde{c})$. Then it holds $u(x,t) \geq 1 - \varepsilon$ in $\{|x| \leq \widetilde{c}t \}$ for all $t \geq t_2$, and the thesis follows from the arbitrariness of $0< \widetilde{c} < c$ and $0 < c < c_{\ast}$.

(ii) We begin with the case $\gamma > 0$. We construct a super-solution which is identically zero on the set $\{x \geq ct\}$ for all $c > c_{\ast}$ and $t$ sufficiently large (note that the same construction can be repeated by ``reflection'' in the set $\{x \leq - ct\}$).

\noindent \emph{Case $u_0(0) = \max_{x\in\RR}u_0(x) < 1$.} Set $\xi = x - c_{\ast}t$ and consider the function
\[
\overline{u}(x,t) := \varphi(\xi) \quad \text{in } \RR\times[0,\infty),
\]
where now $\varphi = \varphi(\xi)$ is the profile of the ``reflected'' \emph{finite} TW found in Theorem \ref{THEOREMEXISTENCEOFTWS}. We showed that there exists $-\infty < \xi_0 < +\infty$ such that $\varphi(\xi) = 0$ for all $\xi \leq \xi_0$, $\varphi' \leq 0$ and $\varphi(-\infty) = 1$. Hence, up to a translation of the $\xi$-axis we can suppose $\overline{u}(x,0) \geq u_0(x)$ for all $x \in \RR^N$. Consequently, applying the Maximum Principle we deduce
\[
\overline{u}(x,t) \geq u(x,t) \quad \text{in } \RR\times[0,\infty).
\]
Let $x_0$ be the \emph{free boundary} point of $\varphi = \varphi(x)$, i.e., $x_0 := \min\{x > 0 : \varphi(x) = 0\}$. Then, it is simple to deduce that
\[
\overline{u}(x,t) \not= 0 \quad \text{ if and only if } \quad  x \leq x_0 + c_{\ast}t.
\]
Now, fix an arbitrary $c > c_{\ast}$ and define $t_2 := x_0/ (c_{\ast} - c)$. Then for all $t \geq t_2$ we have $ct \geq x_0 + c_{\ast}t$ and so
\[
\overline{u}(x,t) \equiv 0 \quad \text{in } \{ x \geq ct\} \quad \text{for all } t \geq t_2
\]
which implies $u(x,t) \equiv 0$ in $\{ x \geq ct\}$ for all $t \geq t_2$ and, since $c > c_{\ast}$ was taken arbitrarily we get assertion (ii) in the case $u_0(0) < 1$.

\noindent Now suppose $u_0(0) = 1$. In this case, the initial datum cannot be ``placed'' under an admissible TW. Then, fix $\delta > 0$ and consider the function
\[
w(y,s) = (1+\delta)u(x,t) \quad \text{ where }\; s(t) = (1+\delta)t \;\text{ and }\; y(x) = (1+\delta)^{m(p-1)/p}x
\]
which solves the equation
\begin{equation}\label{eq:RESCALEDEQUATIONONEDIMENSION}
\partial_s w = \partial_y(|\partial_yw^m|^{p-2}\partial_yw^m) + f((1+\delta)^{-1}w).
\end{equation}
Thus, Theorem \ref{THEOREMEXISTENCEOFTWS} assures that equation \eqref{eq:RESCALEDEQUATIONONEDIMENSION} possesses a finite TW solution $\varphi = \varphi(\xi)$, with $\varphi(\xi) = 0$ for all $\xi \geq \xi_0$, $\varphi' \leq 0$ and $\varphi(-\infty) = 1 + \delta$ (note that we have only re-scaled the solution). Note that in this case both the moving coordinate $\xi = x - c_{\ast}t$ and $c_{\ast} = c_{\ast}(m,p,\delta)$ depend on $\delta > 0$. Now, since $\varphi(-\infty) = 1 + \delta$, we can suppose $u_0(x) \leq \varphi(x)$ for all $x \in \RR$ and we can repeat the same analysis we did in the case $u_0(0) < 1$. Finally, since  $c_{\ast}(m,p,\delta) \to c_{\ast}(m,p)$ as $\delta \to 0$ (this follows from the fact that the the system of ODEs of the re-scaled equation converges to system \eqref{eq:SYSTEMNONSINGULARTWS} which is derived from the standard equation), we get the assertion (ii) in the case $\gamma > 0$.

(iii) Case $\gamma = 0$. Again we suppose $u_0(0) < 1$ (the case $u_0(0) = 1$ can be treated as we did previously). Fix $\varepsilon > 0$ and let $\varphi = \varphi(\xi)$ be the profile of the positive TW with critical speed $c_{\ast}(m,p)$ found in Theorem \ref{THEOREMEXISTENCEOFTWS1}, where $\xi = x - c_{\ast}t$. We have $\varphi > 0$, $\varphi' < 0$ and $ \varphi(-\infty) = 1$, $\varphi(\infty) = 0$ and, in particular, $\varphi(\xi) < \varepsilon$ for all $\xi \geq \xi_{\varepsilon}$, where $\xi_{\varepsilon}$ is chosen large enough depending on $\varepsilon > 0$. Now,
define $\overline{u}(x,t) = \varphi(\xi)$ and note that we can suppose $\underline{u}(x,0) \geq u_0(x)$ in $\RR$. Consequently, it follows $\underline{u}(x,t) \geq u(x,t)$ in $\RR\times[0,\infty)$ and, furthermore, we obtain
\[
u(x,t) \leq \varepsilon \quad \text{ in } \{x \geq \xi_{\varepsilon} + c_{\ast}t\} \quad \text{for all } t \geq 0.
\]
Now, for all fixed $c > c_{\ast}$ and $t \geq t_2 := \xi_{\varepsilon}/(c - c_{\ast})$ we have $ct \geq \xi_{\varepsilon} + c_{\ast}t$ and so it follows
\[
u(x,t) \leq \varepsilon \quad \text{ in } \{ x \geq ct\} \quad \text{for all } t \geq t_2
\]
i.e., the thesis. Since the same procedure can be repeated for the ``reflected'' TWs and $c > c_{\ast}$ is arbitrary, we end the proof of (ii) case $\gamma = 0$. $\Box$
\paragraph{Remark.} As previously mentioned, when $\gamma > 0$ (``slow diffusion'' assumption) the general solutions of problem \eqref{eq:REACTIONDIFFUSIONEQUATIONPLAPLACIAN} exhibit free boundaries. This fact follows from the proof of the previous theorem (part (ii) case $\gamma > 0$). Indeed, we showed that the solution is identically zero when $c > c_{\ast}$ in the outer set $\{|x| \geq ct\}$ as $t \to \infty$. This fact represents the significant difference respect to the case $\gamma = 0$ (``pseudo-linear'' assumption) in which the general solutions are positive everywhere. Hence, we can conclude that for all $\gamma \geq 0$, the general solutions of problem \eqref{eq:REACTIONDIFFUSIONEQUATIONPLAPLACIAN} expand linearly in time (for large times) with a critical speed $c_{\ast} > 0$ but, in the case $\gamma > 0$, for all fixed time, they are identically zero outside a ball with radius large enough, whilst, when $\gamma = 0$, they are positive everywhere. This is true when $N = 1$ and, in the next section, we will see that it is possible to extend the previous assertions for all $N \geq 1$.
%
%
%
%
%
%
%
%
\section{Asymptotic behaviour with radial waves in several dimensions}\label{SECTIONASYMPTOTICBEHAVIOUR2}
In this section, we focus on the case $N \geq 2$ and we establish our main PDE result (see Theorem \ref{NTHEOREMCONVERGENCEINNEROUTERSETS}). We are going to show that for all $0 < c < c_{\ast}$ ($c_{\ast} := c_{\ast}(m,p)$ is the critical value found in Theorem \ref{THEOREMEXISTENCEOFTWS}), the solution of problem \eqref{eq:REACTIONDIFFUSIONEQUATIONPLAPLACIAN} with initial datum \eqref{eq:ASSUMPTIONSONTHEINITIALDATUM} converges uniformly to the equilibrium point $u = 1$ in the inner set $\{|x| \leq ct \}$ as $t \to \infty$, while, for all $c > c_{\ast}$ to the trivial solution $u = 0$ in the outer set $\{|x| \geq ct \}$ as $t \to \infty$. More precisely, we will prove that, when $c > c_{\ast}$, the solution of problem \eqref{eq:REACTIONDIFFUSIONEQUATIONPLAPLACIAN} with initial datum \eqref{eq:ASSUMPTIONSONTHEINITIALDATUM} is identically zero in the set $\{|x| \geq ct \}$ as $t \to \infty$, i.e., the general solution has a \emph{free boundary}. So, in some sense we obtain a nonlinear version of Theorem \ref{THMASYMPTOTICSTLARGELINEARCASEINTRO} (case $m = 1$ and $p = 2$) and again we deduce that the steady states $u = 1$ is asymptotically stable while the null solution is unstable and the density $u = u(x,t)$ saturates the all the free space with velocity $c_{\ast}$. However, as we explained in the introduction, it is difficult to obtain information on the solution $u(x,t)$ on the mobile coordinate $x = \xi - c_{\ast}t$, $\xi \in \RR$ and describe the precise shape of the solution.

Before proceeding, we point out that it is sufficient to study the asymptotic behaviour of radial solutions of problem \eqref{eq:REACTIONDIFFUSIONEQUATIONPLAPLACIAN} with radial non-increasing initial datum satisfying \eqref{eq:ASSUMPTIONSONTHEINITIALDATUM}. Indeed, we only need to note that for every initial datum $u_0(x)$ satisfying \eqref{eq:ASSUMPTIONSONTHEINITIALDATUM}, it is possible to construct (non-increasing) initial radial ``sub-datum'' and ``super-datum'' which generate radial solutions that act as barriers from below and from above for the general solution generated by $u_0(x)$. Consequently, the asymptotic behaviour of general solutions is completely described by the one of radial solutions.

So we consider solutions of the problem
\begin{equation}\label{eq:KPPRADIALASYMPTOTICBEHAVIOURPGREATER2}
\begin{cases}
\begin{aligned}
\partial_tu = \Delta_{p,r} u^m + f(u) \quad &\text{in } \RR_+\times(0,\infty)\\
u(r,0) = u_0(r) \;\;\;\quad\qquad &\text{in } \RR_+,
\end{aligned}
\end{cases}
\end{equation}
where $\RR_+$ is the set of positive real numbers, $r = |x|$ is the radial coordinate ($x \in \RR^N$) and
\[
\Delta_{p,r} u^m := r^{1-N}\partial_r(r^{N-1}|\partial_r u^m|^{p-2}\partial_r u^m) = \partial_r(|\partial_r u^m|^{p-2}\partial_r u^m) + (N-1)r^{-1}|\partial_r u^m|^{p-2}\partial_r u^m
\]
is the ``radial doubly nonlinear'' operator. Finally, we assume that the initial datum $u_0 \in \mathcal{C}_c(\RR^N)$ is non-trivial, $0 \leq u_0 \leq 1$ and radial with $\partial_ru_0 \leq 0$ (again we assume $u_0(0) = \max_{x \in \RR^n}u_0(x)$).

The main idea consists in using the solutions of problem \eqref{eq:REACTIONDIFFUSIONEQUATIONPLAPLACIAN} posed in dimension one as sub-solutions and super-solutions of the radial solutions of problem \eqref{eq:KPPRADIALASYMPTOTICBEHAVIOURPGREATER2} and deducing their asymptotic behaviour applying the results found in Proposition \ref{THEOREMCONVERGENCEINNEROUTERSETS}.
\paragraph{Proof of Theorem \ref{NTHEOREMCONVERGENCEINNEROUTERSETS}.} We prove separately part (i) and (ii). In (i), we firstly consider the case $\gamma > 0$ and then $\gamma = 0$.

(i) Consider $\gamma > 0$, fix $0 < c < c_{\ast}$ and take $0 < \varepsilon < 1$ such that $c +\varepsilon < c_{\ast}$ (this choice does not represent a problem since $\varepsilon > 0$ will be taken arbitrary small). First of all, we construct a sub-solution for problem \eqref{eq:KPPRADIALASYMPTOTICBEHAVIOURPGREATER2} in a set of the form $\RR_{+}\times[t_{\varepsilon},\infty)$ for $t_{\varepsilon} > 0$ large enough.

\noindent As we showed in Theorem \ref{THEOREMEXISTENCEOFTWS}, the one-dimensional equation
\[
\partial_tv = \partial_y(|\partial_yv^m|^{p-2}\partial_yv^m) + f(v) \quad \text{in } \RR\times[0,\infty)
\]
admits \emph{CS-TWs of type 2} in the form $v(y,t) = \varphi(y - (c+\varepsilon)t)$, $y \in \RR$ and $t \geq 0$, with speed $c+\varepsilon$. In particular, we are considering the TW solution moving to the right direction with $\varphi(0) = \max_{\xi\in\RR}\varphi(\xi) = 1 - \varepsilon$, $\varphi(\xi_1) = 0 = \varphi(\xi_0)$ for some $\xi_1 < 0 < \xi_0$ (note that we change the notation respect to Section \ref{SECTIONEXISTENCEOFTWS} since here we are referring to the ``reflected'' TW). Now, we define
\[
\underline{\varphi}(\xi) =
\begin{cases}
\begin{aligned}
1-\varepsilon \quad &\text{if }  \xi \leq 0 \\
\varphi(\xi) \quad &\text{if } 0 \leq \xi \leq \xi_0 \\
0         \;\qquad &\text{otherwise}
\end{aligned}
\end{cases}
\]
and we observe that the profile $\underline{\varphi}(\cdot)$ satisfies the differential inequality
\begin{equation}\label{eq:PROFILEEQUATIONPSIASYMPTOTICBEHAVIOURPGREATER2}
-c\underline{\varphi}' - \varepsilon\underline{\varphi}' \leq [|(\underline{\varphi}^m)'|^{p-2}(\underline{\varphi}^m)']' + f(\underline{\varphi}),\quad \text{for all } \xi \in \RR
\end{equation}
where we computed the derivative respect to the variable $\xi = y -(c+\varepsilon)t$, i.e., $\underline{\varphi}'(\xi) = d\underline{\varphi}(\xi)/d\xi$.
\paragraph{Case $\boldsymbol{m \geq 1}$ and $\boldsymbol{p > 2}$}. Now, we use \eqref{eq:PROFILEEQUATIONPSIASYMPTOTICBEHAVIOURPGREATER2} to show that the function $\underline{u}(r,t) = \underline{\varphi}(r - ct)$, where $r = |x|$, $x \in \RR^N$, is a sub-solution for problem \eqref{eq:KPPRADIALASYMPTOTICBEHAVIOURPGREATER2} if $t$ is large enough, i.e., it satisfies
\begin{equation}\label{eq:SUBSOLUTIONEQUATIONNLARGERTHAN2}
\partial_t\underline{u} \leq \Delta_{p,r}\underline{u}^m + f(\underline{u}), \quad \text{for } t \gg 0.
\end{equation}
First of all, it is simple to see that the previous inequality is trivially satisfied when $\underline{u}(\cdot)$ is constant. Hence, focusing on the interval $0 < \xi < \xi_0$ and using \eqref{eq:PROFILEEQUATIONPSIASYMPTOTICBEHAVIOURPGREATER2}, we deduce that a sufficient condition so that the previous inequality is satisfied is
\[
\varepsilon \underline{\varphi}' \leq m^{p-1}(N-1)(\xi+ct)^{-1} \underline{\varphi}^{\mu}|\underline{\varphi}'|^{p-2}\underline{\varphi}',
\]
which, since the derivative of $\varphi(\cdot)$ is negative, can be rewritten as
\begin{equation}\label{eq:1SUFFICIENTCONDITIONONPSI}
\xi + ct \geq \frac{m^{p-1}(N-1)}{\varepsilon}\underline{\varphi}^{\mu}|\underline{\varphi}'|^{p-2}.
\end{equation}
Note that we computed the derivative with respect to the variable $\xi = r - ct$. We do not change notation respect to the computation carried out in \eqref{eq:PROFILEEQUATIONPSIASYMPTOTICBEHAVIOURPGREATER2} in order to avoid the introduction of a large number and useless variables.

\noindent \emph{Case $\xi \sim 0$.} Now, note that, since $\underline{\varphi} \sim 1-\varepsilon$ and $|\underline{\varphi}'| \sim 0$ for $\xi \sim 0$ and $p > 2$, it is possible to conclude that the previous inequality is satisfied for $\xi \sim 0$ and for all $t \geq t_{1\varepsilon}' > 0$.

\noindent \emph{Case $\xi \sim \xi_0$.} We begin supposing $m > 1$. Since we showed $|\varphi'| \sim a\varphi^{1-m}$ for $\varphi \sim 0$ (see \eqref{eq:ASYMPTOTICBEHAVIOURDERIVATIVEXZERO}) for some $a > 0$, we deduce $\underline{\varphi}^{\mu}|\underline{\varphi}'|^{p-2} \sim a^{p-2}\varphi^{m-1} \sim 0$ for $\xi \sim \xi_0$ and, consequently, we can conclude that for all $t \geq 0$, \eqref{eq:1SUFFICIENTCONDITIONONPSI} is satisfied when $\xi \sim \xi_0$. If $m = 1$, we have $\underline{\varphi}^{\mu}|\underline{\varphi}'|^{p-2} \sim a^{p-2}$ and inequality \eqref{eq:1SUFFICIENTCONDITIONONPSI} is satisfied if
\[
t \geq t_{1\varepsilon}'':= [a^{p-2}m^{p-1}(N-1)/\varepsilon - \xi_0]/c.
\]

\noindent \emph{Case $0 < \xi < \xi_0$ with $\xi \not\sim 0$ and $\xi \not\sim \xi_0$.} In this case, inequality \eqref{eq:1SUFFICIENTCONDITIONONPSI} holds for all $t \geq t_{1\varepsilon}'''$, where $t_{1\varepsilon}'''$ is chosen large enough, since both $|\underline{\varphi}|$ and $|\underline{\varphi}'|$ are bounded from above and below by positive constants. In particular, a sufficient condition is given by
\[
t \geq t_{1\varepsilon}''' := \max\big(\underline{\varphi}^{\mu}|\underline{\varphi}'|^{p-2}\big)m^{p-1}(N-1)/(\varepsilon c)
\]
where the maximum is chosen is calculated in the ``set'' $0 < \xi < \xi_0$ with $\xi \not\sim 0$ and $\xi \not\sim \xi_0$. Consequently, this analysis allows us to state that, when $m\geq1$ and $p>2$, we have that $\underline{u} = \underline{u}(r,t)$ satisfies \eqref{eq:SUBSOLUTIONEQUATIONNLARGERTHAN2} for all $t \geq t_{1\varepsilon}:=\max\{t_{1\varepsilon}',t_{1\varepsilon}'',t_{1\varepsilon}'''\}$.
\paragraph{Case $\boldsymbol{m>1}$ and $\boldsymbol{1 < p \leq 2}$.} Again, we want to find a function $\underline{u} = \underline{u}(r,t)$ such that \eqref{eq:SUBSOLUTIONEQUATIONNLARGERTHAN2} is satisfied. This time our candidate is $\underline{u}(r,t) = \underline{\varphi}(\delta^{1/p} r - c\delta t)$, with $\delta = 1-\varepsilon$ and $\xi = \delta^{1/p} r - c\delta t$. Following the computations done previously, it is not difficult to see that a sufficient condition so that \eqref{eq:SUBSOLUTIONEQUATIONNLARGERTHAN2} is satisfied is
\begin{equation}\label{eq:2SUFFICIENTCONDITIONONPSI}
\varepsilon|\underline{\varphi}'| \geq \frac{b_{\varepsilon}}{\delta^{-1}\xi + ct} \,\underline{\varphi}^{\mu}|\underline{\varphi}'|^{p-1} - \frac{\varepsilon}{1-\varepsilon}f(\underline{\varphi}),
\end{equation}
where $b_{\varepsilon}:= m^{p-1}(N-1)/(1-\varepsilon)$.

\noindent \emph{Case $\xi \sim 0$.} For all $t \geq t_{2\varepsilon}' > 0$, the previous inequality is trivially satisfied when $\xi \sim 0$ since we have $\underline{\varphi} \sim 1 - \varepsilon$ and $\underline{\varphi}' \sim 0$ and $f(1-\varepsilon) > 0$.

\noindent \emph{Case $\xi \sim \xi_0$.} On the other hand, when $\xi \sim \xi_0$, we have $\underline{\varphi} \sim 0$, $|\underline{\varphi}'| \sim a\underline{\varphi}^{1-m} \sim \infty$ (since $m > 1$) and $\underline{\varphi}^{\mu}|\varphi'|^{p-1} \sim a^{p-1}$. Thus we can deduce that \eqref{eq:2SUFFICIENTCONDITIONONPSI} is satisfied for all $t \geq 0$ (of course when $\xi \sim \xi_0$).

\noindent \emph{Case $0 < \xi < \xi_0$ with $\xi \not\sim 0$ and $\xi \not\sim \xi_0$.} Finally, when $0 < \xi < \xi_0$ with $\xi \not\sim 0$ and $\xi \not\sim \xi_0$, it is possible to note that the condition
\[
\varepsilon ct \geq b_{\varepsilon}\underline{\varphi}^{\mu}|\underline{\varphi}'|^{p-2}
\]
is sufficient to guarantee \eqref{eq:2SUFFICIENTCONDITIONONPSI} and so, since with the current assumptions on $\xi$ we have that $|\varphi'|$ is bounded from above and below, we deduce the existence of a value $t_{2\varepsilon}''$ such that \eqref{eq:SUBSOLUTIONEQUATIONNLARGERTHAN2} is satisfied for all $t \geq t_{2\varepsilon}''$ when $0 < \xi < \xi_0$ with $\xi \not\sim 0$ and $\xi \not\sim \xi_0$. Hence, we have that when $m>1$ and $1 < p \leq 2$, \eqref{eq:SUBSOLUTIONEQUATIONNLARGERTHAN2} is satisfied for all $t \geq t_{2\varepsilon} := \max\{t_{2\varepsilon}',t_{2\varepsilon}''\}$.
\paragraph{Case $\boldsymbol{0 < m < 1}$ and $\boldsymbol{p > 2}$.} This is the most delicate case. We define $\underline{u}(r,t) = \underline{\varphi}(r - c\delta t)$, with $\delta = (1+\varepsilon^n)^{-1}$ and we impose condition \eqref{eq:SUBSOLUTIONEQUATIONNLARGERTHAN2} (the value of $n \in \NN$ is not important now and will be specified later). Using inequality \eqref{eq:PROFILEEQUATIONPSIASYMPTOTICBEHAVIOURPGREATER2} and carrying out some tedious computations, we obtain that a sufficient condition so that \eqref{eq:SUBSOLUTIONEQUATIONNLARGERTHAN2} is satisfied is
\begin{equation}\label{eq:3SUFFICIENTCONDITIONONPSI}
\varepsilon^{1-n}m^{1-p}|\underline{\varphi}'| \geq |\mu|\underline{\varphi}^{\mu-1}|\underline{\varphi}'|^p -(p-1)\underline{\varphi}^{\mu}|\underline{\varphi}'|^{p-2}\underline{\varphi}'' + b_{\varepsilon}(\xi + c\delta t)^{-1}\underline{\varphi}^{\mu} |\underline{\varphi}'|^{p-1} - m^{1-p}f(\underline{\varphi}),
\end{equation}
where this time $b_{\varepsilon} := (1+\varepsilon^n)(N-1)/\varepsilon^n$.

\noindent\emph{Case $\xi \sim 0$.} If $\xi \sim 0$, we have $\underline{\varphi} \sim 1 - \varepsilon$ and $|\underline{\varphi}'| \sim 0$. Hence, recalling relation \eqref{eq:ESTIMATEONSECONDDERIVATIVENEARMAXIMUMPOINT} and applying it with $X = 1 - \varepsilon$, it is not difficult to see that for all $t \geq t_{3\varepsilon} > 0$, \eqref{eq:3SUFFICIENTCONDITIONONPSI} is equivalent to
\[
m^{2-p}f(1-\varepsilon) \geq m^{2-p}(1-\varepsilon)^{\mu + p - 2 -\gamma}f(1-\varepsilon),
\]
which is satisfied (the equality holds) since $\mu = \gamma + 2 -p$.

\noindent\emph{Case $\xi \sim \xi_0$.} When $\xi \sim \xi_0$, we have again $\underline{\varphi} \sim 0$ and $\underline{\varphi}' \sim a\underline{\varphi}^{1-m}$ and, as the reader can easily check, \eqref{eq:3SUFFICIENTCONDITIONONPSI} is equivalent to
\[
0 \geq a^{p-2}\underline{\varphi}^{-m}\big[a^2|\mu| - (p-1)\underline{\varphi}^{2m-1}\underline{\varphi}'' \big] + a^{p-1}b_{\varepsilon}(\xi_0 + c\delta t)^{-1}.
\]
Using relation \eqref{eq:SECONDESTIMATEONSECONDDERIVATIVENEARMAXIMUMPOINT}, it is simple to re-write the previous inequality as
\[
0 \geq -a^pd_{m,p}\underline{\varphi}^{-m} + b_{\varepsilon}(\xi_0 + c\delta t)^{-1},
\]
where $d_{m,p} := (p-2-\gamma)/m^2 + \mu$. Since $0 < m < 1$, we have $d_{m,p}> 0$ and so, using that $\underline{\varphi}^{-m} \sim \infty$ as $\underline{\varphi} \sim 0$, we obtain that for all fixed $t \geq 0$, the last inequality is satisfied and we conclude the analysis of the case $\xi \sim \xi_0$.

\noindent \emph{Case $0 < \xi < \xi_0$ with $\xi \not\sim 0$ and $\xi \not\sim \xi_0$.} It is not difficult to see that \eqref{eq:3SUFFICIENTCONDITIONONPSI} is equivalent to
\begin{equation}\label{eq:4SUFFICIENTCONDITIONONPSI}
|\underline{\varphi}'| + \varepsilon^{n-1}f(\underline{\varphi}) \geq -\varepsilon^{n-1}(|(\underline{\varphi}^m)'|^{p-2}(\underline{\varphi}^m)')' + \varepsilon^{n-1}b_{\varepsilon}(\xi +c\delta t)^{-1}\underline{\varphi}^{\mu}|\underline{\varphi}'|^{p-1}
\end{equation}
and, furthermore, we have that $\underline{\varphi}$ and $|\underline{\varphi}'|$ are bounded from above and below by positive constants. Moreover, combining \eqref{eq:SYSTEMNONSINGULARTWSXI} and \eqref{eq:SECONDDERIVATIVEOFPROFILEFORMULA}, it is straightforward to see that the ``second order term'' $(|(\underline{\varphi}^m)'|^{p-2}(\underline{\varphi}^m)')'$ has not definite sign but is bounded too (from below and above). Moreover, since $t \geq 0$, we have that a sufficient condition so that \eqref{eq:4SUFFICIENTCONDITIONONPSI} is satisfied is
\begin{equation}\label{eq:5SUFFICIENTCONDITIONONPSI}
|\underline{\varphi}'| + \varepsilon^{n-1}f(\underline{\varphi}) \geq -\varepsilon^{n-1}(|(\underline{\varphi}^m)'|^{p-2}(\underline{\varphi}^m)')' + \varepsilon^{n-1}b_{\varepsilon}\xi^{-1}\underline{\varphi}^{\mu}|\underline{\varphi}'|^{p-1}.
\end{equation}
Hence, we take $n \in \NN$ large enough (independently of $t\geq0$) to make the terms in the right side of the previous inequality smaller and proving the validity of the inequality \eqref{eq:5SUFFICIENTCONDITIONONPSI}. Consequently, we can state that also in this last case, there exists $t_{3\varepsilon} > 0$ large enough such that $\underline{u} = \underline{u}(r,t)$ satisfies \eqref{eq:SUBSOLUTIONEQUATIONNLARGERTHAN2} for all $t \geq t_{3\varepsilon}$.

\noindent Hence, for all $\gamma > 0$, we constructed a sub-solution $\underline{u} = \underline{u}(r,t)$ in the set $\RR_{+}\times[t_{\varepsilon},\infty)$, where $t_{\varepsilon} = \max\{t_{1\varepsilon},t_{2\varepsilon},t_{3\varepsilon}\}$. Now, we proceed with the proof of the case $m > 1$ and $p > 2$. We can treat the other ranges of parameters with identical methods, so we skip them leaving only a short comment at the end of the proof.

\noindent Define $\underline{w}(r,t) := \underline{u}(r,t_{\varepsilon}+t)$ and let $r_0>0$ be the ``free boundary point'' of $\underline{w}(r,0)$, i.e., $r_0 := \min\{r > 0 : \underline{w}(r,0) = 0\}$.
%
%
From the analysis done previously, we have that the sub-solution $\underline{w} = \underline{w}(r,t)$ satisfies
\[
\begin{cases}
\begin{aligned}
\partial_t\underline{w} \leq \Delta_{p,r}\underline{w}^m + f(\underline{w}) \;\qquad &\text{in } \RR_{+}\times[0,\infty) \\
\underline{w}(r,t) = \underline{\varphi}(r - c(t_{\varepsilon}+t))  \quad &\text{in } \{r = 0\}\times[0,\infty) \\
\underline{w}(r,0) = \underline{\varphi}(r-ct_{\varepsilon})      \;\quad\qquad &\text{in } \RR_{+}. \\
\end{aligned}
\end{cases}
\]
Now, let $u = u(r,t)$ be a radial solution of problem \eqref{eq:REACTIONDIFFUSIONEQUATIONPLAPLACIAN} and fix $\widetilde{\varrho}_1 \geq \widetilde{\varrho}_{\varepsilon}$ large enough such that $r_0 \leq \widetilde{a}_{\varepsilon}\widetilde{\varrho}_1$, where $0 < \widetilde{a}_{\varepsilon} < 1$ and $\widetilde{\varrho}_{\varepsilon} > 0$ are the values found in Lemma \ref{LEMMACOMPACTCONVERGENCETO1ASYMPTOTICBEHAVIOURPGREATER2}. Consider the function $w(r,t) = u(r,t_1 + t)$ in $\RR_+\times(0,\infty)$, where $t_1 > 0$ is chosen depending on $\varepsilon > 0$ and $\widetilde{\varrho}_1 > 0$ as in Lemma \ref{LEMMACOMPACTCONVERGENCETO1ASYMPTOTICBEHAVIOURPGREATER2}. Thanks to Lemma \ref{LEMMACOMPACTCONVERGENCETO1ASYMPTOTICBEHAVIOURPGREATER2}, we have that $w(r,0) \geq \underline{w}(r,0)$ for all $r \geq 0$ and, since $1-\varepsilon \geq \underline{w}(r,t)$ by construction, it holds $w(0,t) \geq \underline{w}(0,t)$ for all $t \geq 0$. So, we can apply the Maximum Principle deducing
\[
u(r,t_1' + t) \geq \underline{u}(r,t) \quad \text{in } \RR_+\times[0,\infty).
\]
where $t_1' = t_1 - t_{\varepsilon}$. In particular, we have $u(r,t_1 + t) \geq \underline{u}(r,t) = 1 - \varepsilon$ for all $r \leq ct$ and all $t \geq 0$. Finally, proceeding as in the final lines of the proof of Proposition \ref{THEOREMCONVERGENCEINNEROUTERSETS} (part (i) case $\gamma > 0$) we complete the proof of Theorem \ref{NTHEOREMCONVERGENCEINNEROUTERSETS} part (i) case $\gamma > 0$.

\noindent Note that when $\gamma > 0$ but we are not considering the case $m > 1$ and $p > 2$, the real speed of propagation of the function $\underline{u} = \underline{u}(r,t)$ is not $c$ but $c\delta$. Nevertheless, since $0 < \varepsilon < 1$ was taken arbitrarily and $c\delta \to c$ as $\varepsilon \to 0$, we can repeat the procedure followed in the case $m > 1$ and $p>2$ and, a the end, take the limit as $\varepsilon \to 0$.

(i) Case $\gamma = 0$. The analysis is very similar to the case $\gamma > 0$ except for the fact that we use the TWs found in Section \ref{CLASSIFICATIONEXISTENCETW}. As we explained in the proof of Theorem \ref{THEOREMEXISTENCEOFTWS1}, their properties coincide with the ones of the case $\gamma > 0$ near the the maximum point and near the change sign point and so, it is possible to repeat the proof carried out for the case $\gamma > 0$. We leave the details to the reader.

\medskip

(ii) This assertion is simpler. Indeed, consider a radial solution $u = u(r,t)$ of problem \eqref{eq:REACTIONDIFFUSIONEQUATIONPLAPLACIAN} with $r = |x|$ and $x \in \RR^N$. Since $\partial_ru^m \leq 0$, it is straightforward to deduce that $\Delta_{p,r}u^m \leq \partial_r(|\partial_ru^m|^{p-2}\partial_ru^m)$, i.e., the solution $v(r,t)$ of the problem \eqref{eq:REACTIONDIFFUSIONEQUATIONPLAPLACIAN} posed in dimension 1 is a super-solution of the solution of the same problem posed in dimension $N \geq 1$. So, we have $u(r,t) \leq v(r,t)$ in $\RR_+\times(0,\infty)$ and we conclude the proof applying the assertion (ii) of Proposition \ref{THEOREMCONVERGENCEINNEROUTERSETS}. $\Box$
\paragraph{Remark.} The proof of the previous theorem (part (ii)) shows that the observations made at the end of Section \ref{SECTIONASYMPTOTICBEHAVIOURDIMENSION1} on the existence of free boundaries in the case $\gamma >0$, and their absence when $\gamma = 0$, are still true in several dimensions.
\paragraph{Proof of Corollary \ref{CORASYMBEHAMOREGENDATUM}.}
We divide the proof in four short steps.

\emph{Step1.} In this first step we show assertion (i) of Theorem \ref{NTHEOREMCONVERGENCEINNEROUTERSETS}. We simply note that for all initial data $u_0(\cdot)$ satisfying \eqref{eq:ASSUMPTIONSONTHEINITIALDATUMEXPONENTIALGGGEQ0} and/or \eqref{eq:ASSUMPTIONSONTHEINITIALDATUMEXPONENTIALGG0}, there exists a ``sub-initial datum'' $\underline{u}_0(\cdot)$ satisfying \eqref{eq:ASSUMPTIONSONTHEINITIALDATUM}, i.e. with compact support, such that $\underline{u}_0(x) \leq u_0(x)$ for all $x \in \RR^N$. Hence, the solution $\underline{u}(x,t)$ of problem \eqref{eq:REACTIONDIFFUSIONEQUATIONPLAPLACIAN} with initial datum $\underline{u}_0(\cdot)$ satisfies $\underline{u}(x,t) \leq u(x,t)$ in $\RR^N\times(0,\infty)$ by the Maximum Principle. Here, $u(x,t)$ stands for the solution of problem \eqref{eq:REACTIONDIFFUSIONEQUATIONPLAPLACIAN} with initial data \eqref{eq:ASSUMPTIONSONTHEINITIALDATUMEXPONENTIALGGGEQ0} or \eqref{eq:ASSUMPTIONSONTHEINITIALDATUMEXPONENTIALGG0} depending on $\gamma \geq 0$. Thus, since $\underline{u}(x,t)$ satisfies statement (i) of Theorem \ref{NTHEOREMCONVERGENCEINNEROUTERSETS}, we deduce that $u(x,t)$ has to satisfy it too.

\emph{Step2.} We show assertion (ii) of Theorem \ref{NTHEOREMCONVERGENCEINNEROUTERSETS} when $\gamma > 0$ and $N = 1$. For all $c > c_{\ast}(m,p)$, we can place an admissible (reflected) TW $\varphi(x - ct)$ above an initial datum $u_0(x)$ satisfying \eqref{eq:ASSUMPTIONSONTHEINITIALDATUMEXPONENTIALGGGEQ0}. This is possible thanks to formula \eqref{eq:ASYMBEHASUPERTWGGPOS}, which can be re-written for $t = 0$ as
\[
\varphi(x) \sim a_0 \exp\Big(-c^{-1}f'(0)x\Big) = a_0\exp\Big(-\nu^{-1}f'(0)^{1/p}x\Big), \quad \text{for } x \sim \infty, \quad \text{for some }\, a_0>0.
\]
We recall that $0 < \nu < \nu_{\ast}$ is the speed with $f'(0)=1$. Consequently, proceeding as in the parts (ii) and (iii) of the proof of Proposition \ref{THEOREMCONVERGENCEINNEROUTERSETS} we get our statement and we complete this step.

\emph{Step3.} Take now $\gamma = 0$ and $N = 1$. We fix $c > c_{0\ast}(m,p) = p(m^2f'(0))^{\frac{1}{mp}}$ (see Section \ref{CLASSIFICATIONEXISTENCETW}), and we proceed as in \emph{Step2}. However, in this case, we consider an admissible (reflected) TW $\varphi(x-ct)$ satisfying for $t = 0$
\[
\varphi(x) \sim a_0 \exp\Big(-\frac{\lambda_1}{m}x\Big) \quad \text{for } x \sim \infty, \quad \text{for some } a_0>0,
\]
according to \eqref{eq:ASYMBEHGGGEQ0SUPERTW}. Recall that $\lambda_1 = \lambda_1(c)$ and $\lambda_1 < \lambda_{\ast} = (c_{0\ast}/p)^m = (m^2f'(0))^{1/p}$ for all $c > c_{0\ast}(m,p)$. This fact implies that we can actually place the TW above an initial datum satisfying \eqref{eq:ASSUMPTIONSONTHEINITIALDATUMEXPONENTIALGG0}. Indeed, in \eqref{eq:ASSUMPTIONSONTHEINITIALDATUMEXPONENTIALGG0} we have assumed
\[
u_0(x) \leq a_0\, x^{\frac{2}{p}}\exp\Big(-m^{\frac{2-p}{p}}f'(0)^{\frac{1}{p}}x\Big) = a_0\, x^{\frac{2}{p}}\exp\bigg(-\frac{\lambda_{\ast}}{m}x\bigg) \quad \text{for }\; x \sim \infty.
\]
Once we can employ the TW $\varphi(\cdot)$ as a barrier from above, we conclude the proof of this step with the same procedure carried out before.

\emph{Step4.} We end the proof with some comments. First of all, when $N=1$, we can repeat the previous analysis when $x \sim -\infty$ by using non-reflected TWs with the form $\varphi (x+ct)$ (note that we used $|x|$ in conditions \eqref{eq:ASSUMPTIONSONTHEINITIALDATUMEXPONENTIALGGGEQ0} and \eqref{eq:ASSUMPTIONSONTHEINITIALDATUMEXPONENTIALGG0}). Secondly, in order to prove our statements for $N \geq 2$, it is sufficient to repeat the arguments of the proof of part (ii) of Theorem \ref{NTHEOREMCONVERGENCEINNEROUTERSETS} in which we have studied the asymptotic behaviour of non-increasing radial solutions. Finally, we stress that if the initial datum satisfies $\max_{x\in\RR^N}u_0(x) = 1$, we can proceed as in the proof of part (ii) of Proposition \ref{THEOREMCONVERGENCEINNEROUTERSETS}.
$\Box$
%
%
%
%
%
%
%
%
%
%
\section{Models with ``strong'' reaction}\label{SECTIONSTRONGREACTION}
In order to give to reader a wider vision of the work we have carried out, we focus on a model with a ``strong'' reaction term. Following \cite{DePablo-Vazquez1:art} and \cite{DePablo-Vazquez2:art}, we shortly present an extension of Theorem \ref{THEOREMEXISTENCEOFTWS} and Theorem \ref{THEOREMEXISTENCEOFTWS1}, studying the existence of TWs for the equation
\begin{equation}\label{eq:SLOWDIFFSTRONGREACTPLAPLACIAN}
\partial_tu = \partial_x(|\partial_xu^m|^{p-2}\partial_xu^m) + u^n(1 - u) \quad \text{in } \RR\times(0,\infty),
\end{equation}
where $m > 0$ and $p > 1$ such that $\gamma \geq 0$ and $n \in \RR$ (note that we get the Fisher-KPP equation with doubly nonlinear diffusion choosing $n = 1$). As we explained in Section \ref{SUBSECTIONNONLINEARDIFFUSION} for the Porous Medium case, the problem consists in understanding if equation \eqref{eq:SLOWDIFFSTRONGREACTPLAPLACIAN} has admissible TWs (in the sense of definition \eqref{eq:CONDITIONONPHIADMISSIBLETWINTRO}) when we replace the usual smooth reaction term $f(u) = u(1-u)$ with a non-smooth one, called ``strong'' reaction for the ``singularity'' in the point $u = 0$ (when $n<1$). We start with the case $\gamma > 0$ and then we discuss the case $\gamma = 0$.
\begin{thm}\label{THEOREMSLOWDIFFUSIONSTRONGREACTION}
Let $m > 0$ and $p > 1$ such that $\gamma > 0$, $n \in \RR$ and $q := [\gamma + (n-1)(p-1)]/(p-1)$. Then there exist admissible TWs for equation \eqref{eq:SLOWDIFFSTRONGREACTPLAPLACIAN} if and only if $q \geq 0$.

\noindent Moreover, for all $q \geq 0$, there exists a critical speed $c_{\ast} = c_{\ast}(m,p,n) > 0$ such that equation \eqref{eq:SLOWDIFFSTRONGREACTPLAPLACIAN} possesses a unique admissible TW for all $c \geq c_{\ast}(m,p,n)$ and does not have admissible TWs for $0 < c < c_{\ast}(m,p,n)$. The TW corresponding to the speed $c_{\ast}(m,p,n)$ is finite.

\noindent Finally, we have:

\noindent $\bullet$ If $q = 0$, each TW is finite;

\noindent $\bullet$ If $q > 0$, the TWs corresponding to the values $c > c_{\ast}(m,p,n)$ are finite if and only if $0 < n < 1$.
\end{thm}
\emph{Proof.} The proof is very similar to the one of Theorem \ref{THEOREMEXISTENCEOFTWS} and we sketch it quickly. We start writing the equation of the profile getting to the system
\[
\frac{dX}{d\tau} = (p-1)X|Z|^{p-2} Z, \quad\quad \frac{dZ}{d\tau} = cZ - |Z|^p - mX^q(1-X)
\]
where $X$ and $Z$ are defined as in \eqref{eq:NONSTANDARDCHANGEOFVARIABLESTWS} and the equation of the trajectories
\[
\frac{dZ}{dX} = \frac{cZ - |Z|^p - mX^q(1-X)}{(p-1)X|Z|^{p-2} Z}.
\]
Now, using the same methods of the proof of Theorem \ref{THEOREMEXISTENCEOFTWS} (\emph{Step4}), it is not difficult to see that if $q < 0$, there are not trajectories linking the saddle point $S = (1,0)$ with a point of the type $(0,\lambda)$ for all $\lambda \geq 0$ and so, there are not admissible TWs.

If $q = 0$, i.e. $\gamma = (1-n)(p-1)$ (note that this expression makes sense only if $n < 1$), we have null isoclines qualitatively equal to the ones found in Theorem \ref{THEOREMEXISTENCEOFTWS1} and so, we can show the existence of a critical speed $c_{\ast} = c_{\ast}(m,p,n)$ such that there are not TWs for $c < c_{\ast}$ and there exists exactly one TW for all $c \geq c_{\ast}$. Moreover, since we know that $n < 1$ and using the same methods of \emph{Step3} of Theorem \ref{THEOREMEXISTENCEOFTWS}, it is simple to see these TWs are finite.

If $q > 0$,  the analysis is very similar to the one done in Theorem \ref{THEOREMEXISTENCEOFTWS}. The unique significant difference appears when we study the local behaviour of the trajectories from the critical point $O = (0,0)$ in the case $c > c_{\ast}$. Indeed, it is simple to see that the trajectories from $O$ satisfies $Z(X) \sim (m/c)X^q$ for $X \sim 0$. Hence, integrating the differential relation $d\xi = X^{(\mu-1)/(p-1)}Z^{-1}dX$ we get that the time (measured respect with the variable $\xi$) in which the profile approaches the level $0$ depends on $n$:
\[
\xi_1-\xi_0 = -\int_{X_0}^{X_1} \frac{dX}{X^{(1-\mu)/(p-1)}Z(X)} \sim -c \int_{X_0}^{X_1} \frac{dX}{X^{(1-\mu)/(p-1) + q}} = -c\int_{X_0}^{X_1} \frac{dX}{X^n},
\]
where $0 < X_1 < 1$ is fixed. Then, letting $X_0 \to 0$, it follows that the TW is finite if and only if $0 < n < 1$. $\Box$
\paragraph{Remarks.} (i) We have showed that the existence of admissible TWs depends on the sign of the value $q := [\gamma + (n-1)(p-1)]/(p-1)$ which, in some sense, describes the interaction between the diffusion and the reaction terms. In particular, we have proved a really interesting fact: if $q \geq 0$ and $0 < n < 1$, then all the admissible TWs for equation \eqref{eq:SLOWDIFFSTRONGREACTPLAPLACIAN} are finite. This represents a very important difference with the case $\gamma > 0$ and $n = 1$. Indeed, when $\gamma > 0$ and $n = 1$ only the TW corresponding to the critical value $c_{\ast}$ is finite.

(ii) We point out that the same procedure can be followed when $\gamma = 0$ and $n \in \RR$. In this case, we have $q = q(n) = n - 1$ and system \eqref{eq:SYSTEMNONSINGULARTWSXI2} has the form
\[
\frac{dX}{d\tau} = (p-1)X|Z|^{p-2}Z, \quad\quad \frac{dZ}{d\tau} = cZ - |Z|^p - mX^{n-1}(1-X).
\]
Following the proof of Theorem \ref{THEOREMEXISTENCEOFTWS} and Theorem \ref{THEOREMEXISTENCEOFTWS1}, it is not difficult to show that for all $n \geq 1$, there exists a critical speed of propagation $c_{\ast}(n) > 0$, such for all $c \geq c_{\ast}(n)$, there exists a unique positive TW for equation \eqref{eq:SLOWDIFFSTRONGREACTPLAPLACIAN}, while there are no admissible TWs for $c < c_{\ast}(n)$. Moreover, if $n < 1$  equation \eqref{eq:SLOWDIFFSTRONGREACTPLAPLACIAN} does not possess TW solutions.

\noindent This means that when $\gamma = 0$, a ``weak/strong'' modification of the reaction term is not sufficient to have finite TWs. So, the previous observations and Theorem \ref{THEOREMSLOWDIFFUSIONSTRONGREACTION} explain us the exact combination of slow diffusion ($\gamma > 0$) and strong reaction needed in order to ``generate'' \emph{only} finite TW solutions: we need to have $q \geq 0$ and $0 < n < 1$. Consequently, we can conclude that the method of balancing the doubly nonlinear slow diffusion with a strong reaction allows us to separate quantitatively the cases in which all the TWs are positive, there exists at least one finite TW and all the TWs are finite.

(iii) Note that Theorem \ref{THEOREMSLOWDIFFUSIONSTRONGREACTION} extends Theorem \ref{THEOREMSLOWDIFFUSIONSTRONGREACTIONPOROUSMEDIUM} proved in \cite{DePablo-Vazquez1:art} for the Porous Medium case (i.e., $p = 2$ and $m > 1$).

(iv) We point out that the study of the super-level sets and the asymptotic study can be repeated when $f(u) = u^n(1-u)$ and $0 < n < 1$. Indeed, since $u^n(1-u) \geq u(1-u)$, we can employ the solutions of the problem with $n = 1$ as sub-solutions for the problem with $0 < n < 1$. The significative open problem in the study of problem \eqref{eq:REACTIONDIFFUSIONEQUATIONPLAPLACIAN} with a strong reaction term is that the solutions are not unique since the reaction term is not Lipschitz continuous (see also \cite{DP-S:art}). When $n \geq 1$, the problem of the uniqueness of the solution does not appear, but it seems not simple to understand the asymptotic behaviour of the solution since the reaction term is not concave and we cannot adapt the methods used in Proposition \ref{LEMMANONDISCRETEVERSIONLEVELSETS} and Proposition \ref{LEMMANONDISCRETEVERSIONLEVELSETSPSEUDO}.
%
%
%
%
%
%
%
%
%
%
%
\section{Appendix: proof of formula (\ref{eq:ASYMPTOTICTALESTWPSEUDOCAST})}\label{Appendix1}
As we explained in Section \ref{CLASSIFICATIONEXISTENCETW}, studying the asymptotic behaviour of the positive TW with critical speed $c_{0\ast}(m,p) := p(m^2f'(0))^{\frac{1}{mp}}$, when $\gamma = 0$, is more complicated than the case $\gamma > 0$. In what follows, we give the detailed proof of formula \eqref{eq:ASYMPTOTICTALESTWPSEUDOCAST}.

\noindent Before proceeding with the analysis, let's recall some facts we have proved before (see the proof of Theorem \ref{THEOREMEXISTENCEOFTWS1}). We have worked in the $(X,Z)$-plane, where $X = \varphi$ and $Z = m X^{-1}X'$ are defined in \eqref{eq:NONSTANDARDCHANGEOFVARIABLE2} and we showed the existence of a trajectory linking $R_{\lambda_{\ast}} \leftrightarrow S$, where $R_{\lambda_{\ast}} = (\lambda_{\ast},0)$, $S = (1,0)$, and $\lambda_{\ast} := (c_{0\ast}(m,p)/p)^m = (m^2f'(0))^{1/p}$. This trajectory corresponds to an admissible positive TW with critical speed $c_{0\ast}(m,p)$.

\noindent (i) Now we can begin with the proof of \eqref{eq:ASYMPTOTICTALESTWPSEUDOCAST}. Let's suppose $f \in C^2([0,1])$ and let $Z = Z(X)$ be the analytic expression of the trajectory of this TW. We consider the approximation $Z \sim \lambda_{\ast} - \zeta_0(X)$ as $X \sim 0$, where $\zeta_0(X) \sim 0$ as $X \sim 0$. Our goal is to compute the remainder $\zeta_0 = \zeta_0(X)$ as $X \sim 0$. Let's consider the equation of the trajectories \eqref{eq:EQUATIONOFTHETRAJECTORIESPSEUDO}. Since
\begin{equation}\label{eq:SECONDORDERTAYLOREXPANSIONFAST}
\begin{aligned}
X^{-1}f(X) &\sim f'(0) + \frac{f''(0)}{2}X, \quad \text{for } X \sim 0,  \\
(\lambda_{\ast} - \zeta_0(X))^p &\sim \lambda_{\ast}^p - p\lambda_{\ast}^{p-1}\zeta_0(X) + \frac{p(p-1)}{2}\lambda_{\ast}^{p-2}\zeta_0(X)^2 \quad \text{for }\; X \sim 0,
\end{aligned}
\end{equation}
we have that equation \eqref{eq:EQUATIONOFTHETRAJECTORIESPSEUDO} with $c = c_{0\ast}$ becomes
\[
\frac{d\zeta_0}{dX} \sim -\frac{c_{0\ast}\lambda_{\ast} - \lambda_{\ast}^p - mf'(0) - (c_{0\ast} - p\lambda_{\ast}^{p-1})\zeta_0 - \frac{p(p-1)}{2}\lambda_{\ast}^{p-2}\zeta_0^2 - m\frac{f''(0)}{2}X}{(p-1)\lambda_{\ast}^{p-1}X} \sim \frac{p}{2\lambda_{\ast}}\frac{\zeta_0^2}{X} + b_0, \quad \text{for } X \sim 0
\]
where $b_0 := mf''(0)/[2(p-1)\lambda_{\ast}^{p-1}] < 0$. We point out that the second approximation holds since both the quantities $c_{0\ast}\lambda_{\ast} - \lambda_{\ast}^p - mf'(0)$ and $c_{0\ast} - p\lambda_{\ast}^{p-1}$ are zero. Now, consider the equation
\[
\frac{d\widetilde{\zeta}_0}{dX} = \frac{p}{2\lambda_{\ast}}\frac{\widetilde{\zeta}_0^{\,2}}{X},
\]
which is obtained by taking $b_0 = 0$ in the previous equation. We have $\widetilde{\zeta}_0(X) \sim -(2\lambda_{\ast}/p)\ln^{-1}(X)$ for $X \sim 0$ while, for all small $\varepsilon > 0$, it is simple to see that the functions
\[
\underline{\zeta}_0(X) \sim -\frac{2\lambda_{\ast}+\varepsilon}{p}\ln^{-1}(X) \quad \text{ and } \quad \overline{\zeta}_0(X) \sim -\frac{2\lambda_{\ast}-\varepsilon}{p}\ln^{-1}(X)
\]
are sub-solution and super-solution of the equation $d\zeta_0/dX = \frac{p}{2\lambda_{\ast}}\frac{\zeta_0^2}{X} + b_0$ for $X \sim 0$, respectively. Thus, using the arbitrariness of $\varepsilon > 0$ and Maximum Principle, we deduce that
\[
\zeta_0(X) \sim \widetilde{\zeta}_0(X) \sim -\frac{2\lambda_{\ast}}{p}\ln^{-1}(X), \quad \text{for }\; X \sim 0.
\]
(ii) Now, we iterate the previous procedure to compute high order terms. So let $Z \sim \lambda_{\ast} - \zeta_1(X)$, for $X \sim 0$, where $\zeta_1(X) = \zeta_0(X) + h.o.t$. We suppose for a moment to have $f \in C^3([0,1])$ and, proceeding as in \eqref{eq:SECONDORDERTAYLOREXPANSIONFAST}, we compute the Taylor expansions:
\[
\begin{aligned}
X^{-1}f(X) &\sim f'(0) + \frac{f''(0)}{2}X + \frac{f'''(0)}{6}X^2, \quad \text{for } X \sim 0,  \\
(\lambda_{\ast} - \zeta_1(X))^p &\sim \lambda_{\ast}^p - p\lambda_{\ast}^{p-1}\zeta_1(X) + \frac{p(p-1)}{2}\lambda_{\ast}^{p-2}\zeta_1(X)^2 - \frac{p(p-1)(p-2)}{6}\lambda_{\ast}^{p-3}\zeta_1(X)^3\quad \text{for }\; X \sim 0,
\end{aligned}
\]
and we substitute in the equation of the trajectories \eqref{eq:EQUATIONOFTHETRAJECTORIESPSEUDO} to find
\[
\frac{d\zeta_1}{dX} \sim \frac{p}{2\lambda_{\ast}}\frac{\zeta_1^2}{X} -\frac{p(p-2)}{6\lambda_{\ast}^2}\frac{\zeta_1^3}{X} + b_0 - b_1X, \quad \text{for } X \sim 0,
\]
where $b_1 :=mf'''(0)/[6(p-1)\lambda_{\ast}^{p-1}]$. Now, we know the less accurate approximation $\zeta_1(X) \sim \zeta_0(X) \sim -\frac{2\lambda_{\ast}}{p}(\ln X)^{-1}$, for $X \sim 0$. So, substituting it in the previous equation we obtain
\[
\frac{d\zeta_1}{dX} \sim \frac{2\lambda_{\ast}}{p}\frac{1}{X\ln^2(X)} + \frac{4(p-2)\lambda_{\ast}}{3p^2}\frac{1}{X\ln^3(X)} + b_0 - b_1X, \quad \text{for } X \sim 0,
\]
and once we integrate with respect to the variable $X$, we get
\[
\begin{aligned}
\zeta_1(X) &\sim -\frac{2\lambda_{\ast}}{p}\frac{1}{\ln(X)} - \frac{2(p-2)\lambda_{\ast}}{3p^2}\frac{1}{\ln^2(X)} + b_0X - \frac{b_1}{2}X^2 \\
& = \zeta_0(X) - \frac{2(p-2)\lambda_{\ast}}{3p^2}\frac{1}{\ln^2(X)} + b_0X - \frac{b_1}{2}X^2, \quad \text{for } X \sim 0.
\end{aligned}
\]
Using this information, the differential equation $X' = (1/m)XZ$ in \eqref{eq:SYSTEMNONSINGULARTWSXI1} becomes
\begin{equation}\label{eq:ASYMBEHGG0CRITSPEED}
X' \sim \frac{\lambda_{\ast}}{m}X\Big[1 + \frac{2}{p}\frac{1}{\ln(X)} + \frac{2(p-2)}{3p^2}\frac{1}{\ln^2(X)} - \frac{b_0}{\lambda_{\ast}}X + \frac{b_1}{2\lambda_{\ast}}X^2 \Big], \quad \text{for } X \sim 0.
\end{equation}
As the reader can easily see, a first approximation of the solution of \eqref{eq:ASYMBEHGG0CRITSPEED} is given by
\[
\ln X(\xi) \sim \frac{\lambda_{\ast}}{m}\xi, \quad \text{for } \xi \sim -\infty.
\]
Consequently, by substituting the previous expression in the square parenthesis of \eqref{eq:ASYMBEHGG0CRITSPEED} we have
\[
\frac{X'}{X} \sim \frac{\lambda_{\ast}}{m} + \frac{2}{p}\frac{1}{\xi} + \frac{2m(p-2)}{3\lambda_{\ast}p^2}\frac{1}{\xi^2} + \frac{b_0}{m}e^{\frac{\lambda_{\ast}}{m}\xi} - \frac{b_1}{2m}e^{\frac{2\lambda_{\ast}}{m}\xi}, \quad \text{for } \xi \sim -\infty,
\]
which, once integrated, can be re-written as
\[
\ln(X) \sim \frac{\lambda_{\ast}}{m}\xi + \frac{2}{p}\ln(|\xi|) - \frac{2m(p-2)}{3\lambda_{\ast}p^2}\frac{1}{\xi} + \frac{b_0}{\lambda_{\ast}}e^{\frac{\lambda_{\ast}}{m}\xi} - \frac{b_1}{4\lambda_{\ast}}e^{\frac{2\lambda_{\ast}}{m}\xi}, \quad \text{for } \xi \sim -\infty.
\]
Hence, we have shown that equation \eqref{eq:ASYMBEHGG0CRITSPEED} is satisfied by taking
\[
X(\xi) = \varphi(\xi) \sim a_0 |\xi|^{\frac{2}{p}}e^{\frac{\lambda_{\ast}}{m}\xi} = a_0 |\xi|^{\frac{2}{p}}\exp\Big(m^{\frac{2-p}{p}}f'(0)^{\frac{1}{p}}\xi\Big) \quad \text{for }\; \xi \sim -\infty,
\]
for some constant $a_0 > 0$, which is exactly \eqref{eq:ASYMPTOTICTALESTWPSEUDOCAST}. With the previous formula we complete the study of the asymptotic behaviour of the TW with critical speed for $\xi \sim -\infty$. We end this paragraph pointing out that neither the assumption $f \in C^2([0,1])$ nor $f \in C^3([0,1])$ are not needed since the terms involving $f''(0)$ and $f'''(0)$ do not influence formula \eqref{eq:ASYMPTOTICTALESTWPSEUDOCAST}.

\paragraph{Remark.} As verification, this estimate is consistent with the results known in the linear case. Indeed, when $m=1$ and $p=2$, the TW $\varphi = \varphi(\xi)$ with critical speed $c_{0\ast}(m=1,p=2) := c_{0\ast}' = 2\sqrt{f'(0)}$ satisfies
\[
\varphi(\xi) \sim |\xi|e^{\frac{c_{0\ast}'}{2}\xi} = a_0|\xi|e^{\sqrt{f'(0)}\xi}, \quad \text{for }\; \xi \sim -\infty,
\]
for some $a_0 > 0$. See for instance \cite{Ham-Nad:art} and the references therein.
%
%
%
%
%
%
%
%
%
%
%
\section{Comments and open problems}\label{SECTIONCOMMENTANDOPENPROBLEMS}
We end the paper by some comments, open problems and information related to our work.
\paragraph{On the exact location of the propagation front.} As mentioned in Section \ref{SECTIONRESULTSFORLINEARDIFFUSION}, determining the properties of the solution $u = u(x,t)$ on the moving coordinate $x = \xi - c_{\ast}t$ with  $c_{\ast}$ the critical speed is a very challenging problem.  Important steps forward have been made. For example, Bramson showed in \cite{B1:art} and \cite{B2:art}, using probabilistic techniques, an interesting property of the level sets $ E_{\omega}(t) = \{x > 0: u(x,t) = \omega \}$, $\omega \in (0,1)$, of the solution $u = u(x,t)$ of equation \eqref{eq:ONEDIMENSIONALKPPEQUATIONINTRO} with reaction term satisfying \eqref{eq:ASSUMPTIONSONTHEREACTIONTERM}. In particular, he proved that for all $\omega \in (0,1)$ there exist constants $x_{\omega}$, $a > 0$ and $C_{\omega} > 0$ such that
\begin{equation}\label{eq:BRAMSONCORRECTION1INTRO}
E_{\omega}(t) \subset \Big[c_{\ast}t - \frac{3}{2\omega_{\ast}}\ln t - x_{\omega} - \frac{a}{\sqrt{t}}- \frac{C_{\omega}}{t}, c_{\ast}t - \frac{3}{2\omega_{\ast}}\ln t - x_{\omega} - \frac{a}{\sqrt{t}} + \frac{C_{\omega}}{t}\Big]
\end{equation}
for $t$ large enough, where $\omega_{\ast} = c_{\ast}/2$. The previous formula is interesting since it allows to estimate the ``delay'' of the solution $u = u(x,t)$ from the positive TW with critical speed $c = c_{\ast}$ which, according to \eqref{eq:BRAMSONCORRECTION1INTRO}, grows in time and consists in a logarithmic deviance. More recently, a similar result have been proved in \cite{Hamel-N-R-R:art} with PDEs techniques. In particular, the authors showed that there exists a constant $C \geq 0$ such that
\[
E_{\omega}(t) \subset \Big[c_{\ast}t - \frac{3}{2\omega_{\ast}}\ln t - C, c_{\ast}t - \frac{3}{2\omega_{\ast}}\ln t + C] \quad \text{for } t \text{ large enough},
\]
which is less precise than \eqref{eq:BRAMSONCORRECTION1INTRO} but, anyway, it contains the most important information on the ``delay'' of the solution, i.e., the logarithmic shift (for more work on this issue see \cite{LA:art, R:art, Uchiyama:art}). Moreover, they gave a interesting proof of the uniform convergence of the general solutions of equation \eqref{eq:ONEDIMENSIONALKPPEQUATIONINTRO} to the TW solution with critical speed $c_{\ast}$ (see Theorem 1.2 of \cite{Hamel-N-R-R:art}).

\noindent It seems quite natural to conjecture that, at least in the ``pseudo-linear'' case, the level sets of the solution of \eqref{eq:REACTIONDIFFUSIONEQUATIONPLAPLACIAN} satisfy similar properties. Nevertheless, studying the problem of the exact location of the propagation front for the doubly nonlinear diffusion seems a really difficult task and we pose it as an interesting open problem.
\paragraph{Exponential propagation for fast diffusion.} In \cite{King-M:art}, King and McCabe studied the radial solution of the equation
\[
\partial_tu = \Delta u^m + u(1-u) \quad \text{in } \RR^N\times(0,\infty)
\]
for the fast diffusion case, $0 < m < 1$. Taking $L^1$ initial data, they showed that if $(N-2)_+/N := m_c < m < 1$, the solution $u = u(x,t)$ of the previous equation converges pointwise to 1 for large times with an exponential rate. This is a severe departure from the TW behavior of the standard Fisher-KPP model since there are no TW solutions and the asymptotic spread rate is faster than linear.

\noindent Evidently, it is possible to make a parallel study when the diffusion term is doubly nonlinear. In this case (see for instance \cite{V1:book}), the assumption on the parameters $m>0$ and $p>1$ can be written as
\[
-p/N < \gamma < 0
\]
and we refer to it as ``fast diffusion'' assumption. Thus, taking into account the results showed in \cite{King-M:art} for the Porous Medium equation, it seems plausible to conjecture that when the ``fast diffusion'' is satisfied, the general solutions of problem \eqref{eq:REACTIONDIFFUSIONEQUATIONPLAPLACIAN} expands exponentially in space for large times and so, the paradigm of the TW behaviour is broken. In particular, we make the following conjecture which stands for a ``fast diffusion'' version of Theorem \eqref{NTHEOREMCONVERGENCEINNEROUTERSETS}.
\begin{conj}
Let $m > 0$ and $p > 1$ such that $-p/N < \gamma < 0$, $N \geq 1$ and set $\sigma_{\ast} := -\frac{\gamma}{p}f'(0) > 0$.

\noindent (i) For all $\sigma < \sigma_{\ast}$, the solution $u = u(x,t)$ of problem \eqref{eq:REACTIONDIFFUSIONEQUATIONPLAPLACIAN} with initial datum \eqref{eq:ASSUMPTIONSONTHEINITIALDATUM} satisfies
\[
u(x,t) \to 1 \quad \text{ uniformly in } \{|x| \leq e^{\sigma t} \} \; \text{ as } t \to \infty.
\]
\noindent (ii) Moreover, for all $\sigma > \sigma_{\ast}$, it satisfies
\[
u(x,t) \to 0 \quad \text{ uniformly in } \{|x| \geq e^{\sigma t} \} \; \text{ as } t \to \infty.
\]
\end{conj}
Exponential propagation happens also with fractional diffusion, both linear and nonlinear, see \cite{C-R2:art} and \cite{S-V:art} and their references.  We will not enter here into the study of the relations of our paper with nonlinear fractional diffusion, though it is an interesting topic.

\noindent Finally, we recall that infinite of speed propagation depends not only on the diffusion operator but also on the initial datum. In particular, in \cite{Hamel-R:art}, Hamel and Roques found that the solutions of the Fisher-KPP problem with linear diffusion i.e., ($m=1$ and $p=2$) propagate exponentially fast for large times if the initial datum has a power-like spatial decay at infinity.
\paragraph{An interesting limit case.} In Section \ref{CLASSIFICATIONEXISTENCETW}, we studied the ``pseudo-linear'' case ($\gamma = 0$, i.e., $m(p-1) = 1$) finding an explicit formula for the function $c_{\ast}$ and we wrote it as a function of $m>0$ and/or $p > 1$. In particular, we found
\[
\lim_{m\to\infty}c_{\ast}(m) = \lim_{p\to1}c_{\ast}(p) = 1.
\]
This fact allows us to conjecture the existence of admissible TWs for the limit case $m\to\infty$ and $p\to1$. Now, we can formally compute the limit of the doubly nonlinear operator:
\[
\begin{aligned}
\Delta_pu^m &= m^{p-1}\nabla\cdot(u^{\mu}|\nabla u|^{p-2}\nabla u) \\
&= m^{1/m} \nabla\cdot(u^{2-p}|\nabla u|^{p-2}\nabla u) \to \nabla\cdot\Bigg(u\frac{\nabla u}{|\nabla u|}\Bigg)  \quad \text{as } m \to \infty \text{ and } p \to 1,
\end{aligned}
\]
keeping $m(p-1) = 1$. We ask the reader to note that, assuming $m(p-1) = \theta$ with $\theta > 0$, we can repeat the previous computations and deduce that
\[
\Delta_pu^m \to \nabla\cdot\Bigg(u^{\theta}\frac{\nabla u}{|\nabla u|}\Bigg)  \quad \text{as } m \to \infty \text{ and } p \to 1,
\]
keeping $m(p-1) = \theta$ fixed. Consequently, a very interesting open problem is the study of the existence of admissible TWs for the equation
\begin{equation}\label{eq:SIMILARTOFLUXLIMITED}
\partial_t u = \partial_x\big(u^{\theta}|\partial_x u|^{-1}\partial_x u\big) + f(u) \quad \text{in } \RR\times(0,\infty)
\end{equation}
for different values of the parameter $\theta \geq 0$. The case $\theta = 1$ has been studied by Andreu et al. in \cite{A-Cas-M:art}, where the authors showed the existence of discontinuous TWs for equation \eqref{eq:SIMILARTOFLUXLIMITED} (with $\theta = 1$). Hence, it seems reasonable to conjecture that in the limit case $m\to \infty$ and $p \to 1$ with $m(p-1) = \theta > 1$ there are admissible TWs but with less regularity.

\noindent Finally, we recommend the papers \cite{C-C-Cas-S-S:art, C-G-S-Sol:art, Cam-Sol:art} for more work on TW solutions (in general discontinuous) to a nonlinear flux limited Fisher-KPP equation, which seem to be related with the limit TWs of our work.

\vskip 1cm

\noindent {\textbf{\large \sc Acknowledgments.} Both authors have been partially funded by Projects MTM2011-24696 and  MTM2014-52240-P (Spain). Work partially supported by the ERC Advanced Grant 2013 n.~339958 ``Complex Patterns for Strongly Interacting Dynamical Systems - COMPAT''.

%
%
%
\vskip 1cm

\

2000 \textit{Mathematics Subject Classification.}
35K57,  
35K65, 
35C07,  	
35K55, 

\medskip

\textit{Keywords and phrases.}  Fisher-KPP equation, Doubly nonlinear  diffusion, Propagation of level sets.

\

\end{document}